\documentclass[11pt,leqno]{article}
\usepackage{amsmath,amssymb}

\renewcommand{\thesection}{\arabic{section}}

\newcommand{\email}[1]{{\it  E-mail address:}\  \textsf{#1}\\*\protect}
\newcommand{\address}[1]{\small \textsc{#1}\\*\protect}

\title{The dimension of the boundary of super-Brownian motion}
\bigskip
\author{Leonid Mytnik $\mbox{}^{1}$ \hspace{1cm} 
 Edwin Perkins $\mbox{}^{2}$  }
\date{}
\begin{document}
\maketitle
\begin{center}
\address{Faculty of Industrial Engineering and Management,}
\address{Technion -- Israel Institute of Technology, Haifa 32000, Israel}
\email{leonid@ie.technion.ac.il}
\mbox{}\\
\address{ Department of Mathematics, The University of British Columbia,}
\address{ 1984 Mathematics Road, Vancouver, B.C., Canada V6T 1Z2}
\email{perkins@math.ubc.ca}
\end{center}

\vskip 0.3in 
{\bf Abstract.} We show that  the Hausdorff dimension of the boundary of $d$-dimensional super-Brownian motion is $0$, if $d=1$, $4-2\sqrt2$, if $d=2$, and $(9-\sqrt{17})/2$, if $d=3$.

\vfill
{\footnoterule
\footnotesize
\noindent 
\today\\
AMS 2000 {\it subject classifications}. Primary 60H15, 60G57. 
Secondary 28A78, 35J65, 60J55, 60H40, 60J80. 
  \\
{\it Running head}. The boundary of Super-Brownian motion\\
1. Supported in part by the Israel Science Foundation 
(grant No. 1325/14). \\ 
2. Supported by an NSERC Discovery grant.\\}
\pagebreak
\vspace{1cm}

\newcommand{\mysection}[1]{\section{#1}
\markright{{\protect\footnotesize Chapter \thechapter, Section \thesection}}}
\newcommand{\resc}[2]{{#1}^{n}(\qn)={#2}
\lbr{\displaystyle \frac{\qn}{n}} \rbr}
\newcommand{\rescb}[2]{\bar{#1}^{n}(\qn)={#2}
\lbr{\displaystyle \frac{\qn}{n}} \rbr}
\newcounter{archapter}[section]
\newtheorem{theorem}{Theorem}[section]
\newtheorem{lemma}[theorem]{Lemma}
\newtheorem{definition}[theorem]{Definition}
\newtheorem{remark}[theorem]{Remark}
\newtheorem{proposition}[theorem]{Proposition}
\newtheorem{case}{Case}[archapter]
\newtheorem{corollary}[theorem]{Corollary}
\newtheorem{fact}[theorem]{Fact}
\newtheorem{assumption}[theorem]{Assumption}
\newtheorem{convention}[theorem]{Convention}
\newtheorem{conjecture}[theorem]{Conjecture}
\newcommand{\gdm}{\hfill\vrule  height5pt width5pt \vspace{.1in}}
\newcommand{\itemstr}[1]{\item[{\bf (#1)}]}
\newfont{\msb}{msbm10 scaled \magstep1}
\newfont{\msbh}{msbm7 scaled \magstep1}
\newfont{\msbhh}{msbm5 scaled \magstep1}
\renewcommand{\Re}{\mbox{\msb R}}
\newcommand{\IR}{\mbox{\msb R}}
\newcommand{\IRsm}{\mbox{\msbh R}}
\newcommand{\Z}{\mbox{\msb Z}}
\newcommand{\Ze}{\mbox{\msb Z}}
\newcommand{\N}{\mbox{\msb N}}
\newcommand{\Zesm}{\mbox{\msbh Z}}
\newcommand{\Ce}{\mbox{\msb C}}
\newcommand{\Cesm}{\mbox{\msbh C}}
\newcommand{\Te}{\mbox{\msb T}}
\newcommand{\Tesm}{\mbox{\msbh T}}
\newcommand{\Le}{\mbox{\msb L}}
\newcommand{\Lesm}{\mbox{\msbh L}}
\newcommand{\Ee}{\mbox{\msb E}}
\newcommand{\Pe}{\mbox{\msb P}}
\newcommand{\Resm}{\mbox{\msbh R}} 
\newcommand{\Resmm}{\mbox{\msbhh R}} 
\newcommand{\Rp}{\Re_{+}}
\newcommand{\Rpf}{\Re_{+}}
\newcommand{\Rpsm}{\Resm_{+}}
\newcommand{\Rd}{\Re^{d}}
\newcommand{\Rdsm}{\Resm^{d}} 
\newcommand{\Rdsmm}{\Resmm^{d}} 
\newcommand{\Rdn}[1]{\Re^{#1}}
\newcommand{\Rdnsm}[1]{\Resm^{#1}}
\newcommand{\Zcal}{{\cal Z}}
\newcommand{\Ecal}{{\cal E}}
\newcommand{\esssup}{\mathrm{ess\, sup}} 
\newcommand{\lsc}{\left\langle}
\newcommand{\rsc}{\right\rangle}
\newcommand{\lbr}{\left(}
\newcommand{\rbr}{\right)}
\newcommand{\lfi}{\left\{}
\newcommand{\rfi}{\right\}}
\newcommand{\lsbr}{\left[}
\newcommand{\rsbr}{\right]}
\newcommand{\rmed}{\right|}
\newcommand{\Ged}{\left|}
\newcommand{\lnor}{\left\|}
\newcommand{\rnor}{\right\|}
\newcommand{\lnorm}{\lmed \lmed}
\newcommand{\rnorm}{\rmed \rmed}
\newcommand{\cndp}{\stackrel{\rm P}{\rightarrow}}
\newcommand{\mf}{M_{F}(\Rd)}
\newcommand{\mfo}{{{\cal M}_{F}}}
\newcommand{\scf}[2]{\lsc #1,\,#2 \rsc} 
\newcommand{\Nat}{{\cal N}} 
\newcommand{\rmp}{{\mathrm{p}}}
\newcommand{\ep}{\epsilon}
\newcommand{\ve}{\varepsilon}
\newcommand{\la}{\langle}
\newcommand{\ra}{\rangle}
\newcommand{\cF}{{\cal F}}
\newcommand{\cE}{{\cal E}}
\newcommand{\cB}{{\cal B}}
\newcommand{\cL}{{\cal L}}
\newcommand{\cR}{{\cal R}}
\newcommand{\cI}{{\cal I}}
\newcommand{\cW}{{\cal W}}
\newcommand{\tcI}{\tilde{\cal I}}
\newcommand{\nn}{\nonumber}
\newcommand{\om}{\overline m}
\newcommand{\oln}{\overline{l_n}}
\newcommand{\mcF}{\mathcal{F}}
\newcommand{\Ft}[1]{\mbox{${\cal F}_{#1}$}}
\newcommand{\mn}{\medskip\noindent}
\newcommand{\vsp}{\vspace{4.5 mm}} 
\newcommand{\mvsp}[1]{\vspace{#1 mm}} 
\newcommand{\hsp}{\mbox{\hspace{10 mm}}} 
\newcommand{\veps}{\varepsilon}
\newcommand{\eop}{\gdm}
\newcommand{\noi}{\noindent}
\newcommand{\mhsp}[1]{\mbox{\hspace{#1 mm}}}
\newcommand{\lb}{\mathsf{L}}
\newcommand{\rb}{\mathsf{R}}
\newcommand{\y}{y_0}
\newcommand{\z}{\eta}
\newcommand{\xii}{\xi}
\newcommand{\intobl}[1]{
 \int\!\!\int_{{}_{{}_{\!\!\!\!\!\!\!\!\!\mbox{\scriptsize $#1$}}}}\;}
\section{Introduction}\label{intro}
\setcounter{equation}{0}
\setcounter{theorem}{0}

\medskip

We consider a $d$-dimensional super-Brownian motion $(X_t,t\ge0)$ starting at $X_0$ under $\Pe_{X_0}$ with $d\le 3$.  Here $X_0\in M_F(\IR^d)$, the space of finite measures on $\IR^d$ with the weak topology, $X$ is continuous $M_F(\IR^d)$-valued Markov process, and $\Pe_{X_0}$ denotes any probability under which $X$ is as above.  We write $X_t(\phi)$ for the integral of $\phi$ with respect to $X$, and take our branching rate to be one, so that for any non-negative bounded Borel functions $\phi,f$ on $\IR^d$, 
\begin{equation}\label{LF}
\Ee_{X_0}\Bigl(\exp\Bigl(-X_t(\phi)-\int_0^t X_s(f)ds\Bigr)\Bigr) =\exp(-X_0(V_t(\phi,f)),
\end{equation}
where $V_t(x)=V_t(\phi,f)(x)$ is the unique solution of the mild form of 
\begin{equation}\label{SLLE}
\frac{\partial V}{\partial t}=\frac{\Delta V_t}{2}-\frac{V_t^2}{2}+f,\quad V_0=\phi,
\end{equation}
that is,
\begin{equation}
V_t=P_t(\phi)+\int_0^tP_s\Bigl(f-\frac{V_{t-s}^2}{2}\Bigr)\,ds.
\end{equation}
In the above $P_t$ is the semigroup of standard $d$-dimensional Brownian motion. See Chapter II 
of \cite{bib:per02} for the above and further properties.  
For $d$ as above, $X$ has a jointly lower semi-continuous local time $L^x_t$ which is monotone increasing in $t$ for all $x$, and satisfies
$$\int_0^t X_s(f)\,ds=\int_{\IR^d} f(x)L^x_t\,dx\text{ for all $t\ge 0$ and non-negative measurable }f.$$
Moreover, $L$ is jointly continuous on $\{(t,x): (t,x) \text{ is a continuity point of }X_0q_t(x)\}$, where $q_t(x)=\int_0^tp_s(x)ds$, $p_t$ is the Brownian density and $X_0q_t(x)=\int q_t(x-y)X_0(dy)$ (see Theorem 3 in \cite{bib:sug89} and Theorem~\ref{LTmodulus} in the next section.)   The fact that $X$ has an a.s. finite extinction time $\zeta=\inf\{t:X_t(1)=0\}$, means that $L^x=L^x_\infty=L^x_\zeta$ is also lower semicontinuous and is continuous on $\text{Supp}(X_0)^c$, where $\text{Supp}(X_0)$ is the closed support of $X_0$.  If $d=1$, then $L^x_t$ and $L^x$ are globally continuous. 

We will largely be concerned with the case $X_0=\delta_0$.  If $d=2$ or $3$, then $L^x_t$ is a.s. jointly continuous on $\IR_+\times\{x\in\IR^d:x\neq 0\}$, $L^x$ is continuous on $\{x\neq 0\}$
and 
\begin{equation}\label{LTexplodes}
L^x_t\to\infty\text{ as }x\to 0\text{ for any }t>0\ P_{\delta_0}-\text{a.s.}
\end{equation}
 (see \cite{bib:hong17} for a precise rate of explosion).  In this case we also define the boundary of $X$ to be 
\begin{equation}\label{frontdef}
F=\partial\{x:L^x>0\}.
\end{equation}

In this work we will find the Hausdorff dimension of $F$, $\text{dim}(F)$.
  Note that if a small $B$ intersects $F$, then $B$ will contain points $x$ where  $L^x$ is small but positive and so
a related question is to find $\alpha$ so that
\begin{equation}\label{lefttail} \Pe_{\delta_0}(0<L^x<a)\sim a^\alpha\text{ as }a\downarrow 0,
\end{equation}
where $\sim a^\alpha$ means bounded below and above by $ca^\alpha$ with different positive constants $c$. Not surprisingly, the resolution of \eqref{lefttail} will play an important role in finding $\text{dim}(F)$ (e.g. in the implicit derivation of Theorem~\ref{polarsetslb} below).  Note that $F$ is a delicate set as it depends on population behaviour when the population is 
quite sparse and so is prone to instabilities.  Moreover it is not monotone in the initial mass.  This complicates many of our arguments.  

Motivation for the study of $F$ arose, in part, from a natural interest in the interface between visited and unvisited sites in a population.  We know that even for low dimensions $d\le 3$, $X$ arises as a scaling limit for the long-range contact process(\cite{bib:MT95},\cite{bib:DP99}), the voter model  and Lotka-Volterra model (\cite{bib:cdp00}, \cite{bib:cp05}), and long range percolation (\cite{bib:lz10}).  So, modulo the obvious problems with interchanging limits, $F$ may describe the large scale behaviour of the interface between infected and non-infected sites in an epidemic, or the boundary between two competing and coexisting species.  Another point of entry was the analysis in \cite{bib:mmp16} of $BZ_t=\partial\{x:X(t,x)>0\}$
where $X(t,x)$ is the density of $X_t$ in one spatial dimension.  There $\dim(BZ_t)\in(0,1)$ was found in terms of the lead eigenvalue of a killed Ornstein-Uhlenbeck operator and the problem had ties
to pathwise uniqueness in SPDE's for $X$ and related processes.  It is natural to apply some of the methods used there to its elliptical counterpart considered here.  A number of new tools in our present setting, including the theory of exit measures of $X$ and a stochastic analysis result of Yor (see Proposition~\ref{yorthm} below), will in fact allow us to explicitly calculate the dimension in this setting.  A final motivation was a comment of Itai Benjamini that the boundary of the range of a scaling limit of tree-indexed random walks  seemed to exhibit interesting fractal behaviour.  

To state our main results we define
$$p=p(d)=\begin{cases}3 &\text{ if }d=1\\
2\sqrt 2 &\text{ if }d=2\\
\frac{1+\sqrt {17}}{2} &\text{ if }d=3,
\end{cases}
$$
and
\[\alpha=\alpha(d)=\frac{p(d)-2}{4-d}=\begin{cases} 1/3&\text{ if }d=1\\
 \sqrt 2-1&\text{ if }d=2\\
\frac{\sqrt{17}-3}{2}&\text{ if }d=3.
\end{cases}\]
\begin{theorem}\label{dimthm} With $\Pe_{\delta_0}$-probability one, 
$$\textnormal{dim}(F)= d+2-p=\begin{cases}0&\text{ if } d=1\\
4-2\sqrt2\approx1.17&\text{ if }d=2\\
\frac{9-\sqrt{17}}{2}\approx 2.44&\text{ if }d=3.
\end{cases}$$
\end{theorem}

We will also consider $L^x$ and $F$ under the canonical measures $\N_{x_0}$ for $X$ starting at
$x_0$, which governs the evolution of a single super-Brownian cluster 
evolving from a single ancestor at $x_0$ at time $0$ and so perhaps is a more natural
setting for our questions.  (We will be using the same notation for the excursion measure of the Brownian snake from $\{x_0\}$.)
 $\N_{x_0}$ is a $\sigma$-finite measure on the space of continuous $M_F(\IR^d)$-valued paths such that
 \begin{equation}\label{PPPcan} X_t=\int \nu_t \Xi(d\nu) \quad\text{under $\Pe_{X_0}$},
 \end{equation}
 where $\Xi$ is a Poisson point process with intensity $\N_{X_0}=\int \N_{x_0}(\cdot)X_0(dx_0)$ (see, e.g., Theorem II.7.3 of \cite{bib:per02}). The existence of $(L^x,x\in\Re^d)$ under $\N_{x_0}$ follows from the above Poisson decomposition and its existence under $\Pe_{\delta_{x_0}}$.  The facts that there are finitely many clusters contributing to $L^x$
 for $x\neq x_0$ and that $(L^x,x\neq x_0)$ has a continuous version under  $\Pe_{\delta_{x_0}}$ easily imply that $x\to L^x$ has a continuous version on $\Re^d\setminus\{x_0\}$ $\N_{x_0}$-a.e. The details are standard. In fact, it is not hard to see that $(t,x)\to L^x_t$ is globally continuous $\N_{x_0}$-a.e. (see \cite{bib:hong17}) but we will not use this result here. 

\begin{theorem}\label{dimN0} $\textnormal{dim}(F)= d+2-p$ $\N_0$-a.e.
\end{theorem}

 Turning to \eqref{lefttail}, we have:

\begin{theorem}\label{lefttailthm}
(a) There is a $C_{\ref{lefttailthm}}$ such that 
\[\Pe_{\delta_0}(0<L^x\le a)\le C_{\ref{lefttailthm}}|x|^{-p}a^{\alpha}\ \forall x\in\IR^d,\ a\ge 0.\]
(b) For any $\veps_0>0$ there is a $c_{\ref{lefttailthm}}(\veps_0)>0$ such that
\[\Pe_{\delta_0}(0<L^x\le a)\ge c_{\ref{lefttailthm}}(\veps_0)|x|^{-p}a^\alpha\ \forall |x|\ge\veps_0,\ a\in[0,1].\]
(c) The same result (a) holds if $\Pe_{\delta_0}$ is replaced with $\N_0$.\\
(d)  There is a $\underline c_{\ref{lefttailthm}}>0$ such that
\[\N_{0}(0<L^x\le a)\ge \underline c_{\ref{lefttailthm}}(|x|^{-2}\wedge |x|^{-p})a^\alpha\ \forall |x|>0,\ a\in[0,1].\]
\end{theorem}
\medskip
The key to the above result (via Laplace transforms) will be a rate of convergence result for solutions of some semilinear pde's which may be of independent interest (Proposition~\ref{Vlambdarate} below). 

Next consider the dimension question for a general finite initial condition $X_0$. We let $\text{conv}(X_0)$ be the closed convex hull of $\text{Supp}(X_0)$.

\begin{theorem}\label{dimX0}
With $\Pe_{X_0}$-probability one:\\
\noindent(a) $\textnormal{dim}(F\cap(\textnormal{Supp}(X_0)^c)\le d+2-p$.\\
\noindent(b)
\begin{equation}\label{dimFX_0}
(\textnormal{conv}(X_0))^c\cap F\neq\emptyset\Rightarrow \textnormal{dim}((\textnormal{Supp}(X_0))^c\cap F)= \textnormal{dim}((\textnormal{conv}(X_0))^c\cap F)=d+2-p.
\end{equation}
\end{theorem}

The presence of $\text{conv}(X_0)$ in the hypothesis of \eqref{dimFX_0} is surely an artifice of our method which uses exit measures on hyperplanes--the above result should hold with $\text{Supp}(X_0)^c$ in place of $\text{conv}(X_0)^c$.
The hypothesis in \eqref{dimFX_0} is needed as the following simple example shows (see  Section~\ref{secproofs}).

\begin{proposition}\label{noexit}
Assume $X_0$
has support $\overline B_1=\{x:|x|\le 1\}$ and\\
 $\int_{\overline B_1}(1-|x|)^{-2}X_0(dx)<\infty$. Then there exists $c_d>0$  such that 
\[ \Pe_{X_0}(F\subset{\textnormal{Supp}}(X_0)=\textnormal{conv}(X_0))\ge \exp\Bigl(
-c_d
\int_{\overline{B}_1} (1-|x|)^{-2}X_0(dx)\Bigr)>0.\]
\end{proposition}

The behaviour of $F\cap\text{Supp}(X_0)$ can be quite different than that of $F\cap\text{Supp}(X_0)^c$, as the following simple Proposition shows in $d=3$. The proof is given in Section~\ref{secproofs}.

\begin{proposition}\label{FLeb} Assume $d=3$ and $X_0$ has a bounded density $X_0(x), x\in \Re^3$ which is not identically $0$.  Then with positive $\Pe_{X_0}$ probability, $F\cap\{X_0>0\}$ has positive Lebesgue measure.  
\end{proposition}
 
 Basically the dynamics underlying the behaviour of $F\cap\text{Supp}(X_0)$ are those of instantaneous extinction at $t=0$, and are quite different from those determining $F$ under $\N_0$ or $\Pe_{\delta_0}$.
 
 If $I(A)=\int_0^\infty X_s(A)ds$, we will define the range $\cR$ of $X$ to be 
 \[\cR=\textnormal{Supp}(I)=\overline{\{x:L^x>0\}}.\]
 A slightly smaller set is usually used in the literature (see \cite{bib:dip89} or Corollary~9 in Ch. IV of \cite{bib:leg99}) but the two definitions coincide under $\Pe_{\delta_0}$ or $\N_0$ and more generally produce the same outcomes for $\cR\cap\text{Supp}(X_0)^c$ and $\partial \cR\cap \text{Supp}(X_0)^c$.
 The topological boundary, $\partial \cR$, of $\cR$ is clearly closely related to $F$ and it is easy to check that
 \begin{equation}\label{boundaryofR} \partial \cR\subset F.
 \end{equation}
 Note that any isolated zeros of $L$ will be in $F$ but not $\partial \cR$ but we do not know if such points exist (they don't $\Pe_{\delta_0}$-a.s. if $d=1$ by the next Theorem).  More generally $x\in F$ will be in $\partial \cR$ iff there are open sets $\{U_n\}$ converging to $\{x\}$ so that $L|_{U_n}$ is identically $0$.  
 In general we do not know if $F=\partial \cR$ $\Pe_{\delta_0}$-a.s. but this is the case for $d=1$ by the 
 following theorem which also refines Theorem~\ref{dimthm} for $d=1$, and is proved in Section~\ref{sec1d}.

 \begin{theorem}\label{onedims}
 If $d=1$ then $\Pe_{\delta_0}$-a.s. there are random variables $\lb<0<\rb$ such that \[\{x:L^x>0\}=(\lb,\rb)\]
 and 
 in particular $F=\partial \cR=\{\lb,\rb\}$.
 \end{theorem}
 
 \begin{remark}\label{oneddens}(a)  By making minor changes in the proof of the above result one can, for example, show that if $X_0$ has a 
 continuous density such that $\{x:X_0(x)>0\}=(\ell_0,r_0)$ is a finite interval, then $\{x:L^x>0\}=(\lb,\rb)$ for some finite random variables satisfying $\lb\le \ell_0<r_0\le \rb$. Moreover $\Pe_{X_0}(\rb=r_0)=\exp\Bigl(-\int 6(r_0-x_0)^{-2}dX_0(x_0)\Bigr)$, which may or may not be $0$. There is some simplification in the argument as now $L^x$ is globally continuous. 
 \medskip
 
\noindent(b) In the proof of Proposition~\ref{FLeb} we will show that Supp$(X_0)\subset\cR$ $\Pe_{X_0}$-a.s. (see \eqref{SinR} in Section~\ref{secproofs}) and the proof applies equally well to the slightly smaller definition of range mentioned above. Therefore if, in Proposition~\ref{FLeb}, we also assume the initial density $X_0(x)$ is continuous, then $\{X_0>0\}\subset \text{Int(Supp}(X_0))\subset\text{Int}(\cR)$ which implies that $\{X_0>0\}\cap\partial\cR=\emptyset$.  In view of Proposition~\ref{FLeb}, we see that 
$F$ and $\partial \cR$ can be quite different inside $\{X_0>0\}$ for $d=3$.

 \end{remark}
 
We do conjecture that for $d=2$ or $3$,
\begin{equation}\label{dimssame}
\textnormal{dim}(\partial \cR)=\textnormal{dim}(F)\quad\Pe_{\delta_0}-a.s.\text{ and }\N_0-a.e.,
\end{equation}
and hope to return to this in a future work. 

The connection between super-Brownian motion and the exponent $p(d)$ was first established by Abraham and Werner  in their very interesting work on intersection exponents \cite{bib:AW97}.  Among other things they
prove:
\begin{equation}\label{AW}  \text{For $d\le 3$, $\N_\veps\Bigl(\sup_{x\in\cR}\Vert x\Vert\ge 1, 0\notin\cR\Bigr)\sim\veps^{p(d)+2-d}$ as $\veps\downarrow 0$,}
\end{equation}
 where $\N_\veps$ denotes the canonical measure
starting at a point a distance $\veps$ from the origin. Although such estimates are not entirely unrelated
from the probability that a small ball overlaps with $F$ we were not able to make direct use of these bounds.  On the other hand, the source of their exponents was a reformulation of a Girsanov theorem for Bessel processes due to Marc Yor( \cite{bib:yor92}), described below in Proposition~\ref{yorthm}, and as Proposition~1 of \cite{bib:AW97}. As in \cite{bib:AW97}, this result will play a central role in the derivation of our probabiliity
bounds.  

We next give some heuristics on how $p(d)$ enters in Theorems~\ref{dimthm} and \ref{lefttailthm}.  From \eqref{LF} one readily shows (see Lemma~\ref{LTVL}) that for $\lambda\ge 0$,
\begin{equation}\label{LVlambda}
\Ee_{\delta_0}(e^{-\lambda L^x})=\exp(-V^\lambda(x)),
\end{equation}
where $V^\lambda$ is the unique solution to 
\begin{equation}\label{vlequation}
\frac{\Delta V^\lambda}{2}=\frac{(V^\lambda)^2}{2}-\lambda\delta_0,\ \ V^\lambda>0\text{ on }\IR^d,\end{equation}
Letting $\lambda\uparrow\infty$ in the above we see that $V^\lambda(x)\uparrow V^\infty(x)$ and 
\begin{equation}\label{VinfL}
\Pe_{\delta_0}(L^x=0)=\exp(-V^\infty(x)).
\end{equation}
The fact that we know $V^\infty(x)=\frac{2(4-d)}{|x|^2}$ (see \eqref{vinftylimit}), is one of the reasons our results are more explicit than those in \cite{bib:mmp16} where the role of $V^\infty$ is played by a solution, $F$, of a nonlinear second order ordinary differential equation (see (3.1) in \cite{bib:mmp16}). Then $p=p(d)$ is the unique positive power such that $f(x)=|x|^{-p}$ is harmonic for Brownian motion with $V^\infty$ killing, that is
\begin{equation}\label{killharm}
\frac{\Delta f}{2}-V^\infty f=0 \text{ on }\Re^d\setminus\{0\},
\end{equation}
as one can easily check.  To see how this leads to Theorem~\ref{lefttailthm}, a Tauberian theorem will reduce the bounds in Theorem~\ref{lefttailthm}(a,b) to 
\begin{equation}\label{LTboundsI}
\Ee_{\delta_0}(e^{-\lambda L^x}1(L^x>0))\sim|x|^{-p}\lambda^{-\alpha}\text{ as }\lambda\uparrow\infty.
\end{equation}
The left-hand side of the above behaves  like $d^\lambda(x):=V^\infty(x)-V^\lambda(x)$, and so we need a rate of convergence of $V^\lambda$ to $V^\infty$.  Clearly we have $\frac{\Delta d^\lambda}{2}=\frac{V^\lambda+V^\infty}{2}d^\lambda$, and so by Feynman-Kac, if $|B_0|=|x|>1$ and $T_1$ is the first time the Brownian motion $B_t$ enters the closed unit ball, 
\begin{align}\label{FKI}
d^\lambda(x)&=E_x\Bigl(d^\lambda(B_{T_1})\exp\Bigl(-\int_0^{T_1}\frac{V^\infty+V^\lambda}{2}(B_s)ds\Bigr)\Bigr)\\
\nonumber&\sim d^\lambda(1)E_x\Bigl(f(B_{T_1})\exp\Bigl(-\int_0^{T_1} V^\infty (B_s)ds\Bigr)\Bigr)\quad(\text{for }\lambda\text{ large}).
\end{align}
In the last line we used radial symmetry of $d^\lambda$, $f(x)=1$ for $|x|=1$, and $d^\lambda(\infty)=0$ (in case $T_1=\infty$ for $d=3$). We have also conveniently  replaced $V^\lambda$ with $V^\infty$ for $\lambda$ large, where of course it is the difference of these functions that we are trying to bound. So by \eqref{killharm} we see that \eqref{FKI} implies 
\[d^\lambda(x)\sim d^\lambda(1)f(x)\text{ say as }x\to\infty\text{ for }\lambda\ge\lambda_0.\]
If $r=r(\lambda)$ is such that $\lambda r^{4-d}=\lambda_0$, then a simple scaling argument shows that as $\lambda\to\infty$ (and so $r\downarrow 0$),
\[
d^\lambda(x)=r^{-2}d^{\lambda r^{4-d}}(x/r)\sim r^{-2}d^{\lambda_0}(1)|x/r|^{-p}=C|x|^{-p}\lambda^{-\alpha}.\]
To handle the heuristic asymptotic in \eqref{FKI} we will use Yor's Girsanov Theorem (Proposition~\ref{yorthm}).  The careful derivation of Theorem~\ref{lefttailthm}(a,b) is in Section~\ref{seclefttail} while parts (c,d) are proved in Section~\ref{secproofs}.

To see how Theorem~\ref{lefttailthm} might give the upper bound on  dim$(F)$ in Theorem~\ref{dimthm},
at least for $d=3$, recall first that $x\to L^x$ is locally H\"older of index $1/2-\eta$ on $\{x\neq 0\}$ (see Theorem~\ref{LTmodulus}).  One can improve this modulus to H\"older $1-\eta$ if one of the endpoints is in the zero set of $L$ (see Theorem~2.3 in \cite{bib:mp11} for a similar result for the density of $X$ if $d=1$). 
So, ignoring the reduction by $\eta$ for convenience, we see that for $|x|\ge \veps_0$ and $\veps>0$ small enough,
\begin{align*}\Pe_{\delta_0}(F\cap\{y:|y-x|<\veps\}\neq\emptyset)&\le \Pe_{\delta_0}(0<L^y<\veps\text{ for some }|y-x|<\veps)\\
&\sim P(0<L^x<\veps)\le C\veps^\alpha=C\veps^{p-2},\end{align*}
where the first inequality is by the improved modulus of continuity, the second asymptotic is wishful
thinking, the last inequality is by Theorem~\ref{lefttailthm}, and the final equality uses $d=3$. 
A standard covering argument would then show that dim$(F)\le d-(p-2)$ a.s.  Although it is in fact
possible to make this argument rigorous for $d=3$, the situation for $d=2$ seems more difficult.  Instead we will establish the upper bound on dim$(F)$ in Section~\ref{secupperbound} using Dynkin's exit measures $X_{G_\veps}$ from $G_\varepsilon$, the complement of a closed ball of radius $\veps$. See Proposition~\ref{exittail} for the analogue of Theorem~\ref{lefttailthm} for such exit measures and Section~\ref{prel} for information on exit measures in general. 

In Section~\ref{lbh} we show the lower bound on dim$(F)$ in Theorem~\ref{dimthm} holds {\it with positive probability} by an energy calculation.  The standard Frostman method would construct a random measure $\cL$ supported by $F$ so that 
\begin{equation}\label{EC}
E\Bigl(\int\int1(|x_1|\le K, |x_2|\le K)|x_1-x_2|^{-(d+2-p-\eta)}\cL(dx_1)\cL(dx_2)\Bigr)<\infty,\text{ for all }K,\eta>0
\end{equation}
and so conclude dim$(F)\ge d+2-p$ on $\{\cL\neq 0\}$.  Define 
$$\cL^\lambda(\phi)=\lambda^{1+\alpha}\int\phi(x)L^x e^{-\lambda L^x}\,dx\equiv\int\phi(x)\ell^\lambda(x)\,dx,$$
so that $\cL^\lambda$ becomes concentrated on $F$ as $\lambda\to\infty$.  It is not hard to use Theorem~\ref{lefttailthm} to see that for $|x|\ge \veps_0$,
\begin{equation}\label{1stmomentL}
c(\veps_0)|x|^{-p}\le \Ee_{\delta_0}(\ell^\lambda(x))\le C|x|^{-p},
\end{equation}
and Proposition~\ref{L2upperbound} will imply that for $|x_1|,|x_2|\ge\veps_0$,
\begin{equation}\label{2ndmomentL}
\Ee_{\delta_0}(\ell^\lambda(x_1)\ell^\lambda(x_2))\le C(\veps_0)[1+|x_1-x_2|^{2-p}].
\end{equation}
We conjecture that $\cL^\lambda$ converges in probability in the space $M_F(\Re^d)$ to a finite measure $\cL$ which necessarily is supported on $F$ and satisfies \eqref{EC} (by \eqref{2ndmomentL}). Although \eqref{1stmomentL} and \eqref{2ndmomentL} are not sufficient for this convergence, they do allow one to establish $\Pe_{\delta_0}(\text{dim}(F)\ge d+2-p)>0$ by first obtaining a lower bound on $\Pe_{\delta_0}(F\cap A\neq\emptyset)$ in terms of the $p-2$ capacity of $A$ (Theorem~\ref{polarsetslb}).

In Section~\ref{lbh2} we complete the proof of Theorem~\ref{dimthm} by showing the lower bound on dim$(F)$ in fact holds {\it with probability one} (Theorem~\ref{genlowerboundl}).  The lack of monotonicity of $F$ in the initial mass makes this step surprisingly delicate.  It uses the exit measures from half spaces and their special Markov property   (Corollary 8 of \cite{bib:leg95}) to analyze the size of $F$ in a half-space as
the bounding hyperplane moves across $\cR$.  

In Section~\ref{secproofs} we finish the proofs of Theorems~\ref{dimN0}, \ref{lefttailthm}, and \ref{dimX0} which now proceed
rather quickly, and also establish Propositions~\ref{noexit} and \ref{FLeb}. In Section~\ref{sec:9} we provide the proof of Proposition~\ref{L2upperbound}.
 In Section~\ref{prel} we introduce a number of our tools including some properties of $V^\lambda$, Yor's Girsanov theorem, and some properties of exit measures such as Le Gall's special Markov property and a slight modification thereof (Proposition~\ref{SpecialMProp}). 
\paragraph{\bf Convention on Functions and Constants.}
Constants whose value is unimportant and may change from line to line are denoted $C, c, c_d, c_1,c_2,\dots$, while constants whose values will be referred to later and appear initially in say, Lemma~i.j are denoted $c_{i.j},$ or $ \underline c_{i.j}$ or $C_{i.j}$. 
\medskip
\paragraph{Acknowledgements.} LM thanks Paul Balan\c{c}a for enjoyable and very helpful discussions on this problem.
 We thank L. Ryzhik for showing us the proof of Proposition~\ref{convrate}, and I. Benjamini whose comments about the boundary of the range of a scaling limit of tree-indexed random walks gave us a strong motivation to carry out  this work.
  EP thanks the Technion for hosting him during a visit where some of this research was carried out.   
 LM thanks UBC and EP  for hosting him during his visits. 

\section{Preliminaries}\label{prel}
\setcounter{equation}{0}
\setcounter{theorem}{0}

For each $\lambda>0$ there is a unique solution, $V^\lambda(x)$ to 
\eqref{vlequation},
where the equation is interpreted in the distributional sense.  Moreover $V^\lambda$ is $C^2$ on $\{x\neq 0\}$ and satisfies the (strong form) of \eqref{vlequation} on $\{x\neq 0\}$.  See p. 187 of \cite{bib:brez86}), and the references given there, for the above results.  Set 
$$g_0(x)=\begin{cases}1&\text{ if }d=1\\
\log^+(1/|x|)&\text{ if }d=2\\
|x|^{-1}&\text{ if }d=3,
\end{cases}$$
and let
$$p_t^{x}(y)=p_t(y-x).$$
Then for $d\ge 2$, (see p. 187 in \cite{bib:brez86}, or the more precise results in Remark 1 in \cite{bib:brezos87} and  \cite{bib:hong17})
\begin{equation}\label{Vlasymp}
\lim_{x\to 0}\frac{V^\lambda(x)}{\lambda g_0(x)}=c_d>0.
\end{equation}
and for $d=1$, $V^\lambda$ is continuous at $x=0$ (see Theorem 3.1 of \cite{bib:ver81}).
A simple extension of \eqref{LF} with $t=\infty$, $f=\lambda\delta_x$, and $\phi=0$ leads to the following result whose proof gives a self-contained proof of the existence of solutions to \eqref{vlequation} which are smooth on $\{x\neq0\}$.

\begin{lemma} \label{LTVL}For any $X_0\in M_F(\IR^d)$ and $\lambda>0$,
\begin{equation}\label{LTLT} 
\Ee_{X_0}(\exp(-\lambda L^x))=\exp\Bigl(-\int V^\lambda(x-x_0)X_0(dx_0)\Bigr).
\end{equation}
The above is strictly positive for all $x$ which are continuity points of $x\to \int g_0(y-x)dX_0(y)$, and in particular for $x\notin \text{Supp}(X_0)$. Moreover there is a 
$c_{\ref{LTVL}}$ so that, 
\begin{equation}\label{Vlpotential}
V^\lambda(x)\le c_{\ref{LTVL}}(\lambda g_0(x)+1)\ \forall x\in\IR^d\setminus\{0\}.
\end{equation}
\end{lemma}
\paragraph{Proof.} The strict positivity of the left-hand side of \eqref{LTLT} under the given condition on the potential of $X_0$ is immediate from the continuity of $L$ in \cite{bib:sug89} noted above since the above condition on $x$ easily implies (by dominated convergence) the joint continuity of $(t',x')\to X_0q_{t'}(x')$ at $(t,x)$ for any $t\ge 0$.  

Let $\{r_\veps:\veps\in(0,1)\}$ be a smooth, radially symmetric approximate identity so that $\{r_\veps>0\}\subset B(0,\sqrt\veps)$ and $r_\veps\le c_0p_\veps$ (the construction of $r_\veps$ is elementary).  Set $r_\veps^x(y)=r_\veps(y-x)$.  
We have from \eqref{LF} and for all $x\in\IR^d$,
\begin{equation}\label{prepreLTE}
\Ee_{X_0}\Bigl(\exp\Bigl(-\lambda\int_0^t X_s(r_\veps^x)ds\Bigr)\Bigr)=\exp\Bigl(-X_0(V_t(0,\lambda r_\veps^x))\Bigr).
\end{equation}
Theorem 3.3 of \cite{bib:iscoe86}  and symmetry imply
\begin{equation}
\label{eq:11_10_1}
V_t(0,\lambda r_\veps^x)(x_0)=V_t(0,\lambda r_\veps)(x-x_0)\uparrow V^{\lambda,\veps}(x-x_0)\text{ as }t\uparrow\infty,
\end{equation}
uniformly in $x$, where $V^{\lambda,\veps}$  is smooth and satisfies
\begin{equation}\label{preLTE}
\frac{\Delta V^{\lambda,\veps}}{2}=\frac{(V^{\lambda,\veps})^2}{2}-\lambda r_\veps, \quad V^{\lambda,\veps}>0 \text{ on }\IR^d.
\end{equation}

Let $t\to\infty$ in \eqref{prepreLTE} to conclude that 
\begin{equation}\label{preLTEb}
\Ee_{X_0}\Bigl(\exp\Bigl(-\lambda\int L^y r_\veps(y-x)dy\Bigr)\Bigr)=\Ee_{X_0}\Bigl(\exp\Bigl(-\lambda\int_0^\infty X_s(r_\veps^x)ds\Bigr)\Bigr)=\exp\Bigl(-\int V^{\lambda,\veps}(x-x_0)dX(x_0)\Bigr).
\end{equation}
Now set $X_0=\delta_{x_0}$ in the above.
Note that for $x\neq x_0$, the a.s. continuity of $L^y$ at $x$ implies 
the lefthand side of \eqref{preLTEb} converges to $\Ee_{\delta_{x_0}}(\exp(-\lambda L^x))\in(0,1)$ as $\veps\to 0$. Therefore the convergence of the righthand side implies that $V^{\lambda,\veps}(x) $ converges pointwise on $\{x\neq 0\}$ to a finite limit, $V^{\lambda,0}(x)>0$, as $\veps\to 0$.  So letting $\veps\downarrow 0$ in \eqref{preLTEb} with $x\neq x_0$ gives 
\begin{equation}\label{PreLTEc}
\Ee_{\delta_{x_0}}\Bigl(\exp(-\lambda L^x)\Bigr)=\exp(-V^{\lambda,0}(x-x_0))\text{ for }x_0\neq x.
\end{equation}
We have $\Delta V^{\lambda,\veps}=(V^{\lambda,\veps})^2$ for $|x|>\sqrt\veps$, and so Proposition 9 of Chapter V in \cite{bib:leg99} implies that $V^{\lambda,0}$ is $C^2$ on $\IR^d-\{0\}$ and solves $\Delta V^{\lambda,0}=(V^{\lambda,0})^2$.  

Next we show that $V^{\lambda,0}$ satisfies the distributional form of \eqref{vlequation}. 
For this it suffices to show 
\begin{equation}\label{vlebound}
V^{\lambda,\veps}(x)\le c(\lambda g_0(x)+1)\text{ for all }x\in\IR^d,\veps\in(0,1),
\end{equation}
because $g_0^2$ is locally integrable and so \eqref{vlebound} allows us to  take  limits as $\veps\downarrow 0$ in the distributional form of \eqref{preLTE}. 
Turning to \eqref{vlebound}, first note that by the semigroup property,  for any $t>1$
 \begin{eqnarray}
\nonumber
V_t(0,\lambda r_\veps^x)&=&V_1(V_{t-1}(0,\lambda r_\veps^x),\lambda r_\veps^x)\\
\nonumber
&\leq& \lim_{n\rightarrow\infty} V_1(n,\lambda r_\veps^x)\\
\label{eq:11_10_2}
&\leq& \lim_{n\rightarrow\infty} V_1(n,0)+V_1(0,\lambda r_\veps^x),
\end{eqnarray}
where the the first inequality follows by the boundedness of $V_{t-1}(0,\lambda r_\veps^x)$ and monotonicity of $V$ 
in its boundary conditions, and the last inequality follows by the subadditivity of $V$ in its boundary conditions (for the 
derivation of a similar result see, e.g.,~Lemma~2.6(b) in~\cite{bib:m02}). Here $n$ denotes the constant function $n$. Now recalling that $V_1(n,0)(\cdot)=\frac{n}{1+n/2}\leq 2$ uniformly in $n$, and 
$$V_1(0,\lambda r_\veps^x)(0)\leq \int_0^1P_s(\lambda r_\veps^x)(0)ds\leq \lambda c_0 \int_0^1 p_{s+\veps}(x)ds\leq c\lambda g_0(x),\;\forall \veps\in (0,1),$$
 we may use the above with~\eqref{eq:11_10_1} and \eqref{eq:11_10_2} to get~\eqref{vlebound}.

 The uniqueness of $V^\lambda$ now shows that $V^\lambda=V^{\lambda,0}$, and so \eqref{Vlpotential} is immediate from \eqref{vlebound} and
the obvious mononicity of $V^\lambda$.  If $x=x_0$ and $d\ge 2$, then both sides of \eqref{PreLTEc} are zero (recall \eqref{Vlasymp}) and so \eqref{LTLT} holds for all $x, x_0$ if $X_0=\delta_{x_0}$.  If $d=1$, we may let $x\to x_0$ on both sides of \eqref{PreLTEc} to obtain
the equality for $x=x_0$ by the continuity of $V^\lambda$ at $0$, and $L^x$ at $x_0$.

For \eqref{LTLT} with general $X_0$, the decomposition \eqref{PPPcan} also  applies to the total local time at $x$ and so 
\begin{equation}\label{LTCM}
\Ee_{X_0}(\exp(-\lambda L^x))=\exp\Bigl(-\int\int1-e^{-\lambda L^x(\nu)} \N_{x_0}(d\nu)X_0(dx_0)\Bigr).
\end{equation}
Taking $X_0=\delta_{x_0}$ in the above and the established \eqref{LTLT} gives 
\begin{equation}\label{vcanonical}V^\lambda(x-x_0)=\int1-e^{-\lambda L^x} d\N_{x_0}.
\end{equation}
Next use the above in \eqref{LTCM} to derive \eqref{LTLT} under $\Pe_{X_0}$.
This completes the proof.
\eop

We next give a modulus of continuity result for $L^x$ which follows easily from the moment bounds in \cite{bib:sug89}, the standard proof of Kolmogorov's continuity criterion, and tail bounds on the extinction time.
\begin{theorem}\label{LTmodulus} 
If $\veps_0\in(0,1)$ and $0<\gamma<((4-d)/2)\wedge 1$, there is a $\rho_{\veps_0,\gamma}(\omega)=\rho(\omega)>0$ $\Pe_{\delta_0}$-a.s. such that 
\[|L^{x_1}-L^{x_2}|\le |x_1-x_2|^\gamma\text{ whenever }|x_1-x_2|\le \rho,\text{ and }\veps_0/2<|x_i|<2\veps_0^{-1}.\]
Moreover, there are positive $\kappa(\gamma)$ and $C_{\ref{LTmodulus}}(\veps_0,\gamma)$ such that 
for all $0<t<1$, $\Pe_{\delta_0}(\rho<t)\le C_{\ref{LTmodulus}} t^\kappa$.
\end{theorem}
In fact one can take $\kappa(\gamma)=\frac{(4-d)/2-\gamma}{4-d-\gamma}$.

A simple scaling argument either for $L^x$ under $\Pe_{\delta_0}$, or directly for solutions of \eqref{vlequation}, shows that 
\begin{equation}\label{Vlscaling} 
V^{\lambda}(x)=V^\lambda(|x|)=r^{-2}V^{\lambda r^{4-d}}(x/r)\quad\forall x\in\IR^d,\ r>0.
\end{equation}
We may let $\lambda\uparrow\infty$ in \eqref{LTLT} under $X_0=\delta_0$ to see that $V^\lambda(x)\uparrow V^\infty(x)$ for all $x$ where 
\begin{equation}\label{Vinfty}
e^{-V^\infty(x)}=\Pe_{\delta_0}(L^x=0)>0,\ \text{ the latter for all } x\neq0,
\end{equation}
and for general $X_0$ that 
\begin{equation}\label{VinftygenX0}
\Pe_{X_0}(L^x=0)=\exp\Bigl(-\int V^\infty(x-x_0)dX_0(x_0)\Bigr).
\end{equation}
Take $r=|x|$ in \eqref{Vlscaling} and let $\lambda\uparrow\infty$ to conclude
$$V^\infty(x)=|x|^{-2}V^\infty(1) \text{ for all }x, \text{ where }0<V^\infty(1)<\infty.$$
Now it follows from \eqref{vlequation} and Proposition V.9(iii) of \cite{bib:leg99} that 
\begin{equation}\label{vequation}
\Delta V^\infty=(V^\infty)^2\text{ on } \{x\neq 0\},
\end{equation}
(or see \cite{bib:brez86} where $V^\infty=W$ in that work by (i) on p. 187).  The only possible positive value of $V^\infty(1)$ which will lead to a solution  of the above is $2(4-d)$ and so may conclude
\begin{equation}\label{vinftylimit}
V^\lambda(x)\uparrow V^\infty(x)=\frac{2(4-d)}{|x|^2}\text{ as }\lambda\to\infty\ \forall x.
\end{equation}

\begin{lemma}\label{Vlambdaprop}

\noi(a) $r^{d-1}(V^\lambda)'(r)$ is strictly increasing on $(0,\infty)$.

\noi(b) $V^{\lambda}(r)$ is strictly decreasing in $r>0$ to $0$.
\end{lemma}
\paragraph{Proof.} Let $v(r)=V^{\lambda}(r)$ for $r>0$.\\ 
\noi(a) Multiply the radial form of \eqref{vlequation} by $r^{d-1}$ to conclude $(r^{d-1}v')'=r^{d-1}v^2(r)>0$.\\
\noi(b) That $v(r)\to 0$ as $r\to\infty$ is immediate from $v(r)\le V^\infty(r)$. Suppose $v'(r_0)\ge0$ for some $r_0>0$. Then by (a), $v'(r)>0$ for all $r>r_0$ which, together with $v>0$, contradicts $\lim_{r\to\infty}v(r)=0$.
\eop

 Let $\rho(t)$ denote a $\delta$-dimensional Bessel process starting at $r>0$ under $P_r^{(\delta)}$, and $(\Ft t)$ be the filtration generated by $\rho$. For $R>0$, set $\tau_R=\inf\{t\ge 0:\rho(t)\le R\}$. The following result will be used frequently.
\begin{proposition}\label{yorthm}
Let $\lambda\ge 0, \mu> -1/2, r>0$ and $\nu=\sqrt{\lambda^2+\mu^2}$.  If $\Phi_t\ge 0$ is $(\Ft t)$-adapted, then for all $R<r$,
\[E_r^{(2+2\mu)}\Bigl(\Phi_{t\wedge \tau_R}\exp\Bigl(-\frac{\lambda^2}{2}\int_0^{t\wedge\tau_R}\frac{1}{\rho_s^2}\,ds\Bigr)\Bigr)=r^{\nu-\mu}E_r^{(2+2\nu)}\Bigl((\rho_{t\wedge\tau_R})^{-\nu+\mu}\Phi_{t\wedge\tau_R}\Bigr).\]
\end{proposition}

This follows from a simple application of Girsanov's theorem (see Section 2 and especially equation (2.b')  of \cite{bib:yor92}) to find $[dP^{(2+2\nu)}/dP^{(2)}]|_{\Ft t}$,
and $[1(T_0>t)dP^{(2+2\mu)}/dP^{(2)}]|_{\Ft t}$.
Once this  is established for the constant time $t$, a simple conditioning argument using optional stopping gives the above result at the stopping time $t\wedge\tau_R$.

We give a brief description of Le Gall's construction of $X$ using the Brownian snake.  Consider an initial condition $X_0\in M_F(\Re^d)$.  We set $\cW=\cup_{t\ge 0}C([0,t],\Re^d)$, call $\zeta(w)=t$ the lifetime of $w\in C([0,t],\Re^d)\subset\cW$, and metrize $\cW$ in the natural manner (see p. 54 of \cite{bib:leg99}).
The Brownian snake $W$ is the continuous $\cW$-valued strong Markov process constructed in Ch. IV of \cite{bib:leg99} and $\N_x$ is the excursion measure of $W$ starting at the path $x\in\Re^d$ with $0$ lifetime.   The construction of $X=X(W)$, first under $\N_x$ and then $\Pe_{X_0}$ is described in Theorem~4 of Ch. IV of \cite{bib:leg99}. Fortunately for our notation, the law of $X(W)$ under $\N_x$ is the canonical measure of super-Brownian motion described in the previous section (and also denoted by $\N_x$). For our purposes it will be important to note that if $\Xi=\sum_{i\in I}\delta_{W_i}$ is a Poisson point process on  $C([0,\infty),\cW)$ (the space of continuous $\cW$-valued paths) with intensity $\N_{X_0}(dW)=\int \N_x(dW)X_0(dx)$, then 
\begin{equation}\label{exdecompX}
X_t=\sum_{i\in I}X(W_i)_t=\int X_t(W)d\Xi(W),\ t>0
\end{equation} 
is a super-Brownian motion with initial state $X_0$, whose local time $L^x$ may therefore be decomposed as
\begin{equation}\label{exdecompL}
L^x=\sum_{i\in I} L^x(W_i)=\int L^x(W)d\Xi(W),
\end{equation}
where $x_i=W_{i,0}(0)$ and $(L^x(W_i))_x$ is a.s. continuous on $\{x\neq x_i\}$. 
We will refer to this construction as the standard setup for $X$ under $\Pe_{X_0}$ in what follows.

Under $\N_{X_0}$, if $\sigma=\inf\{t:\zeta_t=0\}$, $\zeta_t=\zeta(W_t)$ and $\hat W_t=W_t(\zeta(t))$ is the ``tip" of the snake
then (see, e.g. p. 70 of \cite{bib:leg99}) for $\phi\ge 0$, 
\begin{equation}\label{manhours}
\int _0^\infty X_t(\phi) dt=\int_0^\sigma \phi(\hat W_t)dt.
\end{equation}

Now fix an open set $G$ such that $d(\text{Supp}(X_0),G^c)>0$ and a Brownian path starting from any $x\in \partial G$ will exit $G$ immediately--in fact in what follows $G$ will be an open half-space or the complement of a closed ball.  We denote the exit measure of $X$ from $G$ by $X_G$. We refer the reader to Ch. V of \cite{bib:leg99} for the construction of $X_G$. Recall that $X_{G}$ is a random finite measure supported on $\partial G$,  which intuitively corresponds to the mass started at $X_0$ which is stopped at the instant it leaves $G$ (see p. 77 of \cite{bib:leg99}). Its Laplace functional is given by 
\begin{equation}\label{GEMLT}
\Ee_{X_0}\Bigl(\exp(-X_{G}(g))\Bigl)=\exp\Bigl(-\int U^{g}(x)X_0(dx)\Bigr),
\end{equation}
where $g:\partial G\to[0,\infty)$ is bounded and continuous, and $U^{g}\ge0$ is the unique continuous function on $\overline G$, which is $C^2$ on $G$,
and solves
\begin{equation}\label{GEMDE}
\Delta U^{g}=( U^{g})^2\text{ on }G,\  U^{g}=g \text{ on }\partial G.
\end{equation}
For this, see Theorem 6 in Chapter V of \cite{bib:leg99}, and the last exercise on p. 86 for uniqueness.
Note here, and elsewhere, we have taken our branching rate for $X$ to be one and so our constants will differ from those in \cite{bib:leg99}.  Although Chapter V of \cite{bib:leg99} deals with the excursion measure, $\N_x$, for the Brownian snake $W$, the extensions to $\N_{X_0}$ are immediate.  The same definitions also apply under $\Pe_{X_0}$, or equivalently, just set 
\begin{equation}\label{excdefexit}X_G=\sum_{i\in I}X_G(W_i)=\int X_G(W)d\Xi(W),
\end{equation}
where $\Xi$ is as in \eqref{exdecompX}.

We next state a version of the Special Markov Property for $W$ from \cite{bib:leg95}.  
First working under $\N_{X_0}$, define (as in \cite{bib:leg95})
\[S_G(W_u)=\inf\{t\le \zeta_u:W_u(t)\notin G\},\]
\begin{equation}\label{etaiGdef}
\eta^{G}_s=\inf\Bigl\{t:\int_0^t 1(\zeta_u\le S_G(W_u))du>s\Bigr\},
\end{equation}
and
\begin{equation}
\cE_G=\sigma\Bigl(W_{\eta^{G}_s},s\ge 0\Bigr)\vee\{\N_{X_0}-\text{null sets}\}.
\end{equation}
The above time-changed snake is in fact continuous in $s$ (see p. 401 of \cite{bib:leg95}).
Let \[\{u\ge 0:S_G(W_u)<\zeta_u\}=\cup_{i\in I}(a_i,b_i),\] 
for some countable set $I$.  Then $S_G(W_u)=S^i_G$ for all $u\in [a_i,b_i]$.  Let 
\[W^i_s(t)=W_{(a_i+s)\wedge b_i}(S^i_G+t)\text{ for }0\le t\le \zeta_{(a_i+s)\wedge b_i}-S^i_G,\]
so that $W^i\in C(\Re_+,\cW)$ are the excursions of $W$ outside $G$.

Proposition~2.3 of \cite{bib:leg95} implies that (under $\N_{X_0}$) $X_G$ is $\cE_G$-measurable.  Recall from \cite{bib:sug89} that $L=(L^y,y\in \overline{G}^c)$ is $\Pe_{X_0}$-a.s., and hence also $\N_{X_0}$-a.e., in the space $C( \overline{G}^c)$ of continuous functions on $\overline{G}^c$. 

\begin{proposition}\label{SpecialMProp}
Let $\Psi$ be a bounded measurable function on $C(\overline{G}^c)$ and $\Phi$ be a bounded measurable function on $M_F(\Re^d)$. Then

\noindent(a) Under $\N_{X_0}$ and conditional on $\cE_G$, $\sum_{i\in I}\delta_{W^i}$ is a Poisson point process with intensity $\N_{X_G}$.

\noindent(b) $\N_{X_0}(\Psi(L)|\cE_G)(\omega)=\Ee_{X_G(\omega)}(\Psi(L))\quad\N_{X_0}$-a.e., and therefore $$\N_{X_0}(\Phi(X_G)\Psi(L))=\N_{X_0}(\Phi(X_G)\Ee_{X_G}(\Psi(L))).$$

\noindent(c) $\Ee_{X_0}(\Phi(X_G)\Psi(L))=\Ee_{X_0}(\Phi(X_G)\Ee_{X_G}(\Psi(L)))$.
\end{proposition}
\paragraph{Proof.} (a) This is Corollary~2.8 of \cite{bib:leg95}, the extension from $\N_x$ to $\N_{X_0}$ being trivial.

\noindent(b) We only need show the first equality as the second is then immediate from the $\cE_G$-measurability of $X_G$.  If $S\subset {\overline G}^c$ (Borel) then $\N_{X_0}$-a.e.
\[\int_A L^xdx=\int_0^\sigma 1_A(\hat W_s)ds=\sum_{i\in I} \int_0^\infty 1_A(\hat W^i_s)ds,\]
where the last equality follows from an elementary calculation (using $\hat W^i_{b_i}\in \overline G$)
and the first equality is by \eqref{manhours}. It follows that 
\[ \text{for all $x\in{\overline G}^c$, }L^x=\sum_{i\in I}L^x(W^i),\]
where the existence of a continuous version of $L^x(W^i)$ on ${\overline G}^c$ follows from (a) and our 
earlier comments on $L$ under the excursion measure.  If $\Xi^{X_G(\omega)}$ is a Poisson point process on $C(\Re_+,\cW)$ with intensity $X_G(\omega)$ and integration is componentwise,
then (a) implies $\N_{X_0}$-a.e.
\[\N_{X_0}(\psi(L)|\cE_G)(\omega)=E\Bigl(\psi\Bigl(\int L(W)\Xi^{X_G(\omega)}(dW)\Bigr)\Bigr).\]
Comparing this to \eqref{exdecompL}, we see that the right-hand side of the above equals $\Ee_{X_G(\omega)}(\psi(L))$, and the proof of (b) is complete.

\noindent (c)   Use \eqref{exdecompL}, \eqref{excdefexit}, and the fact that $X_{G}=0$ implies $L=0$ $\N_{X_0}$-a.e. (by (b)) to see that $(X_{G},L)$ is equal in law to $\sum_{i=1}^M(X^i_{G},L_i)$, where $M$ is Poisson with mean $\N_{X_0}(X_{G}\neq0)<\infty$ (recall that $d(\text{Supp}(X_0), G^c)>0$), $(X^i_{G},L_i)_{i\ge 1}$ are iid with law $\N_{X_0}((X_{G},L)\in\cdot|X_{G}\neq0)$, and $M$ is independent of this iid sequence.  The result then follows easily, once we note that $\N_{X_0}(L\in\cdot\ |X_{G})=\Pe_{X_{G}}(L\in \cdot)$ (by (b)) and so by the multiplicative property of superprocesses,
\begin{equation*}
\Ee\Bigl(\Psi\Bigl(\sum_1^ML_i\Bigr)\Bigl|(X^i_{G})_{i\le M}, M\Bigr)=\Ee_{\sum_1^MX^i_{G}}\Bigl(\Psi(L)\Bigr),
\end{equation*}
which in turn implies (condition on $\sum_{i=1}^M X^i_{G}$) 
\[\Ee_{X_0}(\Psi(L)|X_G)=\Ee_{X_G}(\Psi(L)),\]
as required.
\eop

It is not hard to prove a full analogue of (a) under $\Pe_{X_0}$ but (c) will suffice for our purposes.

\section{The Upper Bound on the Dimension} \label{secupperbound}
\setcounter{equation}{0}
\setcounter{theorem}{0}
The goal of this section is to prove the following:
\begin{theorem}\label{dimub2} With $\Pe_{\delta_0}$-probability one, 
$\textnormal{dim}(F)\le d+2-p$.
\end{theorem}

\medskip
\noi
First let us introduce additional notation. 

\noi{\bf Notation} $B_\veps=\{y\in\IR^d:|y|<\veps\}$, $G_\veps=\overline{B_\veps}^c$.

\medskip

By \eqref{GEMLT} and \eqref{GEMDE} if $g:\partial B_\veps\to[0,\infty)$ is bounded and continuous
the Laplace functional of the exit measure from $G_\veps$ is given by
\begin{equation}\label{EMLT}
\Ee_{\delta_x}\Bigl(\exp(-X_{G_\veps}(g))\Bigl)=\exp(-U^{g,\veps}(x)),\ \forall|x|>\veps,
\end{equation}
where $U^{g,\veps}\ge 0$ is the unique continuous function on $ {B^c_\veps}$ which is $C^2$ on $\{|x|>\veps\}$
and solves
\begin{equation}\label{EMDE}
\Delta U^{g,\veps}=( U^{g,\veps})^2\text{ on }\{|x|>\veps\},\  U^{g,\veps}=g \text{ on }\partial B_\veps.
\end{equation}




For constants $\lambda>0$, we will be particular interested in $U^{\lambda,\veps}$ which, thanks to the uniqueness in \eqref{EMDE}, satisfies
\begin{equation}\label{Uscale}
U^{\lambda,\veps}(|x|)=U^{\lambda,\veps}(x)=\veps^{-2}U^{\veps^2\lambda,1}(x/\veps)\text{ for }|x|\ge\veps.
\end{equation}
By \eqref{EMLT}, with $g=\lambda$, we have $U^{\lambda,\veps}\uparrow U^{\infty,\veps}$ as $\lambda\uparrow\infty$, where
\begin{equation}\label{UinfLT}
\Pe_{\delta_x}(X_{G_\veps}(1)=0)=\exp(-U^{\infty,\veps}(x))\ \ \forall|x|>\veps.
\end{equation}
Elementary properties of $X$ show that the left-hand side of \eqref{EMLT} is in $(0,1)$ and converges to $1$ as $|x|\to\infty$, and so
\begin{equation}\label{Uinfprop}
0<U^{\infty,\veps}(x)<\infty\ \ \forall|x|>\veps,\ \lim_{|x|\to\infty}U^{\infty,\veps}(x)=0,\text{ and so }\lim_{|x|\to\infty}U^{\lambda,\veps}(x)=0.
\end{equation}
Proposition~9(iii) of \cite{bib:leg99} implies that $U^{\infty,\veps}$ is $C^2$ and 
\begin{equation}\label{Uinfequation}
\Delta U^{\infty,\veps}=(U^{\infty,\veps})^2\text{ on }G_\veps,\ \ \lim_{|x|\to\veps,|x|>\veps}U^{\infty,\veps}(x)=+\infty,
\end{equation}
where the last limit is immediate from the definition of $U^{\infty,\veps}$.
A simple application of the maximum principle implies that
\begin{equation}\label{uinftylb}
U^{\infty,1}(x)\ge V^\infty(x)=\frac{2(4-d)}{|x|^2}\ \ \forall|x|>1.
\end{equation}
For $1\ge \delta_0>0$, let $D^{\delta_0}=U^{\infty,1}-U^{\delta_0,1}\ge 0$.  Then \eqref{Uinfequation} and \eqref{EMDE} imply $D^{\delta_0}$ is a $C^2$ solution of 
\begin{equation*}
\frac{\Delta }{2}D^{\delta_0}=\Bigl(\frac{U^{\infty,1}+U^{\delta_0,1}}{2}\Bigr)D^{\delta_0}\text{ on }\{|x|>1\}.
\end{equation*}
We let $B$ denote a $d$-dimensional Brownian motion starting at $x$ under $P_x$.  The above equation and the Feynman-Kac formula (e.g. p. 114 of \cite{bib:SV79}) implies that if $\tau_R=\inf\{t\ge 0:|B_t|\le R\}$ and $|x|\ge R>1$, then
\begin{equation*}
D^{\delta_0}(x)=E_x\Bigl(D^{\delta_0}(B_{t\wedge\tau_R})\exp\Bigl(-\int_0^{t\wedge\tau_R}\Bigl(\frac{U^{\infty,1}+U^{\delta_0,1}}{2}\Bigr)(B_s)\,ds\Bigr)\Bigr).
\end{equation*}
We may let $t\to\infty$ in the above and, noting that $D^{\delta_0}(B_{t\wedge\tau_R})\to 0$ as $t\to\infty$ when $\tau_R=\infty$ (for $d=3$), conclude
\begin{equation}\label{DFK}
D^{\delta_0}(x)=D^{\delta_0}(R)E_x\Bigl(1(\tau_R<\infty)\exp\Bigl(-\int_0^{\tau_R}\Bigl(\frac{U^{\infty,1}+U^{\delta_0,1}}{2}\Bigr)(B_s)\,ds\Bigr)\Bigr)\text{ for }|x|\ge R>1.
\end{equation}

\begin{lemma}\label{Ulambdaprop}

\noi(a) $r^{d-1}(U^{\lambda,\veps})'(r)$ is strictly increasing.

\noi(b) $U^{\lambda,\veps}(r)$ is strictly decreasing in $r\ge\veps$.

\noi(c) If $d=2$ or $3$, then $\lim_{r\to\infty} r(U^{\lambda,\veps})'(r)=0$.
\end{lemma}
\paragraph{Proof.} Let $u(r)=U^{\lambda,\veps}(r)$ for $r>\veps$.\\ 
\noi(a),(b). These are proved as for $V^\lambda$ in Lemma~\ref{Vlambdaprop}.\\
\noi(c) It follows from (a) and (b) that $r^{d-1}u'(r)\uparrow -c\le 0$ as $r\uparrow\infty$.  If $d=3$ we therefore have $ru'(r)=r^{-1}(r^{d-1}u'(r))\to 0$ as $r\to\infty$.  Assume now $d=2$, so we must show $c=0$. Assume not. Then $u'(r)\le -c/r$ for all $r\ge \veps$.  This implies 
$-u(\veps)=\int_\veps^\infty u'(r)dr\le -c\int_\veps^\infty r^{-1}dr=-\infty$, a contradiction.
\gdm

If $d=2$, $\delta_0>0$, and $y(t)=U^{\delta_0,1}(e^t)$ for $t\ge 0$, then Lemma~\ref{Ulambdaprop} and a simple calculation using the radial form of \eqref{EMDE} shows
\begin{equation}\label{yprop}y''(t)=e^{2t}y(t)^2,\ y(0)=\delta_0,\ y \text{ is decreasing},\ y'(t)\to 0\text{ as }t\to\infty.
\end{equation}
A direct calculation gives the same conclusions for $z(t)=V^\infty(e^t)$, where now $z(0)=2(4-d)=4$.  Theorem~1.1 of \cite{bib:tal78} therefore implies that
\begin{equation}\label{ratio1}
\lim_{|x|\to\infty}\frac{U^{\delta_0,1}(x)}{V^\infty(x)}=\lim_{t\to\infty}\frac{y(t)}{z(t)}=1.
\end{equation}
Similar reasoning applies if $d=3$, but now with $y(t)=(t+1)U^{\delta_0,1}(t+1)$ and $z(t)=(t+1)V^\infty(t+1)$ for $t\ge 0$.  This leads to $y''(t)=(t+1)^{-1}y(t)^2>0$ and $y'(t)=(t+1)(U^{\delta,1})'(t+1)+U^{\delta_0,1}(t+1)\to 0$ as $t\to\infty$ by Lemma~\ref{Ulambdaprop}(c). Note also that if $y'(t_0)\ge0$ for some $t_0\ge 0$, then $y''>0$ clearly contradicts the fact that $y'(t)\to 0$ and so \eqref{yprop} again holds with $(t+1)^{-1}$ in place of $e^{2t}$.  Again, a direct calculation gives the same conclusions for $z$ with the modified initial condition.  Hence \eqref{ratio1} again follows from Theorem~1.1 of \cite{bib:tal78}. Finally for $d=1$ it follows from Lemma~\ref{Ulambdaprop}(a) and the fundamental theorem of calculus that $\lim_{r\to\infty} (U^{\delta_0,1})'(r)=0$.  Therefore one can conclude \eqref{ratio1} directly from Theorem~1.1 of \cite{bib:tal78} with $y=U^{\delta_0,1}$ and $z=V^\infty$.

We will, however, need to quantify the convergence in \eqref{ratio1}.  

\begin{proposition}\label{convrate} Let $\delta_0\in(0,1)$.\\
\noi(a) For $|x|>1$, $U^{\delta_0,1}(x)\le V^\infty(x)$.\\
\noi(b) For any $\delta\in(0,1)$ there is a $C_\delta>1$ so that 
$U^{\delta_0,1}(x)\ge (1-\delta)V^\infty(x)$ for all $|x|\ge C_\delta/\delta_0$.
\end{proposition}
\paragraph{Proof.}(a) As $\delta_0<1$, clearly the required inequality holds if $|x|=\delta_0$. Equations \eqref{vequation}, \eqref{EMDE} and the last part of \eqref{Uinfprop}  now allow us to apply the maximum principle and conclude the inequality for all $|x|\ge 1$. 

\noi(b) In view of (a) and the monotonicity of $U^{\delta_0,1}$ in $\delta_0$ (recall \eqref{EMLT}), it clearly suffices to consider $\delta_0\le \delta_1$ for some fixed value of $\delta_1\in(0,1)$.  We will write $u(r)$ for $U^{\delta_0,1}(r)$ and $v(r)$ for $V^\infty(r)$. 

Consider first $d=2$.  For $t\ge 0$, let 
$$q(t)=\frac{u(e^{t/4})}{v(e^{t/4})}=\frac{1}{4}u(e^{t/4})e^{t/2}\in(0,1]\text{ (by (a))}.$$
(See the proof of Theorem 1.1 of \cite{bib:tal78} for the motivation for this change variables.)
A simple calculation using \eqref{EMDE} gives
\[ q''-q'+\frac{1}{4}(q-q^2)=0,\ q(0)=\delta_0/4,\ \lim_{t\to\infty}q(t)=1,\]
where the last limit is by \eqref{ratio1}.
Therefore if $w(t)=1-q(-t)$ for $t\le0$, we have
\begin{equation}\label{travwave}
\frac{1}{2}w''+\frac{1}{2}w'+\frac{1}{8}(w^2-w)=0\text{ for }t\le 0,\ \lim_{t\to-\infty}w(t)=0,\ w(0)=1-\frac{\delta_0}{4}, \ w\in[0,1).
\end{equation}
The above steps can be reversed and as the equation defining $u$ has a unique solution it follows that so does \eqref{travwave}. Solutions to the above give travelling wave solutions to the KPP
equation and, using the notation on p. 55 of \cite{bib:kyp04}), $c=\frac{1}{2}=\sqrt{2\beta}$ where $\beta=\frac{1}{8}$, so we are in the critical case in the above reference.  On p. 55 of the above reference an increasing solution $\phi$  to \eqref{travwave} is given, but with $\phi(0)=\frac{1}{2}$, defined on the line and satisfying 
\begin{equation}\label{twasymptotics}
\lim_{t\to\infty}(1-\phi(t))t^{-1}e^{t/2}=c_0>0.
\end{equation}
Define $R_{\delta_0}>0$ by $\phi(R_{\delta_0})=1-(\delta_0/4)$ and $L_\delta<0$ by $\phi(L_\delta)=\delta$. Use \eqref{twasymptotics} to conclude that for some $\delta_1>0$,
\[R_{\delta_0}<4\log(8c_0/\delta_0)\text{ for }0<\delta_0<\delta_1.\]
Uniqueness in \eqref{travwave} implies that $w(t)=\phi(R_{\delta_0}+t)$ for $t\le 0$ and therefore
$$1-\delta\le q(t)=1-w(-t)=1-\phi(R_{\delta_0}-t))\text{ iff }R_{\delta_0}- L_\delta\le t.$$
Recalling the definition of $q$, we see that this implies $u(r)\ge(1-\delta)v(r)$ providing 
$$r\ge \exp(\log(8c_0/\delta_0)-L_\delta/4)=C_\delta\delta_0^{-1}\text{ for }\delta_0<\delta_1.$$

Turning now to $d=3$, for $t\ge \log(1/12)$, let
$$q(t)=\frac{u({12}^{1/3}e^{t/3})}{v({12}^{1/3}e^{t/3})}=\frac{{12}^{2/3}}{2}u({12}^{1/3}e^{t/3})e^{2t/3}\in(0,1].$$
Then as before we have,
$$q''-q'+\frac{2}{9}(q-q^2)=0,\ q(\log(1/12))=\delta_0/2,\ \lim_{t\to\infty}q(t)=1,$$
and $w(t)=1-q(-t)$, for $t\le \log(12)$, satisfies
\begin{equation}\label{travwave3}\frac{1}{2}w''+\frac{1}{2}w'+\frac{1}{9}(w^2-w)=0,\ t\le \log 12,\ w(\log 12)=1-\frac{\delta_0}{2},\ w\in[0,1),\ \lim_{t\to-\infty}w(t)=0.
\end{equation}
Proceeding as before, now with $\beta=1/9$ and $c=1/2>\sqrt{2\beta}$, we are in the supercritical case in \cite{bib:kyp04} and the solution of \eqref{travwave3} on the line, which equals $1/2$ at $t=0$ (rather than the terminal condition in \eqref{travwave3}), $\phi$, satisfies
$$\lim_{t\to\infty}(1-\phi(t))e^{t/3}=c_1>0,$$
where the value of the exponential rate in the above may be found in Theorem 2.1 of \cite{bib:har99}. The argument is now completed for $d=3$, just as it was in $d=2$.

The argument for $d=1$ is similar to that for $d=3$.  For $t\ge \log(1/180)$, set 
$q(t)=\frac{u(180^{1/5}e^{t/5})}{v(180^{1/5}e^{t/5})}$,
so that $w(t)=1-q(-t)$, $t\le \log(180)$, satisfies
\[\frac{1}{2}w''+\frac{1}{2}w'+\frac{3}{25}(w^2-w)=0,\quad w(\log(180))=1-\frac{\delta_0}{6}, \lim_{t\to-\infty}w(t)=0.\]
We are again the supercritical case for the KPP equation with $\beta=3/25$ and $c=1/2>\sqrt{2\beta}$.  Now (Thm. 2.1 of \cite{bib:har99}) $\lim_{t\to \infty}(1-\phi(t))e^{2t/5}=c_2>0$.
The result now follows as above. In fact we get the desired conclusion under the weaker condition that $|x|\ge C_\delta/\sqrt{\delta_0}$.
\gdm

\noi{\bf Throughout the rest of this Section we fix $\veps_0\in (0,1)$.}

\begin{proposition}\label{exittail}
If $0<p'<p$, there is an $\eta=\eta(p')>0$ and $C_{\ref{exittail}}=C_{\ref{exittail}}(\veps_0,p')<\infty$ so that for all $\veps_0\le |x|$ and $\veps\in(0,\veps_0)$,
$$\Pe_{\delta_x}(0<X_{G_\veps}(1)<\veps^{2-\eta})\le C_{\ref{exittail}}\veps^{p'-2}.$$
\end{proposition}
\paragraph{Proof.} Choose $p'\in(0,p)$ and then $\eta,\delta\in(0,1)$.  Below we will choose $\eta,\delta$ sufficiently small, depending on $p'$. If $C_\delta$ is as in Proposition~\ref{convrate}, consider
\begin{equation}\label{epschoice}
0<\veps<(\veps_0/C_\delta)^{1/(1-\eta)}(<\veps_0),
\end{equation}
and set $\lambda'=\veps^{\eta-2}$. Then for $|x|>\veps_0$,
\begin{align}\label{exitmeasbnd1}
\nn\Pe_{\delta_x}(0<X_{G_\veps}(1)<\veps^{2-\eta})&\le e\Ee_{\delta_x}\Bigl(\exp(-\lambda' X_{G_\veps}(1))1(X_{G_\veps}(1)>0)\Bigr)\\
\nn&=e\Bigl(e^{-U^{\lambda',\veps}(x)}-e^{-U^{\infty,\veps}(x)}\Bigr)\ \ (\text{by \eqref{EMLT} and \eqref{UinfLT}})\\
&\le e\veps^{-2}(U^{\infty,1}-U^{\veps^{\eta},1})(x/\veps)=e\veps^{-2}D^{\veps^\eta}(x/\veps),
\end{align}
where in the last line we use the scaling relation \eqref{Uscale} and the definition of $D$. Now use the Feynman-Kac representation in \eqref{DFK} with $R=C_\delta\veps^{-\eta}$ ($C_\delta$ as in Proposition~\ref{convrate}). Note here that \eqref{epschoice} implies $|x|/\veps>\veps_0/\veps>R>1$, so that \eqref{DFK} applies and we conclude that
\begin{align*} 
\Pe_{\delta_x}(0<X_{G_\veps}(1)<\veps^{2-\eta})&\le e\veps^{-2}D^{\veps^\eta}(R)
E_{x/\veps}\Bigl(1(\tau_R<\infty)\exp\Bigl(-\int_0^{\tau_R}\frac{U^{\infty,1}+U^{\veps^\eta,1}}{2}(B_s)\,ds\Bigr)\Bigr)\\
&\le e\veps^{-2}D^{\veps^\eta}(R)E_{x/\veps}\Bigl(1(\tau_R<\infty)\exp\Bigl(-\int_0^{\tau_R}\Bigl(1-\frac{\delta}{2}\Bigr)V^\infty(B_s)\,ds\Bigr)\Bigr),
\end{align*}
the last by \eqref{uinftylb} and Proposition~\ref{convrate}(b). Next, apply Proposition~\ref{yorthm} to see that if 
\begin{equation}\label{mu}
\mu=\begin{cases} -1/2&\text{ if }d=1\\
0&\text{ if }d=2\\
\frac{1}{2}&\text{ if }d=3,
\end{cases}
\end{equation}
and 
\begin{equation}\label{nu}
\nu=\nu(\delta)=\Bigl(\mu^2+\Bigl(1-\frac{\delta}{2}\Bigr)4(4-d)\Bigr)^{1/2}\equiv(\mu^2+\lambda^2)^{1/2},
\end{equation}
then
\begin{align*}
\Pe_{\delta_x}(0<X_{G_\veps}(1)<\veps^{2-\eta})&\le  e\frac{2(4-d)}{R^2}\veps^{-2}
\lim_{t\to\infty}E^{(2+2\mu)}_{|x|/\veps}\Bigl(1(\tau_R<t)\exp\Bigl(-\frac{\lambda^2}{2}\int_0^{\tau_R\wedge t}\frac{1}{\rho_s^2}\,ds\Bigr)\Bigr)\\
&\le  e\frac{2(4-d)}{R^2}\veps^{-2}(|x|/\veps)^{\nu-\mu}
\lim_{t\to\infty}E^{(2+2\nu)}_{|x|/\veps}\Bigl(1(\tau_R<t)\rho_{t\wedge\tau_R}^{-\nu+\mu}\Bigr)\\
&=e\Bigl[2(4-d)R^{-\nu+\mu-2}\Bigr]\veps^{\mu-\nu-2}|x|^{\nu-\mu}P_{|x|/\veps}^{(2+2\nu)}(\tau_R<\infty)\\
&=e\Bigl[2(4-d)R^{-\nu+\mu-2}\Bigr]\veps^{\mu-\nu-2}|x|^{\nu-\mu}(|x|/\veps)^{-2\nu}R^{2\nu}\\
&\le e2(4-d)C_\delta^{\nu+\mu-2}\veps_0^{-\mu-\nu}\veps^{(\mu+\nu-2)(1-\eta)},
\end{align*}
where we also used the following result for  the hitting probabilities for Bessel processes:
\begin{equation}\label{besselexit}
P_r^{(2+2\xi)}(\tau_R<\infty)=(R/r)^{2\xi}\quad\forall\xi\ge0,\ r\ge R>0.
\end{equation}

A bit of arithmetic shows that $\lim_{\delta\downarrow 0}\mu+\nu(\delta)=p$ and so we can choose $\delta,\eta>0$ sufficiently small, depending on $p'$ so that $(\mu+\nu(\delta)-2)(1-\eta)\ge p'-2$.  This gives the required bound for $\veps<(\veps_0/C_\delta)^{1/(1-\eta)}\equiv\veps_1(\veps_0,p')$, and the result then follows for all $0<\veps<\veps_0$ by adjusting $C_{\ref{exittail}}=C_{\ref{exittail}}(\veps_0,p')$ accordingly.
\gdm

\begin{theorem}\label{dimub} If $p'\in(0,p)$, there is a $C_{\ref{dimub}}=C_{\ref{dimub}}(\veps_0,p')<\infty$ such that 
\[\Pe_{\delta_0}(B(x,\veps)\cap F\neq\emptyset)\le C_{\ref{dimub}}\veps^{p'-2}\text{ whenever }\veps_0\le |x|\le\veps_0^{-1},\ \veps\in(0,\veps_0/2).\]
\end{theorem}
\paragraph{Proof.} Let $0<p'<p$ and select $\eta>0$ as in Proposition~\ref{exittail}.  Assume $\veps_0$, $\veps$ and $x$ are as above.  Note that $B_\veps\cap F\neq\emptyset$ iff $\exists |x'|<\veps$ and $|x_n|<\veps$ so that $x_n\to x'$, $L^{x'}=0$ and $L^{x_n}>0$.  Also Proposition~\ref{SpecialMProp}(c) shows that 
\[\Pe_{\delta_x}(\exists|x'|<\veps\text{ such that }L^{x'}>0,X_{G_\veps}(1)=0)=0.\]
By translation invariance, followed by another application of Proposition~\ref{SpecialMProp}(c), we have
\begin{align}
\nonumber\Pe_{\delta_0}(B(x,\veps)\cap F\neq\emptyset)&=\Pe_{\delta_x}(B_\veps\cap F\neq\emptyset)\\
\nonumber&\le \Pe_{\delta_x}(1(X_{G_\veps}(1)>0)\Pe_{X_{G_\veps}}(B_\veps\cap F\neq\emptyset))\\
\nonumber&\le \Pe_{\delta_x}(0<X_{G_\veps}(1)<\veps^{2-\eta})+\Ee_{\delta_x}\Bigl(1(X_{G_\veps}(1)\ge\veps^{2-\eta})\Pe_{X_{G_\veps}}(B_\veps\cap F\neq\emptyset)\Bigr)\\
\label{coverbnd}&\le C_{\ref{exittail}}(\veps_0,p')\veps^{p'-2}+\Ee_{\delta_x}\Bigl(1(X_{G_\veps}(1)\ge\veps^{2-\eta})\Pe_{X_{G_\veps}}(\exists|x'|<\veps\text{ such that }L^{x'}=0)\Bigr),
\end{align}
the last by Proposition~\ref{exittail}. For the second term above, if $\beta>(4-d)4$ consider a collection of points $\{x_i:1\le i\le M_\veps\}\subset B_\veps$ such that $B_\veps\subset\cup_{i=1}^{M_\veps}B(x_i,\veps^\beta)$, $M_\veps\le c_d\veps^{-d\beta}$.  Next set $\gamma=1/4$ in Theorem~\ref{LTmodulus} and let $C_{\ref{LTmodulus}}(\veps_0,1/4)\equiv C_{\ref{LTmodulus}}(\veps_0)$ and $\kappa(\veps_0,1/4)\equiv\kappa(\veps_0)$ be as in that result. If $\lambda=\veps^{-(4-d)}$, then for each $i$,
\begin{align}
\nonumber\Pe_{X_{G_\veps}}(L^{x_i}\le \veps^{4-d})&\le e\Ee_{X_{G_\veps}}(e^{-\lambda L^{x_i}})\\
\nonumber&=e\exp\Bigl(-\int V^\lambda(x-x_i)dX_{G_\veps}(x)\Bigr)\quad(\text{by Lemma~\ref{LTVL}})\\
\nonumber&=e\exp\Bigl(-\veps^{-2}\int V^1((x-x_i)/\veps)dX_{G_\veps}(x)\Bigr)\quad(\text{by \eqref{Vlscaling}})\\
\label{smallLbnd}&\le e\exp(-\veps^{-2}V^1(2)X_{G_\veps}(1)),
\end{align}
the last since $X_{G_\veps}$ is supported on $\partial B_\veps$, $|x-x_i|/\veps\le2$ for $|x|=\veps$, and $V^1(r)$ is decreasing in $r$ (Lemma~\ref{Vlambdaprop}). Use \eqref{smallLbnd} above and then the definition of $\rho_{\veps_0,1/4}\equiv \rho_{\veps_0}$ in Theorem~\ref{LTmodulus} to bound the second term in \eqref{coverbnd} by
\begin{align}
\nonumber\Ee_{\delta_x}(&1(X_{G_\veps}(1)\ge\veps^{2-\eta})\Pe_{X_{G_\veps}}(L^{x_i}>\veps^{4-d}\ \forall i\le M_\veps, L^{x'}=0\ \exists|x'|<\veps)\\
\nonumber&\qquad+M_\veps\Ee_{\delta_x}(1(X_{G_\veps}(1)\ge\veps^{2-\eta})e\exp(-\veps^{-\eta}V^1(2)))\\
\nonumber&\le \Pe_{\delta_x}(L^{x_i}>\veps^{4-d}\ \forall i\le M_\veps,\ L^{x'}=0\ \exists|x'|<\veps)\\
\nonumber&\qquad+c_d\veps^{-\beta d}e\exp(-\veps^{-\eta}V^1(2))\quad(\text{by Proposition~\ref{SpecialMProp}(c) again})\\
\label{term2bnd}&\le\Pe_{\delta_0}(\rho_{\veps_0}\le\veps^\beta)+c_de\veps^{-\beta d}\exp(-\veps^{-\eta}V^1(2)),
\end{align}
where the last line follows after a bit of arithmetic after translation by $-x$. Now bound the first term in \eqref{term2bnd} using Theorem~\ref{LTmodulus}, and hence bound \eqref{coverbnd} by
\[C_{\ref{exittail}}(\veps_0,p')\veps^{p'-2}+C_{\ref{LTmodulus}}(\veps_0)\veps^{\beta\kappa}+c_de\veps^{-\beta d}\exp(-\veps^{-\eta}V^1(2)).\]
Finally choose $\beta=\beta(\veps_0)>4(4-d)$ so that $\beta\kappa(\veps_0)>p-2$.  The result is then immediate from the above bound.
\eop


Now we are ready to complete 
\paragraph{Proof of Theorem~\ref{dimub2}}
Let $p'\in (0,p)$, $d_f>d+2-p'$, and fix $\veps_0\in(0,1)$.  For $\veps\in(0,\veps_0/2)$, cover $A=\{\veps_0<|x|<\veps_0^{-1}\}$ with $N_\veps\le C(\veps_0)\veps^{-d}$ open balls $\{B_i\}$ of radius $\veps>0$ centered in $A$.  If $\veps=\veps_n\downarrow 0$, then by Fatou's Lemma and Theorem~\ref{dimub}, 
\[\Ee_{\delta_0}\Bigl(\liminf_{n\to\infty}\sum_{i=1}^{N_\veps}\veps_n^{d_f}1(B_i\cap F\neq\emptyset)\Bigr)\le \liminf_{n\to\infty}C'(\veps_0,p')\veps_n^{d_f-d+p'-2}=0.\]
So letting $p'\uparrow p$, $d_f\downarrow d+2-p$, and finally $\veps_0\downarrow 0$, gives the required upper bound on the Hausdorff dimension of $F$.
\eop

\section{
Proof of Theorem~\ref{onedims}} \label{sec1d}
\setcounter{equation}{0}
\setcounter{theorem}{0}
We assume $d=1$ throughout this section and prove Theorem~\ref{onedims}, and hence that $F$ is a two-point set $\Pe_{\delta_0}$-a.s.
It is well known that the range of super-Brownian motion $\cR\equiv\overline{\{x:L_x>0\}}$ is an interval (see, e.g. Theorem 7 of Chapter IV of \cite{bib:leg99}) but this does not imply Theorem~\ref{onedims}, i.e., that $L$ is strictly positive on the interior of the range.  

We will work with a one-dimensional super-Brownian motion $X$ with initial state $\y\delta_0$ defined from a Brownian snake as in Section~\ref{secupperbound}.  For $r>0$ let $Y_r\delta_r$ denote the exit measure from $(-\infty,r)$, and set $Y_0=\y$. The same notation is used for $r>0$ when working under the excursion measure, $\N_0$, for the snake $W$. It follows from the Markov property 1.3D on p. 36 of \cite{bib:d02}, and the ensuing construction in Chapter 4 of the same reference, that $\{Y_r:r>0\}$ is a Markov process. That is, if $\psi:\Re_+\to \Re$ is bounded and measurable and $0<r_1<r_2$, then 
\begin{equation}\label{exitmarkov}
\Ee_{\y\delta_0}(\Psi(Y_{r_2})|(Y_r,\ r\le r_1))(\omega)=\Ee_{Y_{r_1}(\omega)\delta_0}(\Psi(Y_{r_2-r_1}))\
\Pe_{\y\delta_0}-a.s.
\end{equation}
With a bit of work one can also derive this from  Proposition~\ref{SpecialMProp}(b).  (One starts by decomposing $Y_r$ into the sum of the contributions from the excursions $W_i$ from $0$ as in \eqref{excdefexit}.)

Recall (see, e.g., Section II.1 of  \cite{bib:leg99}) that a stable continuous state branching process (SCSBP) with parameter $p\in(1,2)$ and scaling constant $c_0>0$ is a $[0,\infty)$-valued cadlag strong Markov process, $Z=\{Z_t:t\ge 0\}$, whose transition kernel $p_t(x,dy)$ satisfies
\begin{equation}\label{CSBP}
\int e^{-\lambda y}p_t(x,dy)=\exp(-xu^{\lambda}(t)), \text{for all }\lambda\ge 0,
\end{equation}
where
\begin{equation}\label{csbpode}
\frac{du^{\lambda}(t)}{dt}=-c_0u^{\lambda}(t)^p,\ \ u^{\lambda}(0)=\lambda.
\end{equation}

The following result is well-known (see, e.g., the comment prior to Proposition~3 in \cite{bib:al17} and
recall that our branching rates differ) but
 we include the elementary proof here for completeness.
\begin{proposition}\label{exitcsbp}
Under $\Pe_{\y\delta_0}$ there is a cadlag version of $Y$ (also denoted $Y$) which is a SCSBP starting at $\y$ with parameter $p=3/2$ and $c_0=\sqrt{6}/3$.  In particular $Y$ is a martingale, and if 
\[R=\inf\{r>0:Y_r\wedge Y_{r-}=0\},\]
then  $Y_r=0$ for all $r\ge R=\inf\{r>0:Y_r=0\}<\infty$ $\Pe_{\y\delta_0}$-a.s.
\end{proposition}
\paragraph{Proof.} 
If
\begin{equation}\label{YLTform}
U^{\lambda,r}(x)=U^{\lambda,0}(x-r)=6(r-x+\sqrt{6/\lambda})^{-2}, \text{ for }x\le r,
\end{equation}
then
\begin{equation}\label{Yemde}
\Delta U^{\lambda,r}(x)=(U^{\lambda,r}(x))^2\text{ on }\{x<r\},\ U^{\lambda,r}(r)=\lambda,
\end{equation}
and so from \eqref{GEMDE} and \eqref{GEMLT} we see that 
\begin{equation}\label{Yemlt}
\Ee_{\y\delta_0}(\exp(-\lambda Y_r))=\exp(-\y U^{\lambda,r}(0)).
\end{equation}
This shows that $Y_r\to \y$ as $r\downarrow 0$, first in law and hence in probability.  By \eqref{exitmarkov} this implies $\{Y_r:r\ge0\}$ is right continuous in probability.  We can now let $r_1\downarrow 0$ in \eqref{exitmarkov} to see that it continues to hold for $r_1=0$--first consider
$\psi(y)=e^{-\lambda y}$ and proceed by a monotone class argument.  
Differentiate \eqref{Yemlt} in $\lambda$ to see that $\Ee_{\y\delta_0}(Y_r)=\y$, and so by \eqref{exitmarkov} $\{Y_r:r\ge 0\}$ is a martingale.  Therefore $Y$ has a cadlag version which we still denote by $Y$ and still satisfies $Y_0=\y$.  The fact that $Y_r=0$ for all $r\ge R$ is a standard result for cadlag non-negative supermartingales
and it implies the second formula given for $R$.  Let $\lambda\to\infty$ and then $r\to\infty$ in \eqref{Yemlt} to see that $R<\infty$ a.s. (recall \eqref{YLTform}). Finally a simple calculation shows that $u^{\lambda}(t)=U^{\lambda,t}(0)$ satisfies \eqref{csbpode} with $p=3/2$ and $c_0=\sqrt 6/3$, thus identifying $Y$ as the appropriate SCSBP (recall \eqref{exitmarkov}). 
\eop

\paragraph{Proof of Theorem~\ref{onedims}.} Let $R_n=n\wedge\inf\{r\ge 0:Y_r\le1/n\}\uparrow R$ as $n\to\infty$.  Choose $\z>1/2$,
 set $\beta=3\z$, and for $i\in\Z_+$ let $x_{i,n}=in^{-\z}$ and $I_{i,n}=(x_{i,n},x_{i+1,n}]$. If $x\in I_{i,n}$, then Proposition~\ref{SpecialMProp}(c)  gives
\begin{align*}
\Pe_{\delta_0}(x<R_n,L^x\le 2n^{-\beta})&\le \Ee_{\delta_0}\Bigl(1(Y_{x_{i,n}}\ge 1/n)\Pe_{Y_{x_{i,n}}\delta_{x_{i,n}}}(L^x\le 2n^{-\beta})\Bigr)\\
&\le e^2\Ee_{\delta_0}\Bigl(1(Y_{x_{i,n}}\ge1/n)\Ee_{Y_{x_{i,n}}\delta_{x_{i,n}}}(\exp(-n^{\beta}L^x))\Bigr)\\
&\le e^2\Ee_{\delta_0}\Bigl(1(Y_{x_{i,n}}\ge1/n)\exp(-V^{n^\beta}(x-x_{i,n})/n)\Bigr)\text{ (by Lemma~\ref{LTVL})}\\
&\leq e^2\exp(-n^{2\z-1}V^{n^{\beta-3\z}}((x-x_{i,n})/n^{-\z}))\text{ (by the scaling relation \eqref{Vlscaling})}\\
&\le e^2\exp(-n^{2\z-1}V^1(1)),
\end{align*}
where in the last line we have used the fact that $V^1$ is decreasing (Lemma~\ref{Vlambdaprop}), and $\beta=3\z$. 
Let $k\in\N$ satisfy $k\ge 2\beta$ (an additional condition requiring $k$ large will appear below). A union bound now implies that if
\[\Lambda_n=\{L^x\le 2n^{-\beta}\text{ for some } x\in[0,R_n)\cap\{j n^{-\z-k}:j\in \Z_+\}\},\]
then
\begin{equation}\label{gridbound}
\Pe_{\delta_0}(\Lambda_n)
\le e^2(n^{\z+k+1}+1)\exp(-n^{2\z-1}V^1(1))\equiv p_n.
\end{equation}
Fix $\veps_0\in(0,1)$ and let $\rho=\rho_{\veps_0,1/2}$ be as in Theorem~\ref{LTmodulus} (with $\gamma=1/2$).
Assume $n\ge N(\veps_0)$ so that $n^{-\z}<\veps_0/2$ and $\omega$ is chosen in $\Lambda_n^c\cap\{\rho\ge n^{-\z-k}\}$. If $x\in[\veps_0,R_n\wedge\veps_0^{-1})$, and $j\in \Z_+$ is chosen so that 
$x\in[jn^{-\z-k},(j+1)n^{-\z-k})$,
then 
\[L^x\ge L^{j n^{-\z-k}}-n^{-(\z+k)/2}>2n^{-\beta}-n^{-\beta}=n^{-\beta}.\]
Therefore by \eqref{gridbound} and Theorem~\ref{LTmodulus} we have
\[\Pe_{\delta_0}(L^x\le n^{-\beta}\text{ for some }x\in[\veps_0,R_n\wedge\veps_0^{-1}))\le p_n+\Pe_{\delta_0}(\rho<n^{-\z-k})\le p_n+Cn^{-\kappa(\z+k)}.\]
Now fix $k$ as above and sufficiently large so that the right-hand side of the above is summable over $n$, and hence
by Borel-Cantelli,
\[\text{w.p.}1 \text{ for }n \text{ large enough }\inf_{x\in[\veps_0,R_n\wedge\veps_0^{-1})}L^x>n^{-\beta}.\]
    As this holds for all $\veps_0\in(0,1)$ and $\lim_{x\to0}L^x=L^0>0$ $\Pe_{\delta_0}$-a.s. (by the continuity of $L$ for $d=1$ and \eqref{Vinfty}),
it follows that
\begin{equation}\label{Lpositiv}
\text{w.p.}1 \ L^x>0\text{ for all }x\in[0,R).
\end{equation}
We next show
\begin{equation}\label{Lzero}
\text{w.p.}1 \ L^x=0\text{ for all }x\ge R.
\end{equation}
If $0<r<q$, then by Proposition~\ref{SpecialMProp}(c)
\[\Pe_{\delta_0}(Y_r=0, L^q>0)=\Ee_{\delta_0}(1(Y_r=0)\Pe_{Y_r\delta_r}(L^q>0))=0.\]
So we may fix $\omega$ outside a $\Pe_{\delta_0}$-null set so that the last statement of Proposition~\ref{exitcsbp} holds, $x\to L^x$ is continuous on $(0,\infty)$, and for all rational $0<r<q$, $Y_r=0$ implies $L^q=0$.  If $q>R$ is rational, then by choosing  a rational $r\in(R,q)$ we have $Y_r =0$ and so $L^q=0$.  The continuity of $L^\cdot$ on $(0,\infty)$ and $R>0$ now implies \eqref{Lzero}.  Combining \eqref{Lpositiv} and 
\eqref{Lzero} we get $\{x\ge 0:L^x>0\}=[0,R)$ a.s., and thus $\rb=R$.  By symmetry we also have a r.v. $-\infty<\lb<0$ a.s. so that 
$\{x\le 0:L^x>0\}=(\lb,0]$ a.s. Therefore we have $\{x:L^x>0\}=(\lb,\rb)$ a.s. and the rest is immediate.\eop

\section{Proof of Theorem~\ref{lefttailthm}(a,b)} \label{seclefttail}
\setcounter{equation}{0}
\setcounter{theorem}{0}
 It follows from \eqref{LTLT} and \eqref{Vinfty} that 
\begin{equation}\label{LTpos}
\Ee_{\delta_0}(e^{-\lambda L^x}1(L^x>0))=\exp(-V^\lambda(x))-\exp(-V^\infty(x))\ \forall x\neq0, \lambda>0.
\end{equation}
So by a  Tauberian theorem, asymptotics for $\Pe_{\delta_0}(0<L^x\le a)$ as $a\downarrow 0$ would follow from good asymptotics on $d^\lambda(x)= V^\infty(x)-V^\lambda(x)\ge0$ ($x\neq0$) as $\lambda\to\infty$.
 Using the Feynmann-Kac formula and arguing as in the derivation of 
\eqref{DFK} we obtain
\begin{equation}\label{DFK2}
d^\lambda(x)=d^\lambda(R)E_x\Bigl(1(\tau_R<\infty)\exp\Bigl(-\int_0^{\tau_R}\Bigl(\frac{V^\infty+V^\lambda}{2}\Bigr)(B_s)\,ds\Bigr)\Bigr)\text{ for }|x|\ge R.
\end{equation}
In the next lemmas we prepare necessary tools for establishing bounds on $d^\lambda$, at least for large $\lambda$. 
We start with some a lower bounds on $V^\lambda$; by \eqref{DFK2} these will help bound $d^\lambda$ from above.
\begin{lemma}\label{lemvlrates}
For any $\eta>0$, there is a $\lambda_1(\eta)>0$, so that for any $R>0$, 
\begin{equation}\label{vllbnd1}V^\lambda(x)\ge\frac{2(4-d)-\eta}{|x|^2}\ \ \forall |x|\ge R\text{ and }\lambda\ge\frac{\lambda_1(\eta)}{R^{4-d}},\end{equation}
and if $\lambda_0(\eta)=\lambda_1(\eta)^{1/(4-d)}$, then for any $\lambda>0$,
\begin{equation}\label{vllbnd2}V^\lambda(x)\ge\frac{2(4-d)-\eta}{|x|^2}\ \ \text{if }|x|\ge r_\lambda\equiv\frac{\lambda_0(\eta)}{\lambda^{1/(4-d)}}.\end{equation}
\end{lemma}
\paragraph{Proof.} By \eqref{vinftylimit} we may chose $\lambda_1(\eta)>0$ so that $V^\lambda(1)\ge 2(4-d)-\eta\equiv c_\eta$, if $\lambda\ge \lambda_1(\eta)$.  If $|x|\ge(\lambda_1(\eta)/\lambda)^{1/(4-d)}$, then by the scaling relation \eqref{Vlscaling} and the monotonicity of $V^\lambda$ in $\lambda$, we have 
\[V^\lambda(x)=|x|^{-2}V^{\lambda|x|^{4-d}}(1)\ge |x|^{-2}V^{\lambda_1(\eta)}(1)\ge c_\eta|x|^{-2}.\]
This gives \eqref{vllbnd2}.  The first bound is then immediate because $|x|\ge R$ and $\lambda\ge \frac{\lambda_1(\eta)}{R^{4-d}}$ imply $|x|\lambda^{1/(4-d)}\ge \lambda_0(\eta)$.
\eop

\medskip

The next two lemmas are important for bounding the expectation in~\eqref{DFK2}. Recall the notation introduced before Proposition~\ref{yorthm}.
\begin{lemma}\label{expbound1} Assume $0<\sqrt{2\gamma}\le \nu$. Then for all $r\ge 1$,
\[E_r^{(2+2\nu)}\Bigl(\exp\Bigl(\int_0^{\tau_1}\frac{\gamma}{\rho_s^2}\,ds\Bigr)\Bigr|\tau_1<\infty\Bigr)=r^{\nu-\sqrt{\nu^2-2\gamma}}.\]
\end{lemma}
\paragraph{Proof.} Let $\mu=\sqrt{\nu^2-2\gamma}$, and assume without loss of generality that $r>1$. By monotone convergence the above expectation is 
\[\lim_{t\to\infty} E_r^{(2+2\nu)}\Bigl(\Phi_{t\wedge\tau_1}\rho_{t\wedge\tau_1}^{-\nu+\mu}\Bigr)/P_r^{(2+2\nu)}(\tau_1<\infty),\]
where $\Phi_u=1(\tau_1\le u)\exp\Bigl(\int_0^u\frac{\gamma}{\rho_s^2}\,ds\Bigr)$ (note that $\rho_{t\wedge\tau_1}^{-\nu+\mu}1(\tau_1\le t)=1(\tau_1\le t)$). 

By Proposition~\ref{yorthm} and the above, we conclude that the expectation we are finding equals
\begin{align*}
\lim_{t\to\infty}&r^{\mu-\nu}E_r^{(2+2\mu)}\Bigl(1(\tau_1\le t)\exp\Bigl(\int_0^{t\wedge\tau_1}\frac{\gamma}{\rho_s^2}+\frac{-\gamma}{\rho_s^2}\,ds\Bigr)\Bigr)r^{2\nu}\\
&=r^{\mu+\nu}\lim_{t\to\infty}P_r^{(2+2\mu)}(\tau_1\le t)\\
&=r^{\nu-\mu},
\end{align*}
the last by~\eqref{besselexit}. 
\eop

\begin{lemma} \label{expbound2} 
Assume $0<\sqrt{2\gamma}\le \nu$ and $q>2$. Then
\[\sup_{r\ge 1}E_r^{(2+2\nu)}\Bigl(\exp\Bigl(\int_0^{\tau_1}\frac{\gamma}{\rho_s^q}\,ds\Bigr)\Bigr|\tau_1<\infty\Bigr)\le C_{\ref{expbound2}}(q,\nu)<\infty.\]
\end{lemma}
\paragraph{Proof.}
Consider $r>1$ (without loss of generality), and choose $N\in \Z_+$, so that $1\le r2^{-N}\le 2$.  Define \[\xii(r)=E_r^{(2+2\nu)}\Bigl(\exp\Bigl(\int_0^{\tau_1}\frac{\gamma}{\rho_s^q}\,ds\Bigr)1(\tau_1<\infty)\Bigr)\]
and
\[\beta(r)=E_r^{(2+2\nu)}\Bigl(\exp\Bigl(\int_0^{\tau_{r/2}}\frac{\gamma}{\rho_s^q}\,ds\Bigr)1(\tau_{r/2}<\infty)\Bigr).\]
Apply the strong Markov property of $\rho$ at $\tau_{r/2}$ to see that for $r>2$,
\begin{equation}\label{iterate}
\xii(r)=\beta(r)\xii(r/2).
\end{equation}
By scaling $\hat\rho_t=\frac{2}{r}\rho_{tr^2/4}$ is a $2+2\nu$-dimensional Bessel process starting at $2$ under $P_r^{(2+2\nu)}$, and $\tau_{r/2}=\frac{r^2}{4}\hat\tau_1$, where $\hat\tau_1=\inf\{t\ge0:\hat\rho_t\le1\}$.  Therefore for $r>2$,
\begin{align}\label{besscale}
\beta(r)&=E_2^{(2+2\nu)}\Bigl(\exp\Bigl(\Bigl(\frac{r}{2}\Bigr)^{2-q}\int_0^{\tau_1}\frac{\gamma}{\rho_u^q}\,du\Bigr)\Bigr|\tau_1<\infty\Bigr)P_2^{(2+2\nu)}(\tau_1<\infty)\\
\nonumber&\le E_2^{(2+2\nu)}\Bigl(\exp\Bigr(\int_0^{\tau_1}\frac{\gamma}{\rho_u^2}\,du\Bigr)\Bigr|\tau_1<\infty\Bigr)^{(r/2)^{2-q}}2^{-2\nu},
\end{align}
where in the last we use \eqref{besselexit} and $(r/2)^{2-q}<1$.  Apply the previous lemma to conclude that 
\[\beta(r)\le 2^{(\nu-\sqrt{(\nu^2-2\gamma)}\,)(r/2)^{2-q}}2^{-2\nu}\le 2^{\nu(r/2)^{2-q}}2^{-2\nu}.\]
Insert the above into \eqref{iterate} to get
\[r^{2\nu}\xii(r)\le(r/2)^{2\nu}\xii(r/2)2^{\nu(r/2)^{2-q}}\ \text{ for }r>2.\]
Iterate this $N$ times (recall $1<r2^{-N}\le 2$) to conclude
\begin{align}\label{iterbound}r^{2\nu}\xii(r)&\le (r2^{-N})^{2\nu}\xii(r2^{-N})2^{\nu\sum_{j=1}^N(r2^{-j})^{2-q}}\\
\nonumber&\le E^{(2+2\nu)}_{r2^{-N}}\Bigl(\exp\Bigl(\int_0^{\tau_1}\frac{\gamma}{\rho_s^q}\,ds\Bigr)\Bigr|\tau_1<\infty\Bigr)\,2^{\nu r^{2-q}2^{q-2}2^{N(q-2)}/(2^{q-2}-1)}\ \text{(by \eqref{besselexit})}\\
\nonumber&\le 2^{\nu-\sqrt{\nu^2-2\gamma}}2^{\nu 2^{q-2}/(2^{q-2}-1)}\\
\nonumber&\le 2^{\nu}2^{\nu 2^{q-2}/(2^{q-2}-1)}=C_{\ref{expbound2}}(q,\nu).
\end{align}
In the next to last line we have used Lemma~\ref{expbound1} and $r2^{-N}\ge 1$. 
Note that (by \eqref{besselexit}) the left-hand side of \eqref{iterbound} is 
\[E_r^{(2+2\nu)}\Bigl(\exp\Bigl(\int_0^{\tau_1}\frac{\gamma}{\rho_s^q}\,ds\Bigr)\Bigr|\tau_1<\infty\Bigr),\]
and so we are done.
\eop

We next use \eqref{DFK2} and the above lemmas to get lower and upper bounds on $d^\lambda(x)$ in terms of 
$\frac{R^p}{|x|^p}d^\lambda(R)$.

\begin{lemma}\label{vlconvrate}
There are universal positive constants $\lambda_{\ref{vlconvrate}}$ and $C_{\ref{vlconvrate}}$ so that if $R>0$, then\\
\noi(a) 
\begin{equation}\label{Vllb}
0<\frac{R^p}{|x|^p}{(V^\infty(R)-V^\lambda(R))}\le{V^\infty(x)-V^\lambda(x)}\quad \forall |x|\ge R, \lambda>0.
\end{equation}
\noi(b) \begin{equation}\label{Vlub}V^\infty(x)-V^\lambda(x)\le C_{\ref{vlconvrate}}\frac{R^p}{|x|^p}(V^\infty(R)-V^\lambda(R))\quad \forall |x|\ge R, \lambda\ge \lambda_{\ref{vlconvrate}}/R^{4-d}.
\end{equation}
\end{lemma}
\paragraph{Proof.} Introduce
\begin{equation}\label{munudef}
\mu=\begin{cases} -1/2&\text{ if }d=1\\
0&\text{ if }d=2,\qquad \text{ and }\nu=\sqrt{\mu^2+4(4-d)}.\\
1/2&\text{ if }d=3
\end{cases}
\end{equation}
(a) Note that if $d^\lambda(x)=0$ for some $x\neq0$, then by \eqref{LTpos} and the a.s. finiteness of $L^x$, we get $L^x=0\ \Pe_{\delta_0}-$a.s. This contradicts \eqref{Vinfty} and
so the first (strict) inequality in \eqref{Vllb} is established. Use $V^\lambda\le V^\infty$ in \eqref{DFK2} to see that for $|x|\ge R$, 
\begin{align*}
d^\lambda(x)&\ge d^\lambda(R)E_x\Bigl(1(\tau_R<\infty)\exp\Bigl(-\int_0^{\tau_R}\frac{2(4-d)}{|B_s|^2}\,ds\Bigr)\Bigr)\\
&=d^\lambda(R)\lim_{t\to\infty}E_{|x|}^{(2+2\mu)}\Bigl(1(\tau_R\le\tau_R\wedge t)\exp\Bigl(-\int_0^{\tau_R\wedge t}\frac{2(4-d)}{|B_s|^2}\,ds\Bigr)\Bigr)\\
&=d^\lambda(R)\lim_{t\to\infty}|x|^{\nu-\mu}E_{|x|}^{(2+2\nu)}(1(\tau_R\le\tau_R\wedge t)(\rho_{t\wedge\tau_R})^{-\nu+\mu})\quad (\text{by Proposition~\ref{yorthm}})\\
&=d^\lambda(R) (R/|x|)^{\mu-\nu}P_{|x|}^{(2+2\nu)}(\tau_R<\infty)\\
&=d^\lambda(R)(R/|x|)^p.
\end{align*}
the last by $p=\mu+\nu$ and \eqref{besselexit}. This gives (a).\\
(b) Fix $R>0$. If $\eta\in(0,2(4-d))$, let $\nu_\eta=\sqrt{\mu^2+4(4-d)-\eta}$, and 
$p_\eta=\nu_\eta+\mu\to p>2$ as $\eta\downarrow0$. Fix $\eta>0$ so that $p_\eta>2$.  Now set $\lambda_{\ref{vlconvrate}}=\lambda_1(\eta)$ ($\lambda_1$ as in Lemma~\ref{lemvlrates}) and assume $\lambda\ge\lambda_{\ref{vlconvrate}}/R^{4-d}$.  Then by \eqref{vllbnd1} and \eqref{DFK2} we have for $|x|\ge R$,
\begin{align*}
d^\lambda(x)&\le d^\lambda(R)E_x\Bigl(1(\tau_R<\infty)\exp\Bigl(-\int_0^{\tau_R}\frac{2(4-d)-(\eta/2)}{|B_s|^2}\,ds\Bigr)\Bigr)\\
&=d^\lambda(R)\lim_{t\to\infty}E_{|x|}^{(2+2\mu)}\Bigl(1(\tau_R\le t\wedge\tau_R)\exp\Bigl(-\int_0^{\tau_R\wedge t}\frac{2(4-d)-(\eta/2)}{|B_s|^2}\,ds\Bigr)\Bigr)\\
&=d^\lambda(R)\lim_{t\to\infty}|x|^{\nu_\eta-\mu}E_{|x|}^{(2+2\nu_\eta)}(1(\tau_R\le t)R^{-\nu_\eta+\mu})\quad(\text{by Proposition~\ref{yorthm}})\\
&=d^\lambda(R)(R/|x|)^{p_\eta}.
\end{align*}
So if $\xii(R)=d^\lambda(R)R^{p_\eta}/2$, then we have shown
\begin{equation}\label{templb}
\Bigl(\frac{V^\lambda+V^\infty}{2}\Bigr)(x)\ge V^\infty(x)-\frac{\xii(R)}{|x|^{p_\eta}}\quad\text{for }|x|\ge R.
\end{equation}
Use this in \eqref{DFK2} to see that for $|x|\ge R$,
\[d^\lambda(x)\le d^\lambda(R)E_x\Bigl(1(\tau_R<\infty)\exp\Bigl(\int_0^{\tau_R}\frac{\xii(R)}{|B_s|^{p_\eta}}\,ds\Bigr)\exp\Bigl(-\int_0^{\tau_R}\frac{2(4-d)}{|B_s|^2}\,ds\Bigr)\Bigr).\]
Now use Fatou's lemma and then Proposition~\ref{yorthm} as in (a) to conclude that for $|x|\ge R$,
\begin{align*}
d^\lambda(x)&\le d^\lambda(R)\liminf_{t\to\infty}E_{|x|}^{(2+2\mu)}\Bigl(1(\tau_R\le\tau_R\wedge t)\exp\Bigl(\int_0^{\tau_R\wedge t}\frac{\xii(R)}{\rho_s^{p_\eta}}\,ds\Bigr)\exp\Bigl(-\int_0^{\tau_R\wedge t}\frac{2(4-d)}{\rho_s^2}\,ds\Bigr)\Bigr)\\
&=d^\lambda(R)\liminf_{t\to\infty}E_{|x|}^{(2+2\nu)}\Bigl(1(\tau_R\le\tau_R\wedge t)\exp\Bigl(\int_0^{\tau_R\wedge t}\frac{\xii(R)}{\rho_s^{p_\eta}}\,ds\Bigr)\Bigr)(R/|x|)^{\mu-\nu}\\
&=d^\lambda(R)E_{|x|}^{(2+2\nu)}\Bigl(\exp\Bigl(\int_0^{\tau_R}\frac{\xii(R)}{\rho_s^{p_\eta}}\,ds\Bigr)\Bigr|\tau_R<\infty\Bigr)(R/|x|)^p,
\end{align*}
the last by monotone convergence and \eqref{besselexit}. A scaling argument, as in \eqref{besscale}, shows that the above equals
\[d^\lambda(R)(R/|x|)^pE_{|x|/R}^{(2+2\nu)}\Bigl(\exp\Bigl(\int_0^{\tau_1}\frac{\xii(R)R^{2-p_\eta}}{\rho_u^{p_\eta}}\,du\Bigr)\Bigr|\tau_1<\infty\Bigr).\]
To apply Lemma~\ref{expbound2} we note that
\[2\gamma\equiv2\xii(R)R^{2-p_\eta}\le V^\infty(R)R^{p_\eta}R^{2-p_\eta}=2(4-d)<\nu^2,\]
and so Lemma~\ref{expbound2} and the above bound show that 
\[d^\lambda(x)\le d^\lambda(R)(R/|x|)^pC_{\ref{expbound2}}(p_\eta,\nu).\]
This gives the required result since the last constant depends only on $d$.
\eop

Finally, we are ready to 
 establish the rate of convergence of $V^\lambda(x)$ to $V^\infty(x)$ in \eqref{vinftylimit}.
Recall from the Introduction that $\alpha=(p-2)/(4-d)$.

\begin{proposition}\label{Vlambdarate}
\noi (a) There is a constant $C_{\ref{Vlambdarate}}$, depending only on $d$, so that 
$$V^\infty(x)-V^\lambda(x)\le C_{\ref{Vlambdarate}}|x|^{-p}\lambda^{-\alpha}\quad \forall x\neq 0,\ \lambda>0.$$
\noi(b) For all $\veps>0$ there is a $c_ {\ref{Vlambdarate}}(\veps)>0$ so that 
$$V^\infty(x)-V^\lambda(x)\ge c_{\ref{Vlambdarate}}(\veps)|x|^{-p}\lambda^{-\alpha}\quad \forall |x|\ge \veps\lambda^{-1/(4-d)},\ \lambda>0.$$
\noi(c) There is a $\underline c_ {\ref{Vlambdarate}}>0$ so that 
$$V^\infty(x)-V^\lambda(x)\ge \underline c_{\ref{Vlambdarate}}|x|^{-p}\lambda^{-\alpha}\quad \forall \lambda\ge |x|^{-(4-d)},\ |x|>0.$$
\end{proposition}
\paragraph{Proof.} By the scaling property of $V^\lambda$ (recall\eqref{Vlscaling}) we have 
\begin{equation}\label{diffscale}
V^\infty(x)-V^\lambda(x)=r^{-2}(V^\infty(x/r)-V^{\lambda r^{4-d}}(x/r))\ \forall x\neq 0,\ r>0.
\end{equation}
(a) Let $\lambda>0$ and set $r=(\lambda_{\ref{vlconvrate}}/\lambda)^{1/(4-d)}$, so that $\lambda r^{4-d}=\lambda_{\ref{vlconvrate}}$.  If $|x/r|\ge 1$, then applying Lemma~\ref{vlconvrate}(b) with $R=1$ to the right-hand side of \eqref{diffscale} we get,
\begin{equation}\label{Vdiffub1}
V^\infty(x)-V^\lambda(x)\le C_{\ref{vlconvrate}}r^{p-2}|x|^{-p}(V^\infty(1)-V^{\lambda_{\ref{vlconvrate}}}(1))\equiv C_1|x|^{-p}\lambda^{-\alpha}.
\end{equation}
If $|x|<r$, then
\begin{align}\label{Vdiffub2}
|x|^{-p}\lambda^{-\alpha}&\ge r^{-(p-2)}|x|^{-2}\lambda^{-\alpha}\\
\nonumber&=\lambda_{\ref{vlconvrate}}^{-(p-2)/4-d)}|x|^{-2}\\
\nonumber&\ge c(d)(V^\infty(x)-V^\lambda(x)),
\end{align}
for some $c(d)>0$.  Clearly \eqref{Vdiffub1} and \eqref{Vdiffub2} imply (a).\\
(c) Fix $x\not=0$, and assume $\lambda\ge |x|^{-(4-d)}$.
Set $r=\lambda^{-1/(4-d)}$ (thus $|x/r|\ge 1$) and then apply Lemma~\ref{vlconvrate}(a) to the right-hand side of \eqref{diffscale} with $R=1$ to see that 
\begin{align*}
V^\infty(x)-V^\lambda(x)&\ge  r^{p-2}|x|^{-p}(V^\infty(1)-V^1(1))\\
&\equiv \underline c_{\ref{Vlambdarate}}|x|^{-p}\lambda^{-\alpha}.
\end{align*}
(b) now follows by applying (c) to $x/\veps$ and using the scaling relation \eqref{diffscale} with $r=\veps$.
\eop

Now we are ready to complete the 
\paragraph{Proof of Theorem~\ref{lefttailthm}(a,b)}\mbox{}

\noindent
(a)  Apply \eqref{LTpos} and then Proposition~\ref{Vlambdarate}(a) to  conclude that for all $x\neq 0$,
\begin{equation}\label{LTupperboundt}\Ee_{\delta_0}(e^{-\lambda L^x}1(L^x>0))\le V^\infty(x)-V^\lambda(x)\le C_{\ref{Vlambdarate}}|x|^{-p}\lambda^{-\alpha}\ \forall\lambda>0.
\end{equation}
(a) now follows from Markov's inequality (take $\lambda=a^{-1}$ in the above) with $C_{\ref{lefttailthm}}=eC_{\ref{Vlambdarate}}$.

\noindent(b) For the lower bound, note that \eqref{LTpos} and Proposition~\ref{Vlambdarate}(b)
with $\veps=1$ imply that for $|x|\ge\veps_0$,
\begin{align}\label{LTlowerboundt}
\Ee_{\delta_0}(e^{-\lambda L^x}1(L^x>0))&\ge e^{-V^\infty(x)}(V^\infty(x)-V^\lambda(x))\\
\nonumber&\ge e^{-V^\infty(\veps_0)}c_{\ref{Vlambdarate}}(1)|x|^{-p}\lambda^{-\alpha}\quad \forall\lambda\ge\veps_0^{-(4-d)}.
\end{align}
Apply a Tauberian theorem of de Haan and Stadtm\"uller \cite{bib:dhs85} (see Lemma~4.7(b) of \cite{bib:mmp16} for an appropriate quantitative version) to see that \eqref{LTlowerboundt} and \eqref{LTupperboundt} imply there is a $c_{\ref{lefttailthm}}(\veps_0)>0$ so that (b) holds.
\eop

\section{On Non-polar Sets for $F$ and Preliminaries for the Lower Bound on the Dimension} \label{lbh}
\setcounter{equation}{0}
\setcounter{theorem}{0}

To show that the lower bound on dim$(F)$ in Theorem~\ref{dimthm} holds with positive probability we
employ the methodology that was used for the proof of Theorem 5.5 in~\cite{bib:mmp16}.
The next proposition is crucial for carrying out that program: it plays here the role of Proposition~5.1 
in~\cite{bib:mmp16}.

\begin{proposition}\label{L2upperbound} For all $\veps_0>0$ there is a $C_{\ref{L2upperbound}}$ so that for all $\lambda\ge 1$ and all $|x_i|\ge \veps_0$,
\[\lambda^{2+2\alpha}\Ee_{\delta_0}\Bigl(\prod_{i=1}^2L^{x_i}e^{-\lambda L^{x_i}}\Bigr)\le C_{\ref{L2upperbound}}(\veps_0)(1+|x_1-x_2|^{2-p}).\]
\end{proposition}
The proof is involved and hence is deferred to Section~\ref{sec:9}.

\medskip

As we are focusing on lower bound results for $F$, in light of the fact that $F=\{L,R\}$ for $d=1$ (see Section~\ref{sec1d}), we will assume $d=2$ or $3$.

If $\beta>0$ and $g_\beta(r)=r^{-\beta}$ for a finite measure  $\mu$ on $\Re^d$ and Borel subset $A$ of $\Re^d$, let 
\[\langle\mu\rangle_{g_\beta}=\int\int g_\beta(|x-y|)d\mu(x)d\mu(y),\]
and
\[I(g_\beta)(A)=\inf\{\langle\mu\rangle_{g_\beta}:\mu \text{ a probability supported by }A\}.\]
The $g_\beta$-capacity of $A$ is $C(g_\beta)=1/I(g_\beta)(A)$ (see, e.g., Section 3 of \cite{bib:haw79}).

Set 
$$\beta=p-2=\begin{cases}2\sqrt 2-2&\text{ if }d=2\\
\frac{\sqrt{17}-3}{2}&\text{ if }d=3,\end{cases}$$
and note $\beta\in(1/2,1)$.  

\begin{theorem}\label{polarsetslb} Assume $d=2$ or $3$. For every $\veps_0\in(0,1)$ there is a $c_{\ref{polarsetslb}}(\veps_0)>0$ such that for any Borel subset, $A$, of $\{x\in\Re^d:\veps_0\le |x|\le \veps_0^{-1}\}$, 
\[\Pe_{\delta_0}(F\cap A\neq\emptyset)\ge c_{\ref{polarsetslb}}C(g_{\beta})(A).\]
In particular for any Borel subset $A$ of $\Re^d$, $C(g_{\beta})(A)>0$ implies that $\Pe_{\delta_0}(F\cap A\neq\emptyset)>0$. 
\end{theorem}
\paragraph{Proof.} This follows from Theorem~\ref{lefttailthm}(a,b) and Proposition~\ref{L2upperbound} by standard arguments exactly as in the proof of Theorem 5.2 and Corollary 5.3 of \cite{bib:mmp16}.\eop

Let $Z_t=(Z^1_t,\dots,Z^d_t)$, where $(Z^i,i\le d)$ are i.i.d. $\Re$-valued symmetric L\'evy processes with L\'evy measure $\nu(dx)=|x|^{-1-\beta}((\log(1/|x|))\vee 1)^2dx$, starting at zero. This means that if 
\begin{align}
\nonumber\psi(\theta)&=\int_{-\infty}^\infty [1-e^{i\theta x}-i\theta x1(|x|\le 1)]\,\nu(dx)\\
&=2|\theta|^\beta\int_0^\infty(1-\cos u)u^{-1-\beta}((\log(|\theta|/u))\vee 1)^2\,du,
\end{align}
then
$E(e^{i\theta_jZ^j_t})=\exp(-t\psi(\theta_j))$ for all $j\le d\text{ and }\theta_j\in\Re.$ If $\log^+r=\log(r\vee 1)$ for $r\in\Re$, then
a straightforward calculation shows:

\begin{lemma}\label{psibounds}
There are constants $0<c_{\ref{psibounds}}\le C_{\ref{psibounds}}$ so that for all real $\theta$,
\[c_{\ref{psibounds}}|\theta|^\beta[1+(\log^+(|\theta|))^2]\le \psi(\theta)\le C_{\ref{psibounds}}|\theta|^\beta[1+(\log^+(|\theta|))^2].\]
\end{lemma}

\begin{lemma}\label{levypot}
(a) If $A$ is a Borel subset of $\Re^d$ such that $\text{dim}(A)<d-\beta$, then $A$ is polar for $Z$, that is, $P(Z_t\in A\text{ for some }t>0)=0$.

\noindent(b) $C(g_\beta)(\{Z_s:1/2\le s\le 1\})>0$ a.s., and if $B$ is any non-empty open set, then $$P(C(g_\beta)(\{Z_s:1/2\le s\le 1\}\cap B)>0)>0.$$
\end{lemma}
\paragraph{Proof.}
(a) For $\theta=(\theta_1,\dots,\theta_d)\in\Re^d$ and $\beta'\in(\beta,2)$, let $\psi_\ell(\theta)=\sum_{j=1}^d\psi(\theta_j)$ and $\psi_{\beta'}(\theta)=|\theta|^{\beta'}$.  Then $e^{-t\psi_\ell(\theta)}$ and $e^{-t\psi_{\beta'}(\theta)}$ are the characteristic functions of $Z_t$ and $Y_t$, respectively, where $Y$ is a symmetric stable process of index $\beta'$.  It follows from Lemma~\ref{psibounds} that for some $C>0$, and all $\theta\in\Re^d$,
\[0<1+\psi_\ell(\theta)\le C(1+\psi_{\beta'}(\theta)),\]
and so 
\begin{equation}\label{psicomp}\text{Re}\Bigl(\frac{1}{1+\psi_\ell(\theta)}\Bigr)\ge C^{-1}\text{Re}\Bigl(\frac{1}{1+\psi_{\beta'}(\theta)}\Bigr).\end{equation}
Lemma~\ref{psibounds} shows that $\int_{-\infty}^\infty e^{-t\psi(\theta)}d\theta<\infty$ and so by Fourier inversion $Z^j(t)$ has bounded density
\begin{equation}\label{Zjdensity}
f_t(z)=(2\pi)^{-1}\int_{-\infty}^\infty e^{-i\theta z}e^{-t\psi(\theta)}d\theta=(2\pi)^{-1}\int_{-\infty}^\infty \cos(\theta z)e^{-t\psi(\theta)}d\theta.
\end{equation}
Therefore $Z_t$ has a bounded density and hence $Z$ also has a resolvent density. Corollary~15 of Chapter II of \cite{bib:bert96} and \eqref{psicomp} imply that any set which is polar for $Y$ is polar for $Z$.  Here we are using the fact the existence of a 
resolvent density for $Y$ and $Z$ implies that essentially polar sets are polar (by \cite{bib:haw79}) and so we can replace essentially polar with polar in the aforementioned Corollary~15.  If $\text{dim}(A)<d-\beta$ then $\text{dim}(A)<d-\beta'$,\ for some $\beta'\in(\beta,2)$, and by the potential theory for $Y$ (see, e.g., Lemma~10 of \cite{bib:tak64}) $A$ is polar for $Y$, and hence also polar for $Z$. 
\medskip

\noindent(b) Define a probability supported by $\{Z_s:s\in[1/2,1]\}$ by $\mu(A)=\int_{1/2}^11_A(Z_s)\,ds$. Then
\begin{align}
\label{energyboundbeta}\langle\mu\rangle_{g_\beta}&=2E\Bigl(\int_{1/2}^1\int_{s_1}^1|Z_{s_2}-Z_{s_1}|^{-\beta}\,ds_2\,ds_1\Bigr)\\
\nonumber&\le 2E\Bigl(\int_0^1|Z_s|^{-\beta}ds\Bigr)\le 2E\Bigl(\int_0^1|Z^1_s|^{-\beta}\,ds\Bigr).
\end{align}
Using \eqref{Zjdensity} and monotone convergence we have
\begin{align*}E\Bigl(\int_0^1|Z^1_s|^{-\beta}ds\Bigr)&\le 1+E\Bigl(\int_0^1|Z^1_s|^{-\beta}1(|Z^1_s|\le 1)\,ds\Bigr)\\
&=1+\lim_{\veps\to0+}\int_\veps^1
\int_{-1}^1|z|^{-\beta}(2\pi)^{-1}\int_{-\infty}^\infty \cos(\theta z)e^{-s\psi(\theta)}\,d\theta dzds.\end{align*}
For each fixed $\veps>0$ the integrand with $\cos(\theta z)$ replaced by $1$ is integrable, and so we can use the above with \eqref{energyboundbeta} and Fubini
to conclude that 
\begin{align*} 
\langle\mu\rangle_{g_\beta}&\le 2+\lim_{\veps\to 0} \int_{-\infty}^\infty\int_{-1}^1|z|^{-\beta}\cos(\theta z)dz\frac{e^{-\veps\psi(\theta)}-e^{-\psi(\theta)}}{\psi(\theta)}d\theta\\
&\le 2+\int_{-\infty}^\infty\Bigl|\int_{-|\theta|}^{|\theta|}|y|^{-\beta}\cos y\,dy\Bigr||\theta|^{\beta-1}\frac{1-e^{-\psi(\theta)}}{\psi(\theta)}d\theta\\
&\le 2+C\int_{-\infty}^\infty |\theta|^{\beta-1}\frac{1-e^{-\psi(\theta)}}{\psi(\theta)}d\theta,
\end{align*}
where in the last line an elementary calculus argument is used to bound $\Bigl|\int_{-|\theta|}^{|\theta|}|y|^{-\beta}\cos y\,dy\Bigr|$ uniformly in $\theta$.  Lemma~\ref{psibounds} shows the contribution to the above integral  from $|\theta|>e$ is finite and the trivial bound $1-e^{-\psi(\theta)}\le \psi(\theta)$ shows the contribution from $|\theta|\le e$ is also finite. 
This completes the proof of the first statement in (b).  For the second statement it suffices to show
\begin{equation}\label{contain}
P(\{Z_s:1/2\le s\le 1\}\subset B\})>0,
\end{equation}
where $B=B(x_0,r)$ is an open ball.
The Markov property of $Z$ shows that the above probability is at least
$$P(Z_{1/2}\in B(x_,r/2))P(\sup_{s\le 1/2}|Z_s|\le r/2).$$
It is now easy to show each of the above factors is positive, for example by writing $Z$ as 
the sum of two independent L\'evy processes, one with jumps bigger than $\epsilon$ and one with
jumps smaller than $\epsilon$ for sufficiently small $\epsilon$. This completes the proof. 
\eop

We are ready to prove that the lower bound on dim$(F)$ in Theorem~\ref{dimthm} holds with positive probability. 
\begin{proposition}\label{dimensionlowerbound} If $B$ is a non-empty open set, then
$\Pe_{\delta_0}(\textnormal{dim}(F\cap B)\ge d+2-p)>0$.
\end{proposition}
\paragraph{Proof.} We work on the product space under $\Pe_{\delta_0}\times P$ where $P$ is the probability under which $Z$ is the $d$-dimensional L\'evy process considered above.  Let $R(\omega_1,\omega_2)=R(\omega_2)=\{Z_s:s\in[1/2,1]\}$. 
By Theorem~\ref{polarsetslb} and Lemma~\ref{levypot}(b),
\[(\Pe_{\delta_0}\times P)(F(\omega_1)\cap (B\cap R(\omega_2))\neq\emptyset))>0.\]
This implies that 
\[\Pe_{\delta_0}(\{\omega_1:P(\{\omega_2:(F(\omega_1)\cap B)\cap R(\omega_2)\neq\emptyset\}>0)>0.\]
By Lemma~\ref{levypot}(a) this implies that $\Pe_{\delta_0}(\text{dim}(F\cap B)\ge d-\beta)>0$.  As $d-\beta=d+2-p$, the proof is complete.
\eop

It will be useful when extending the above lower bound to an a.s. statement (in Section~\ref{lbh2})
 to have a version of the above bound for the canonical measure of the Brownian snake, $\N_0$. 
 
\begin{corollary}\label{canonicallowerbound} If $B$ is a non-empty open set, then
$\N_0(\textnormal{dim}(F\cap B)\ge d+2-p)\equiv p_{\ref{canonicallowerbound}}(B)>0$.
\end{corollary}
\paragraph{Proof.} By reducing $B$ we may assume $B$ is an open ball such that $0\notin \bar B$.
If $w$ is a continuous $\Re^d$-valued path, let $\tau(w)=\inf\{t\ge 0:w_t\in \bar B\}\le\infty$. 
We work in the standard setup from Section~\ref{prel} with $X_0=\delta_0$. 
Let $\zeta^i_s$ be the lifetime of $W_{i,s}(\cdot)$, let $\hat W_{i,s}=W_{i,s}(\zeta^i_s)$ be the position of the tip of the snake at time $s$.  We abuse notation slightly and set $\tau(W_i)=\tau(\hat W_i)$. If $I_B=\{i\in I:\tau(W_i)<\infty\}$ then $|I_B|$ is a Poisson r.v. with mean $\N_0(\tau<\infty)<\infty$, the last since $0\notin\bar B$. Therefore, given $I_B$, $\{W_i:i\in I_B\}$ are iid with law $\N_0(\cdot|\tau<\infty)$.  Clearly $\tau(W_i)=\infty$ implies that $L^x(W_i)=0$ for $x\in\overline B$, and so by \eqref{exdecompL},
\[\text{for }x \in \bar B,\ L^x=\sum_{i\in I_B} L^x(W_i)\equiv\sum_{i\in I_B}L^{x,i}.\]
If $F_i=\partial\{x:L^{x,i}>0\}$, then an elementary argument (left for the reader) shows that 
\[F\cap B\subset\cup_{i\in I_B}F_i\cap B.\]
It now follows from Proposition~\ref{dimensionlowerbound} that 
\begin{align*}
0<p_B&=\Ee_{\delta_0}(\Pe_{\delta_0}(\text{dim}(\cup_{i\in I_B}(F_i\cap B))\ge d+2-p|I_B))\\
&\le \Ee_{\delta_0}(\sum_{i\in I_B}\Pe_{\delta_0}(\text{dim}(F_i\cap B)\ge d+2-p|I_B))\\
&=\Ee_{\delta_0}(|I_B|)\N_0(\text{dim}(F\cap B)\ge d+2-p|\tau<\infty)\\
&=[\Ee_{\delta_0}(|I_B|)/\N_0(\tau<\infty)]\N_0(\text{dim}(F\cap B)\ge d+2-p)\\
&=\N_0(\text{dim}(F\cap B)\ge d+2-p).
\end{align*}
The result follows.
\eop

\section{The  Lower Bound on the Dimension} \label{lbh2}
\setcounter{equation}{0}
\setcounter{theorem}{0}
Recall that $\cR=\overline{\{x:L^x>0\}}$ and
 $\text{conv}(X_0)$ denotes the closed convex hull of Supp$(X_0)$.
This section is devoted to the proof of the following theorem.
\begin{theorem}\label{genlowerboundl} Assume $\textnormal{conv}(X_0)\neq \Re^d$.  Then $\Pe_{X_0}$-a.s.
\[ \textnormal{conv}(X_0)^c\cap \cR\neq\emptyset\text{ implies }\textnormal{dim}(\textnormal{conv}(X_0)^c\cap F)\ge d+2-p.\]
\end{theorem}

The key step will be obtaining a lower bound on dim$(F)$ under the canonical measure in the next result.

\begin{proposition}
\label{prop:can1}
Let $ H$ be an open half-space such that $0\in \partial H$, and for $r>0$,  let  ${H_r}$ be 
${H}$ translated by $r$ (perpendicular to $\partial {H}$) so that it is increasing in $r$. Under $\N_0$ there is a cadlag version of the total exit measure mass, $X_{H_r}(1)$, and for this version, 
\begin{equation*}
\N_0(\exists r>0:\;\textnormal{dim}(F\cap \overline{H}_r^c)<2+d-p, X_{{H}_r}(1)>0)=0.
\end{equation*}
\end{proposition}
\paragraph{Proof}

Fix $\epsilon>0$. By translation and  rotation, and considering $r>\veps$ in the Proposition,  we may assume that the process is given under the excursion measure of the snake,
$\N_{-\epsilon}\equiv  \N_{(-\epsilon,0,\ldots,0)}$,
 ${H}={H_0}=\{x\in\IR^d:x_1< 0\}$,
 and for $r\ge 0$, $H_r=\{x:x_1<r\}$. For $r\ge 0$ define $Z_r=X_{H_r}$ and $Y_r=Z_{H_r}(1)$. 
Hence the objective is to show  that there is a cadlag version of $Y$ satisfying
\begin{equation}
\label{eq:25_3}
\N_{-\epsilon}(\exists r>0:\;\textnormal{dim}(F\cap \overline{H_r}^c)<2+d-p, Y_r>0)=0.
\end{equation}

We define the snake $W_t=(W^j_{t},j\le d)$ under $ \N_{-\epsilon}$  and let $\hat W_t=W_t(\zeta_t)\equiv(\hat W^j_{t},j\le d)$ be the tip of the 
snake $W_t$.  We also denote  $\tau_r(W)=\inf\{t\ge 0:\hat W^1_{t}\ge r\}$.
Following Section 2.4 of \cite{bib:leg95}, for $r,s,u\ge 0$ define
\[S_r(W_u)=\inf\{t\le \zeta_u:W^1_{u}(t)\ge r\},\]
\[\eta^{r}_s=\inf\{t:\int_0^t 1(\zeta_u\le S_r(W_u))\,du>s\},\]
and
\[\cE_r=\sigma(W_{\eta_s^{r}},s\ge 0)\vee \{ \N_{-\epsilon}-\text{null sets}\}.\]
Note that the inequality in the definition of $\eta_s^{r}$ is attained for $t$ sufficiently large, so that $\eta^{r}_s<\infty$ for all $s\ge 0$. One can check that $\cE_r$ is non-decreasing in $r$ (this will also follow from \eqref{timecheta5} below with $T\equiv r' \le r$). Intuitively $\cE_r$ is the $\sigma$-field
generated by the excursions of $W$ in $\overline{ H_r}$.
We set $\cE^{+}_r=\cap_{r'>r}\cE_{r'}$.  It follows from Proposition 2.3 of \cite{bib:leg95} that $Z_r$ is $\cE_r$-adapted.

An elementary argument shows that for $r\ge 0$ there is a measurable map $\psi:C(({\overline{H_r}})^c,\IR)\to \{0,1\}$ so that 
\begin{equation} \label{dimmeas_canon}
1(\text{dim}(F\cap({\overline{H_r}})^c)<2+d-p)=\psi((L^x,x\in ({\overline{H_r}})^c)).
\end{equation}
For this, note that this easily reduces to considering canonical {\it compact} subsets of $F\cap({\overline{H_r}})^c$ for which the optimal open coverings by open balls reduces to finite open covers by ``rational balls".  One also needs to note that $F\cap({\overline{H_r}})^c$ is the boundary of $\{x\in ({\overline{H_r}})^c:L^x>0\}$ {\it in the space } $({\overline{H_r}})^c$. The details are routine. The above allows us to apply the special Markov property (Proposition~\ref{SpecialMProp}(b))  to conclude that for $r\geq s \ge  0$, $ \N_{-\epsilon}$-a.e.,
\begin{align}
\nonumber M^s_r(\omega)&\equiv\N_{-\epsilon}(\text{dim}(F\cap({\overline{H_s}})^c)\ge 2+d-p|\cE_r)(\omega)\\
\nonumber&\ge\N_{-\epsilon} (\text{dim}(F\cap({\overline{H_r}})^c)\ge 2+d-p|\cE_r)(\omega)\\
\label{mrlb_canon}&=\Pe_{Z_r(\omega)}(\text{dim}(F\cap({\overline{H_r}})^c)\ge 2+d-p).
\end{align}

At this point we need to extend Proposition~\ref{dimensionlowerbound} to a more general class of initial measures. The following is a simple consequence of Corollary~\ref{canonicallowerbound} and scaling.

\begin{lemma}\label{iclowerboundondim} There is a universal constant $p_{\ref{iclowerboundondim}}>0$ so that if $S(X'_0)\subset\{x:x_1=0\}$ and $X'_0(1)>0$, then
\begin{equation*}
\Pe_{X'_0}\Bigl(\textnormal{dim}\Bigl( F\cap\Bigl\{x:x_1>\sqrt{X'_0(1)}\Bigr\}\Bigr)\ge d+2-p\Bigr)\ge p_{\ref{iclowerboundondim}}.
\end{equation*}
\end{lemma}
\paragraph{Proof.} Let $\delta=X'_0(1)>0$, and set $X_0^{(\delta)}(A)=\delta^{-1}X'_0(\sqrt{\delta} A)$ and $L^{(\delta),x}=L^{x/\sqrt{\delta}}\delta^{2-(d/2)}$.
By scaling we have 
\[\Pe_{X'_0}((L^x,x_1>0)\in\cdot)=\Pe_{X_0^{(\delta)}}((L^{(\delta),x},x_1>0)\in\cdot),\]
where we only restricted to $x_1>0$ to ensure continuity of $L^x$ (by \cite{bib:sug89}). So under $\Pe_{X^{(\delta)}_0}$,  if $$F^{(\delta)}=\partial\{x:L^{(\delta),x}>0\},$$ then
\begin{align*}
F^{(\delta)}\cap\{x:x_1>\sqrt{\delta}\}&=(\partial\{x:L^{x/\sqrt\delta}>0,x_1>0\})\cap\{x:x_1>\sqrt\delta\}\\
&=\sqrt\delta[(\partial\{x:L^x>0,x_1>0\})\cap\{x:x_1>1\}].
\end{align*}
Therefore
\begin{align}
\nonumber\Pe_{X'_0}(\text{dim}(F\cap\{x:x_1>\sqrt\delta\})\ge d+2-p)&=\Pe_{X_0^{(\delta)}}(\text{dim}(F^{(\delta)}\cap \{x:x_1>\sqrt\delta\})\ge d+2-p)\\
\label{rescaledim}&=\Pe_{X_0^{(\delta)}}(\text{dim}(F\cap\{x:x_1>1\})\ge d+2-p).
\end{align}
We continue to use the notation of our standard setting, that is,  $\{W_i:i\in I\}$ is a Poisson point process on $C([0,\infty),\cW)$  with intensity $\N_{X_0}(dW)$. 
We let $W_{i,t}=(W^{j}_{i,t},j\le d)$ and $\hat W_{i,t}=W_{i,t}(\zeta^i_t)\equiv(\hat W_{i,t}^j,j\le d)$ be the tip of the $i$th
snake $W_i$.  We also set $\tau_r(W_i)=\inf\{t\ge 0:\hat W_{i,t}^1\ge r\}$.
Let $N_1$ be the number of $i\in I$ such that $\tau_1(W_i)<\infty$, so that under $\Pe_{X_0^{(\delta)}}$, $N_1$ is Poisson with mean $\N_{X_0^{(\delta)}}(\tau_1<\infty)=\N_0(\tau_1<\infty)<\infty$.  This last equality holds because $X_0^{(\delta)}(1)=1$, and for any $x$ such that $x_1=0$, $\N_{x}(\tau_1<\infty)=\N_0(\tau_1<\infty)$ by translation invariance. Similar reasoning shows that 
\begin{align}\label{Transinv}
\N_{X_0^{(\delta)}}&(\text{dim}(F\cap\{x:x_1>1\})\ge d+2-p|\tau_1(W)<\infty)\\
\nonumber&=\N_0(\text{dim}(F\cap\{x:x_1>1\})\ge d+2-p|\tau_1(W)<\infty).
\end{align}
We may assume that given $N_1$, $\{W_i:\tau_1(W_i)<\infty\}$ are iid with law $\N_{X_0^{(\delta)}}(W\in\cdot|\tau_1(W_i)<\infty)$. Using this, we see from \eqref{rescaledim} and \eqref{Transinv} that
\begin{align*}
\Pe_{X'_0}&(\text{dim}(F\cap\{x:x_1>\sqrt\delta\})\ge d+2-p)\\
&\ge\Pe_{X_0^{(\delta)}}(N_1=1)\N_{X_0^{(\delta)}}(\text{dim}(F\cap\{x:x_1>1\})\ge d+2-p|\tau_1(W)<\infty)\\
&=\exp(-\N_0(\tau_1<\infty))\N_0(\tau_1<\infty)\N_0(\text{dim}(F\cap\{x:x_1>1\})\ge d+2-p)/\N_0(\tau_1<\infty)\\
&=\exp(-\N_0(\tau_1<\infty))p_{\ref{canonicallowerbound}}(\{x:x_1>1\}),
\end{align*}
the last by Corollary~\ref{canonicallowerbound}. This completes the proof.
\eop

We return to the derivation of~\eqref{eq:25_3}. 
First note that  the projection of the $d$-dimensional snake under $\N_{-\epsilon}(\cdot)$ onto the first coordinate is a $1$-dimensional snake. Let $\N^1_{-\epsilon}$ denote the excursion measure of the corresponding one-dimensional snake. Then, using Proposition~\ref{SpecialMProp}(a), we  get that under 
$\N^1_{-\epsilon}$, and conditional on $\cE_0$, $\{W^i, i\in I\}$ is a Poisson point process with intensity $Y_0\N^1_0(\cdot)$.  However for our standard setup, under $\Pe_{\y\delta_0}$, $\{W_i, i\in I\}$
is also a Poisson point process with intensity $\y\N^1_0(\cdot)$. Thus $\N^1_{-\epsilon}(Y\in \cdot|\cE_0)(\omega)=
\Pe_{Y_0(\omega)\delta_0}(Y\in \cdot)$, and $Y$ constructed here has the same finite-dimensional distributions as 
$Y$ in Proposition~\ref{exitcsbp}, where $\y=Y_0(\omega)$. 
 In particular we may work with the cadlag version of $(Y_r,r\ge 0)$ obtained there and define an $(\cE_r^{+})$-stopping time by $T_0=\inf\{r\ge 0:Y_r=0\}<\infty$ a.s. (the finiteness by Proposition~\ref{exitcsbp}). Fix $s_0\geq 0$. By \eqref{mrlb_canon}, Lemma~\ref{iclowerboundondim} and translation invariance we have
\begin{equation}\label{martlb}M^{s_0}_r\ge p_{\ref{iclowerboundondim}}\text{ on }\{r<T_0\}\;\N_{-\epsilon}-\text{a.e. for all }r\ge s_0.\end{equation}
Now let $G_{s_0}=\{\text{dim}(F\cap({\overline{H_{s_0}}})^c)\ge 2+d-p\}$, and work under the probability $Q_{s_0}(\cdot)= \N_{-\epsilon}(\cdot\ |T_0>s_0)$. Then
the definition of conditional expectation and $\{T_0>s_0\}{\stackrel{{\text{a.s.}}}{=}}\{Y_{s_0}>0\}\in\cE_{s_0}$ imply that $M^{s_0}_r=Q_{s_0}(G_{s_0}|\cE_r)$ $Q_{s_0}$-a.s. for $r\ge s_0$, and so is an $(\cE_r)$-martingale
under this law, where we are adding the slightly larger class of $Q_{s_0}$-null sets to our filtration.
Now take limits in \eqref{martlb} from above along rationals to see that for each $r\ge s_0$,
\begin{equation*}M^{s_0}_{r+}=Q_{s_0}(G_{s_0}|\cE^+_r)\ge p_{\ref{iclowerboundondim}}\quad Q_{s_0}-\text{a.s. on }\{r<T_0\}.
\end{equation*}
The $(\cE^+_r)$-martingale $M^{s_0}_{r+}, r\ge s_0,$ has a cadlag version, and we will abuse notation and let $M^{s_0}_r$ denote
this cadlag version. Then we can conclude from the above that
\begin{equation}\label{M+lb} M^{s_0}_r\ge p_{\ref{iclowerboundondim}}\text{ for all }r\in [s_0,T_0) \  Q_{s_0}-\text{a.s.}
\end{equation}
The SCSBP $Y$ (recall Proposition~\ref{exitcsbp}) only has non-negative jumps (see e.g. Theorem~1 in~\cite{bib:clb09}) and therefore 
$T_n=\inf\{r\ge 0 :Y_r\le 1/n\}$ increase to $T_0$ and are strictly smaller than $T_0$ a.s. on $\{T_0>0\}$.
By   \eqref{M+lb} we have
\begin{equation}\label{MTn}
Q_{s_0}(G_{s_0}|\cE^+_{T_n\vee s_0})\ge p_{\ref{iclowerboundondim}}\ Q_{s_0}-\text{a.s.}
\end{equation}
We have  $(T_n\vee s_0)<T_0, $
 $Q_{s_0}$-a.s., and so by (17.9)(ii) and (17.10) of Chapter VI of \cite{bib:rw94} we have  
\begin{equation}\label{T0-}
Q_{s_0}(G_{s_0}|\cE^+_{T_0-})=M_{T_0-}\ge p_{\ref{iclowerboundondim}}\ \ Q_{s_0}-\text{a.s., and }\vee_{n}\cE^+_{T_n\vee s_0}=\cE^+_{T_0-}.
\end{equation}
 We claim
that 
\begin{equation}\label{ET0-meas}
G_{s_0}\in\vee_{n}\cE^+_{T_n}.
\end{equation}
Assuming this, we see from \eqref{T0-} that $1_{G_{s_0}}\ge  p_{\ref{iclowerboundondim}}>0$ $Q_{s_0}$-a.s., and 
hence $1_{G_{s_0}}=1$ $Q_{s_0}$-a.s., which implies 
\begin{equation*}
\N_{-\epsilon}(\textnormal{dim}(F\cap \overline{H_{s_0}}^c)<2+d-p,Y_{s_0}>0)=0.
\end{equation*}
for all $s_0\geq 0$. From this we immediately get that
$$\N_{-\epsilon}(\exists\ \text{rational}\ s\ge 0: \textnormal{dim}(F\cap \overline{H_s}^c)<2+d-p,Y_s>0)=0. $$
Use the non-decreasing property of $s\mapsto 1(\textnormal{dim}(F\cap \overline{H_s}^c)<2+d-p)$ and right-continuity of $Y$ to complete the proof of the proposition. 

\eop

It remains to prove \eqref{ET0-meas}. Intuitively this is obvious as $T_0$ will be the rightmost level reached by the first coordinate of the snake (see\eqref{T0max} below for the important one-sided bound for $W$), and so observing the snake $W$ for the first coordinate to the left of $T_0$ means we see all of the $W$ and hence know $1_{G_{s_0}}$.  The reader happy with this explanation should skip the rest of this proof on a first reading.  
 We claim it suffices to show that for all $s\ge 0$, 
\begin{equation}
\label{ETOi-meas}
W_s\text{ is }\vee_n\cE^{+}_{T_n}-\text{measurable}.
\end{equation}
Indeed,  this condition implies that $W$ is $\vee_n\cE^{+}_{T_n}-$-measurable, and recalling \eqref{dimmeas_canon} with $r=s_0$ and that $L^x=L^x(W)$, we conclude that $G_{s_0}$ is $\vee_n\cE^+_{T_n}$-measurable, thus proving \eqref{ET0-meas}, as required.  


\begin{lemma}\label{etaA}
(a) If $A_t^r=\int_0^t 1(\zeta_u\le S_r(W_u))\,du$, then $\N_{-\epsilon}$-a.e. for all $t\ge 0$, 
\[A^r_t=\int_0^t1(\sup_{v\le \zeta_u}W^1_u(v)<r)du\text{ for all }r\ge 0,\]
and $r\to A^r_t$ is left continuous on $(0,\infty)$.

(b) $\lim_{r'\uparrow r}\eta_s^{r'}=\eta^r_s$ for all $r>0$ and $s\ge 0$ $\N_{-\epsilon}$-a.e.

(c) If $T$ is an $(\cE_r^+)$-stopping time then $W_{\eta^T_s}$ is $\cE^+_T$-measurable.
\end{lemma}

Let us assume the above result and first finish the proof of \eqref{ETOi-meas}.  Recall we are working under $\N_{-\epsilon}$.
By Lemma~\ref{etaA}(b,c) we may conclude that 
\[W_{\eta_s^{T_0}}\text{ is }\vee_n\cE^+_{T_n}-\text{measurable}.\]
So to prove \eqref{ETOi-meas} it clearly suffices to show $W_{\eta_s^{T_0}}=W_s$ $\N_{-\epsilon}$-a.e., and this clearly would follow from $A_t^{T_0}=t$ for all $t\ge 0$ $\N_{-\epsilon}$-a.e. or equivalently (by Lemma~\ref{etaA}(a)),
\begin{equation}\label{notatmax}
\int_0^\infty 1\Bigl(\sup_{v\le \zeta_u}W^1_u(v)\ge T_0\Bigr)du=0\ \N_{-\epsilon}-\text{a.e.}
\end{equation}
A simple application of the Special Markov Property under $\N_{-\epsilon}$ (see 
Proposition~\ref{SpecialMProp}(a) with $G=H_r$) shows that on $\{T_0\le r\}\in\cE^+_r$, $\sup_{u,v\le\zeta_u}W^1_u(v)\le r\quad \N_{-\epsilon}-\text{a.e.}$  Take limits over rational $r\downarrow T_0$ to see that 
\begin{equation}\label{T0max}
\sup_{u,v\le\zeta_u}W^1_u(v)\le T_0\quad \N_{-\epsilon}-\text{a.e.}
\end{equation}
Therefore by the above and \eqref{notatmax} it suffices to establish
\begin{equation}\label{notatmax2}
\int_0^\infty1\Bigl(\sup_{v\le \zeta_u}W^1_u(v)=T_0\Bigr)du=0\ \N_{-\epsilon}-\text{a.e.}
\end{equation}
If $y$ is a one-dimensional continuous path defined on $[0,\zeta]$ for some $\zeta<\infty$ or on $[0,\infty)$, set $M(y)=\sup_t y(t)$, where the $\sup$ is over the domain of $y$.  Let $(H_t,t\ge 0)$ be the $1$-dimensional historical Brownian motion constructed from $W^1$ so that (by p. 64 of \cite{bib:leg99})
for any non-negative measurable function $\phi$ on the space of $1$-dimensional paths,
\begin{equation}\label{histproc}
\int_0^\infty \int\phi (y)H_s(dy)ds=\int_0^\infty \phi(W^1_u)du.
\end{equation}
Therefore 
\eqref{notatmax2} is equivalent to 
\begin{equation}\label{notatmax3}
\int_0^\infty\int1(M(y)=T_0)H_s(dy)ds=0\quad \N_{-\epsilon}-\text{a.e.}
\end{equation}
By \eqref{T0max} and \eqref{histproc}, if $X$ is the one-dimensional super-Brownian motion associated with $H$ (or equivalently $W^1$), this clearly would follow from
\begin{equation}\label{histpalm1}
\int\Bigl[ \int1\Bigl(\int_0^\infty  X_u((M(y),\infty))du=0\Bigr)H_s(dy)\Bigr]d\N_{-\epsilon}=0\text{ for each }s> 0.
\end{equation} 
Let $(B(t),t\ge 0)$ denote a one-dimensional Brownian motion starting at $-\epsilon$ under $P_{-\epsilon}$ and let $M_s=\sup_{t\le s}B(t)$.  
The Palm measure formula for $H_s$ under its canonical measure (see Proposition 4.1.5 of \cite{bib:dp91}
with $\beta=1$ and $\gamma=1/2$), shows that the left-hand side of \eqref{histpalm1} equals
\begin{align}\nonumber\lim_{\lambda\to\infty}&\int\Bigl[\int \exp\Bigl(-\lambda\int_0^\infty X_u((M(y),\infty))du\Bigr)H_s(dy)\Bigr]d\N_{-\epsilon}\\
\nonumber&=\lim_{\lambda\to\infty}E_{-\epsilon}\Bigl(\exp\Bigl(-\int_0^s\int\Bigl[1-\exp\Bigl(-\lambda\int_0^\infty X_u((M_s,\infty))du\Bigr)\Bigr]d\N_{B(t)}(X)dt\Bigr)\Bigr)\\
\nonumber&=E_{-\epsilon}\Bigl(\exp\Bigl(-\int_0^s\int1\Bigl(\int_0^\infty X_u((M_s,\infty))du>0\Bigr)d\N_{B(t)}(X)dt\Bigr)\Bigr)\\
\nonumber&\le E_{-\epsilon}\Bigl(\exp\Bigl(-\int_0^s\int 1(L^{M_s}>0)d\N_{B(t)}dt\Bigr)\Bigr)\\
\label{Bint}&=E_{-\epsilon}\Bigl(\exp\Bigl(-\int_0^s\frac{6}{(B_t-M_s)^2}dt\Bigr)\Bigr),
\end{align}
where in the last line we have used \eqref{vcanonical} with $\lambda\to\infty$ and $d=1$. An easy 
application of L\'evy's modulus of continuity shows that $\int_0^s(B_t-M_s)^{-2}dt=\infty$ $P_{-\epsilon}$-a.s. and so \eqref{Bint} is zero. This proves \eqref{histpalm1} and so completes the proof of Proposition~\ref{prop:can1} once we establish Lemma~\ref{etaA}.

\medskip
\noindent{\bf Proof of Lemma~\ref{etaA}.} (a) Note that $\{\zeta_u=S_r(W_u)\}\subset\{\hat W^1_u=r\}$, so that 
\begin{equation}\label{edgemass}
\int_0^\infty1(\zeta_u=S_r(W_u))du\le \int_0^\infty 1(\hat W^1_u=r)du=\int_0^\infty X_u^{(1)}(\{r\})du,
\end{equation}
where $X^{(1)}$ is the one-dimensional super-Brownian motion associated with the projection of the snake, $W^1$.  The existence of local time shows that the right-hand side of \eqref{edgemass} is zero for all $r\ge 0$ $\N_{-\epsilon}$-a.e. This shows that $\N_{-\epsilon}$-a.e. for all $r,t\ge 0$,
\[A^r_t=\int_0^t 1(\zeta_u<S_r(W_u))du=\int_0^t1(\sup_{v\le \zeta_u}W^1_u(v)<r)du,\]
where the second equality is elementary. This gives the first part of (a) and the second part is then immediate by Monotone Convergence. 

\noindent(b)  Work outside the null set so that (a) holds.  Let $s\ge 0$. Then $\eta^r_s$ is decreasing in $r$ since $A^r_t$ is increasing in $r$ and therefore $\eta^{r-}_s\ge \eta^r_s$.  If $r>0$ and $\delta>0$, then $A^r_{\eta_s^r+\delta}>s$ by the definition of $\eta_s^r$. By the left continuity of $r'\to A^{r'}_{\eta_s^r+\delta}$, there is an $\veps>0$ such that $A^{r'}_{\eta_s^r+\delta}>s$ for $r'>r-\veps$, which implies that $\eta_s^{r'}\le \eta_s^r+\delta$ for $r'>r-\veps$. This proves that $\eta_s^{r-}\le \eta_s^r$ and (b) is proved.

\noindent(c) Fix $r>0$ and let $W'_s=W_{\eta_s^r}$ which is continuous in $s$ as noted on p. 401 of \cite{bib:leg95}.  We claim that 
\begin{equation} \label{Wtimech}
\text{on }\{T\le r\},\ W_{\eta^T_s(W)}=W'_{\eta^T_s(W')},
\end{equation}
where we are denoting the snake dependence of $\eta_s^T$ explicitly.  Assuming this claim and noting that $W'_{\eta^T_s(W')}$
is a measurable function of $(W',T)$, we see that for a measurable set $A\subset\mathcal{W}$,
\[\{T\le r\}\cap\{W_{\eta_s^T}\in A\}=\{T\le r\}\cap\{W'_{\eta^T_s(W')}\in A\}\in\cE^+_r,\] 
as required (the result then follows immediately for $r=0$). 

Turning to \eqref{Wtimech}, we may use (a) to see that on $\{T\le r\}$,
\begin{align}
\nonumber\eta^T_s&=\inf\{t:\int_0^t 1(\sup_{v\le \zeta_u} W^1_u(v)<T)1(\sup_{v\le \zeta_u} W^1_u(v)<r)du>s\}\\
\label{Arint}&=\inf\{t:\int_0^t 1(\sup_{v\le \zeta_u} W^1_u(v)<T)dA^r(u)>s\}.
\end{align}
We claim that for any Borel $\psi:[0,\infty)\to[0,\infty)$, 
\begin{equation}\label{timecheta}
\int_0^t\psi(u)dA^r(u)=\int_0^{A^r_t}\psi(\eta_u^r)\,du\quad\forall t\ge 0.
\end{equation}
As usual it suffices to consider $\psi(u)=1(u\le a)$ for $a\ge 0$.  Note that 
\[A_a^r>q\Rightarrow \eta^r_q\le a\Rightarrow A^r_a\ge q,\]
and so,
\[\int_0^{A^r_t}1(\eta_q^r\le a)dq=\int_0^{A^r_t}1(q\le A^r_a)dq=A^r_t\wedge A^r_a=\int_0^t1(u\le a)dA^r(u),\]
thus proving \eqref{timecheta}. 
Now use \eqref{timecheta} in \eqref{Arint} to see that on $\{T\le r\}$,
\begin{align}\nonumber
\eta_s^T=&\inf\{t:\int_0^{A_t^r(W)}1(\sup_{v\le \zeta(\eta_q^r)}W^1_{\eta_q^r}(v)<T)dq>s\}\\
\label{timecheta2}=&\inf\{t:A^T_{A^r_t(W)}(W')>s\}.
\end{align}
Call $t$ a point of increase of $A^r_t(W)$, and write $t\in I^r(W)$, iff for all $\veps>0$, $A^r_{t+\veps}(W)>A^r_t(W)$. Clearly it suffices to take the infimum in \eqref{timecheta2} over $t\in I^r(W)$ and for any such $t$, $t=\eta^r_{A^r_t}$.  Therefore \eqref{timecheta2} implies that on $\{T\le r\}$,
\begin{align}
\nonumber\eta^T_s&=\inf\{\eta^r_{A^r_t(W)}:t\in I^r(W), A^T_{A^r_t(W)}(W')>s\}\\
\label{timecheta3}&=\inf\{\eta_v^r(W):A^T_v(W')>s\}.
\end{align}
In the last line we use the right continuity of $v\to\eta^r_v$ for all $r$ (by an elementary argument), and the fact that $\{A^r_t:t\in I^r(W)\}\supset[0,\infty)\setminus C$, where $C$ is the countable set of values of $A^r$ corresponding to values of $A^r$ at ``flat spots" of $A^r$.  To ease eyestrain we will write $\eta^r(W)(s)$ for $\eta^r_s(W)$. The right continuity of $s\to \eta^r(W)(s)$ and \eqref{timecheta3} imply that
\begin{equation}\label{timecheta4}
\eta^T(W)(s)=\eta^r(W)(\inf\{v:A^T_v(W')>s\})=\eta^r(W)(\eta^T(W')(s))\quad\text{on }\{T\le r\}.
\end{equation}
Therefore on $\{T\le r\}$,
\begin{equation}\label{timecheta5}W_{\eta^T(W)(s)}=W_{\eta^r(W)(\eta^T(W')(s))}=W'_{\eta^T(W')(s)},
\end{equation}
and so \eqref{Wtimech} is proved, thus completing the proof of Lemma~\ref{etaA}.\eop


\begin{corollary}
\label{cor:25_1}
Assume $\textnormal{conv}(X_0)\neq \Re^d$ and
 $H$ is an open half-space such that $\text{conv}(X_0)\subset \overline H$.  For $r>0$,  let  ${H_r}$ be 
${H}$ translated by $r$ (perpendicular to $\partial H$) so that it is increasing in $r$.
Under $\N_{X_0}$ there is a cadlag version of the total exit measure mass, $X_{H_r}(1)$, and for this version, 
\begin{equation*}
\N_{X_0}(\exists r>0:\;\textnormal{dim}(F\cap \overline{H_r}^c)<2+d-p, X_{H_r}(1)>0)=0.
\end{equation*}
\end{corollary}
This is immediate by Proposition~\ref{prop:can1} and the definition of $\N_{X_0}$.

\bigskip

\paragraph{Proof of Theorem~\ref{genlowerboundl} }

By the Hahn-Banach Theorem and separability of $\Re^d$ there is a countable collection of open half-spaces $\{{H_j}:j\in J\}$ such that 
\[\text{conv}(X_0)=\cap_{j\in J} \overline{H_j}.\]
Our condition that $\text{conv}(X_0)^c\cap \cR$ is non-empty easily implies there is an $x_0\notin \text{conv}(X_0)$ so that $L^{x_0}>0$.  We may choose $j_0\in J$ so that $x_0$ lies in the open half-space $\overline{H_{j_0}}^c$.  We may choose a natural number $n$ so that if we translate $\overline{H_{j_0}}$ by $1/n$ (perpendicular to $\partial \overline{H_{j_0}}$), then (denoting the translated open half-space still by ${H_{j_0}}$) $d(x_0,\overline{H_{j_0}})\ge 1/n$ and
\begin{equation}\label{Hcond} d(\text{conv}(X_0), H_{j_0}^c)\ge 1/n.
\end{equation}
By Proposition~\ref{SpecialMProp}(c) w.p. 1
$L^{x_0}>0$  implies that $X_{H_{j_0}}(1)>0$ and so by further increasing $n$ we may
assume  $X_{H_{j_0}}(1)>1/n$.

As there are countably many choices of $({H_{j_0}},1/n)$ it suffices to fix such a pair as in \eqref{Hcond} 
and show that 
\begin{equation*}
\Pe_{X_0}(\text{dim}(F\cap \overline{H_{j_0}}^c)<2+d-p, X_{H_{j_0}}(1)\ge 1/n)=0.
\end{equation*}
By translation and rotation we may assume that ${H_{j_0}}=\{x\in\IR^d:x_1< 0\}$ and so $S(X_0)\subset (-\infty,-1/n]\times\Re^{d-1}$. For $r\ge 0$ define $H_r=\{x:x_1<r\}$, $Z_r=X_{H_r}$ and $Y_r=Z_{H_r}(1)$. 
Hence our objective is to show that for $S(X_0)$ as above,
\begin{equation}\label{dimhalfsp}
\Pe_{X_0}(\text{dim}(F\cap({\overline{H_0}})^c)<2+d-p, Y_0\ge 1/n)=0.
\end{equation}

Working in the standard setup under $\Pe_{X_0}$, we let $W_{i,t}=(W^j_{i,t},j\le d)$ and $\widehat W_{i,t}=W_{i,t}(\zeta^i_t)\equiv(\widehat W_{i,t}^j,j\le d)$ be the tip of the $i$th
snake $W_i$.  If $\tau_r(W_i)=\inf\{t\ge 0:\widehat W_{i,t}^1\ge r\}$, let 
\[I_0=\{i\in I:\tau_0(W_i)<\infty\}=\{i
\in I: \sup_{u}\widehat W_{i,u}^1\ge 0\}.\]
Therefore $|I_0|$ is a Poisson r.v. with mean $\N_{X_0}(\tau_0<\infty)<\infty$ (recall that $S(X_0)$ is bounded away from $H_0^c$). Moreover if $\widetilde W_i, i=1,2,\ldots$ are iid with the 
law  $\N_{X_0}(W\in \cdot| \tau_0<\infty)$   and the sequence $(\widetilde W_i, i=1,2,\ldots)$ is   independent of $|I_0|$, then 
\begin{equation}
\label{eq:25_1}
 \sum_{i\in I_0} \delta_{W_i}  \text{ is equal in law to } \sum_{i=1}^{|I_0|}\delta_{\widetilde W_i}. 
\end{equation}

Given $|I_0|$, define
\begin{eqnarray*}
T_{i,0}&=&\sup_{u} \widehat {{\widetilde W}}^1_{i,u}, \ i\ge 1,\\
\tilde T_0 &=& \inf\{s: \exists! \ i=:\tilde i\leq |I_0|\text{ such that } T_{\tilde i,0}>s, T_{j,0}<s, \forall j\leq |I_0|, j\not= \tilde i   \},\\
F_i &=& F(\widetilde W_i).
\end{eqnarray*}
If 
\begin{equation}
T_0=T_0(W)\equiv \sup_{u} \widehat W^1_u\,,
\end{equation} 
then 
for each $i$, $T_{i,0}$ is distributed according to the law $\N_{X_0}(T_0\in \cdot | T_0\geq 0)$. 
If $\N^1_{x_1}$ is the excursion measure of the first component of the snake, then
\begin{eqnarray}
\N_x(T_0 \ge y)=\N^1_{x_1}(y\in\cR)=\frac{6}{(y-x_1)^2},\; \forall y>x_1.
\end{eqnarray}
The first equality uses (8) on p. 69 of \cite{bib:leg99} and our earlier comments in Section~\ref{intro} on the definitions of $\cR$, and the second then follows from Theorem~1.3 of \cite{bib:dip89} with $d=1$.
From the above and by the definition of $T_{i,0}$ we easily derive that the law of $T_{i,0}$ is absolutely continuous with respect to Lebesgue measure. 
This, and the independence of  the $T_{i,0}$'s imply 
\begin{equation}
\label{eq:25_2}
\Pe_{X_0}(T_{i,0}\not= T_{j,0}, \forall i\not=j\in I_0)=1,  
\end{equation}
and so $\tilde T_0=\vee_{i=1,i\neq\tilde i}^{|I_0|}T_{i,0}<\infty$ $\Pe_{X_0}$-a.s. on $\{I_0\neq \emptyset\}$.
By \eqref{eq:25_2} and  the definitions of $\tilde i, \tilde T_0$ we obtain
\begin{align}
\nonumber 
\Pe_{X_0}(\text{dim}(F\cap({\overline{H_0}})^c)\ge 2+d-p)
&
\ge\Pe_{X_0}(\text{dim}(F_{\tilde i} \cap({\overline{H_{\tilde T_0}}})^c)\ge 2+d-p, I_0\neq \emptyset)\\
\nonumber
&\ge \Pe_{X_0}(\{|I_0|>0\}\cap\bigcap_{i\in I_0} \left\{\text{dim}(F_i\cap({\overline{H_r}})^c)\ge 2+d-p, \forall r<T_{i,0}\right\}).
\end{align}
Then we have 
\begin{align}
\nonumber
&\Pe_{X_0}(\{|I_0|>0\}\cap\bigcap_{i\in I_0} \left\{\text{dim}(F_i\cap({\overline{H_r}})^c)\ge 2+d-p, \forall r<T_{i,0}\right\})\\
\nonumber
&= \Ee_{X_0}\left(1(\{|I_0|>0\})\Pi_{i=1}^{|I_0|}\N_{X_0} \left(\text{dim}(F\cap({\overline{H_r}})^c)\ge 2+d-p, \forall r<T_{0}| \tau_0<\infty\right)\right)
\\
\nonumber
&=\Pe_{X_0}(|I_0|>0)\\
\nonumber
&\geq \Pe_{X_0}(Y_0>0). 
\end{align}
where the first equality follows by~\eqref{eq:25_1}, the second equality follows by  Corollary~\ref{cor:25_1}, and the last line is immediate from the definition of the exit measure, $Z_0$, in Ch. V of \cite{bib:leg99}. 
So if $G_0=\{\text{dim}(F\cap({\overline{H_0}})^c)\ge 2+d-p\}$, we see from the above  that
$\Pe_{X_0}(G_0^c)\leq \Pe_{X_0}(Y_0=0)$ and hence
\begin{align}
\nonumber
\Pe_{X_0}(G_0^c, Y_0\geq 1/n)&= \Pe_{X_0}(G_0^c) - \Pe_{X_0}(G_0^c, Y_0< 1/n) \\
\nonumber
&\leq \Pe_{X_0}(Y_0=0) - \Pe_{X_0}(G_0^c, Y_0< 1/n)
\\
\nonumber
&\rightarrow  \Pe_{X_0}(Y_0=0) - \Pe_{X_0}(G_0^c, Y_0=0)\;\;\text{(as $n\rightarrow\infty$)}\\
\nonumber
&\leq\Pe_{X_0}(G_0, Y_0=0)\\
\nonumber
&=0,
\end{align}
where the last equality is obvious since on $\{Y_0=0\}$, we have $F\cap({\overline{H_0}})^c\subset \cR\cap({\overline{H_0}})^c =\emptyset$ $\Pe_{\delta_0}$-a.s. by Proposition~\ref{SpecialMProp}(c). 
Thus we have derived~\eqref{dimhalfsp} and so are done. 

\eop

\section{Remaining Proofs of Main Results}\label{secproofs}
\setcounter{equation}{0}
\setcounter{theorem}{0}

\noindent{\bf Proof of Theorem~\ref{dimthm}.} This is immediate from Theorem~\ref{dimub2} and Theorem~\ref{genlowerboundl},
the latter with $X_0=\delta_0$ so that $\text{conv}(X_0)^c=\Re^d\setminus\{0\}$.
\eop
\medskip

\noindent{\bf Proof of Theorem~\ref{dimN0}.} For $i=1\dots,d$, and $\veps>0$, let $H^i_{\veps}=\{x:x_i>\veps\}$, and $-H^i_{\veps}=\{x:x_i<-\veps\}$. By the special Markov property (see Proposition~\ref{SpecialMProp}(b))  $X_{e{\overline{H^i_\veps}^c}}(1)=0$ implies $\int _{eH_\veps^i}L^xdx=0$ $\N_0$-a.e. for all $i\le d$ and $e=\pm$.
Therefore Proposition~\ref{prop:can1} shows that for $i$ and $e$ as above,
\[\N_0\Bigl(\text{dim}(F)<d+2-p,\ \ \int_{eH_\veps^i}L^xdx>0\Bigr)=0\quad\forall\veps>0.\]
Take the union over $i$ and $e$ to conclude
\[\N_0\Bigl(\text{dim}(F)<2+d-p,\ \int_{\{|x|>\veps\}} L^xdx>0)=0\quad\forall\veps>0,\]
and hence (let $\veps\downarrow 0$)
\begin{equation}\label{N0lowbnd}\N_0(\text{dim}(F)<2+d-p)=\N_0(\text{dim}(F)<2+d-p, \int_0^\infty X_s(1)ds>0)=0.
\end{equation}

Consider next the upper bound on $\text{dim}(F)$. Fix $\veps>0$ and let $L=(L^x,|x|>\veps)$.  Then \eqref{exdecompL} implies that under $\Pe_{\delta_0}$, $L=\sum_{i=1}^{N_\veps}L_i$ (addition is componentwise), where $N_\veps$ is Poisson with mean $\N_0\Bigl(\int_0^\infty X_s(|x|>\veps)ds>0\Bigr)<\infty$ and given $N_\veps$, $(L_i=(L^ x_i,|x|>\veps))_{i\in\N}$ are iid
with law $\N_0\Bigl(L\in\cdot\ \Bigl|\int_0^\infty X_s(|x|>\veps)ds>0\Bigr)$.  Theorem~\ref{dimub2} implies that 
\begin{align*}0&=\Pe_{\delta_0}(N_\veps=1,\text{dim}(F\cap\{|x|>\veps\})>d+2-p)\\
&=\Pe_{\delta_0}(\N_\veps=1)\N_0\Bigl(\text{dim}(F\cap\{|x|>\veps\})>d+2-p\Bigl|\int_0^\infty X_s(|x|>\veps)ds>0\Bigr).
\end{align*}
Therefore we have $\N_0(\text{dim}(F\cap\{|x|>\veps\})>d+2-p, \int_0^\infty X_s(|x|>\veps)ds>0)=0\text{ for all }\veps>0$. Let $\veps\downarrow 0$ to conclude that $\N_0(\text{dim}(F)>d+2-p)=0$.
This and \eqref{N0lowbnd} imply the result.
\eop
\medskip

\noindent{\bf Proof of Theorem~\ref{lefttailthm}.}
Since (a), (b) of Theorem~\ref{lefttailthm} have been proved 
in Section~\ref{seclefttail} we only need consider (c), (d), i.e., work under $\N_0$.  By \eqref{vcanonical} (including the $\lambda=\infty$ case) we have for $x\neq 0$ (without loss of generality),
\begin{align}
\label{eq:15_10_1}
\int e^{-\lambda L^x}1(L^x>0)d\N_0&=\int 1(L^x>0)d\N_0-\int(1-e^{-\lambda L^x})d\N_0\\
\nonumber
&=V^\infty(x)-V^\lambda(x).
\end{align}
Normalize \eqref{eq:15_10_1} to get the Laplace transform of a probability:
 \begin{equation}
\label{eq:15_10_2}
\frac{\int e^{-\lambda L^x}1(L^x>0)d\N_0}{\N_0(L^x>0)}
=\frac{V^\infty(x)-V^\lambda(x)}{V^\infty(x)}\le C_1|x|^{2-p}\lambda^{-\alpha},
\end{equation}
where the last line holds by Proposition~\ref{Vlambdarate}(a). A simple application of Markov's inequality
now gives (for any $a>0$, $x\neq 0$)
\[\N_0(0<L^x\le a)\le eC_1 V^\infty(x)|x|^{2-p}a^\alpha=c_1|x|^{-p}a^\alpha,\]
proving (c). For 
(d), use the equality in \eqref{eq:15_10_2} and then apply Proposition~\ref{Vlambdarate}(c)  to get for all
$\lambda\ge |x|^{-(4-d)}$,
 \begin{align}
\label{eq:15_10_3}
\frac{\int e^{-\lambda L^x}1(L^x>0)d\N_0}{\N_0(L^x>0)}
&\geq  \frac{\underline c_{\ref{Vlambdarate}}}{2(4-d)}|x|^{2-p}\lambda^{-\alpha}
\\
\nonumber
 &=  C_2|x|^{2-p}\lambda^{-\alpha}.
\end{align}
\eqref{eq:15_10_2} and \eqref{eq:15_10_3} allow use to apply a Tauberian theorem (Lemma~4.7(b) of \cite{bib:mmp16}), and after a short calculation conclude  there exists a constant 
$c_2>0$ such that for all $x\neq 0$ and $a\in[0,1]$, 
 \begin{align}
\label{eq:15_10_4}
\frac{\N_0(0<L^x\leq a)}{\N_0(L^x>0)}
 &\ge  c_2|x|^{2-p} \min\{ 1, |x|^{\alpha(4-d)} \}  a^{\alpha}
\\
\nonumber
&= c_2 \min\{ |x|^{2-p}, 1 \}  a^{\alpha}.
\end{align}
Recalling that $\N_0(L^x>0)=V^\infty(x)=2(4-d)|x|^{-2}$, we obtain the result.

\eop
\medskip

\noindent{\bf Proof of Theorem~\ref{dimX0}.}  (a) Recalling the standard setup from Section~\ref{prel}, we have
\begin{equation}\label{pppdecomp}
X_t=\sum_{i\in I}X_t(W_i),\quad L^x=\sum_{i\in I}L^x(W_i),
\end{equation}
where $\{W_i:i\in I\}$ are the points of a Poisson point process, $\Xi$, with intensity $\N_{X_0}$.  Let \[F_i=\partial\{x:L^x(W_i)>0\}.\]  
Fix $\veps>0$, and define open sets $G_\veps=\{x:d(x,\textnormal{Supp}(X_0))<\veps\}$ and $U_\veps=\{x:d(x,\textnormal{Supp}(X_0))>\veps\}$. Recalling that $\hat W_t$ is the tip of the snake $W$ at time $t$ under $\N_{X_0}$, we set $S_\veps(W)=\inf\{t:\hat W_t\in G_\veps^c\}$ and $L^{U_\veps}(W)=(L^y(W),y\in U_\veps)$. We also use $L^{U_\veps}$ to denote the local time of $X$ restricted to $U_\veps$ under $\Pe_{X_0}$.  Then we have $\N_{X_0}$-a.e.,
\begin{equation*}
S_\veps(W)=\infty\Rightarrow X_{G_\veps}=0\Rightarrow L^{U_\veps}=0,
\end{equation*}
where the first implication follows from the definition of the exit measure $X_{G_\veps}$ (e.g. in Ch. IV of \cite{bib:leg99}) and the second from the special Markov property (Proposition~\ref{SpecialMProp}(b)).
It follows from the above and the decomposition in \eqref{pppdecomp} that 
\begin{equation}\label{LUdec}
L^{U_\veps}=\sum_{i\in I}L^{U_\veps}(W_i)1(S_\veps(W_i)<\infty),
\end{equation}
where the summation is componentwise. Note this represents the local time on $U_\veps$ as the integral of a Poisson point process $\N_{X_0}(\cdot,S_\veps<\infty)$ with finite total intensity $\N_{X_0}(S_\veps<\infty)$. 

We claim that 
\begin{equation}\label{Fdecomp}
F\cap U_\veps\subset \cup_{i\in I, S_\veps(W_i)<\infty}F_i\cap U_\veps,\quad \N_{X_0}-a.e.
\end{equation}
To see this let $x\in F\cap U_\veps$ and first note that $L^x=0$ implies that $L^x(W_i)=0$ for all $i$ by \eqref{pppdecomp}.  Also there is a sequence $x_n\in U_\veps$ converging to $x$ such that $L^{x_n}>0$. 
In view of \eqref{LUdec} and the fact that the summation there is a.e. finite, by taking a subsequence
we may assume there is an $i$ so that $S_\veps(W_i)<\infty$ and $L^{x_n}(W_i)>0$ for all $n$.  
This proves that $x\in F_i\cap U_\veps$, and the claim is established. 

It follows from \eqref{Fdecomp} that 
\begin{align*}
\Pe_{X_0}(\textnormal{dim}(F\cap U_\veps)>d+2-p)&\le \Pe_{X_0}(\exists\,i\in I\text{ so that }S_\veps(W_i)<\infty, \textnormal{dim}(F_i\cap U_\veps)>d+2-p)\\
&=1-\exp\Bigl(-\N_{X_0}(S_\veps<\infty,\textnormal{dim}(F\cap U_\veps)>d+2-p)\Bigr)\\
&\le 1-\exp\Bigl(-\N_{X_0}(\textnormal{dim}(F)>d+2-p)\Bigr)\\
&=0,
\end{align*}
the last by Theorem~\ref{dimN0} (which implies that $\N_x(\text{dim}(F)>d+2-p)=0$ for all $x$). We have shown that $\textnormal{dim}(F\cap U_\veps)\le d+2-p$ $\Pe_{X_0}$-a.s. and letting $\veps\downarrow 0$ completes the proof of (a).

\noindent(b) This is immediate from the upper bound in (a), the lower bound in Theorem~\ref{genlowerboundl}, and
the trivial inclusion $\textnormal{conv}(X_0)^c\subset \textnormal{Supp}(X_0)^c$.

\eop

\medskip

\noindent{\bf Proof of Proposition~\ref{noexit}.}
Let $B_{1+\ep}=B(0, 1+\ep)$ be an open ball centered at zero and radius $1+\ep$.
Clearly 
\begin{eqnarray}\label{step1}
\Pe_{X_0}(F\subset{\rm Supp}(X_0))&=& \Pe_{X_0}(F\subset \overline{B}_1)\\
\nonumber
&=& \lim_{\ep\downarrow 0} \Pe_{X_0}(F\subset \overline{B}_{1+\ep}).
\end{eqnarray}
 Now,
\begin{eqnarray}\label{step2}
 \Pe_{X_0}(F\subset \overline{B}_{1+\ep})&\geq&  \Pe_{X_0}(L^x =0, \;\forall x: |x|> 1+\ep)\\
\nonumber
&=&\Pe_{X_0}\left(\int_{\overline{B}_{1+\ep}^c} L^x\,dx =0\right),
\end{eqnarray}
where the last equality follows by the continuity of $x\mapsto L^x$.
By Theorem 1 of~\cite{bib:iscoe88} (see also the first equality on p. 205 there) we get that 
\begin{eqnarray*}
\Pe_{X_0}\left(\int_{\overline B^c_{1+\ep}} L^x\,dx =0\right)&=& e^{-\int_{ \overline{B}_1} v_{\ep}(x) X_0(dx)}
\end{eqnarray*}
where $v_{\ep}$ is the unique positive solution to 
\begin{eqnarray}\label{pde10}
\left\{ \begin{array}{rcl}
\Delta v_{\ep}&=&(v_\ep)^2,\;\; x\in B_{1+\ep},\\
\nonumber
 v_\ep(x)&\rightarrow& \infty,\;\;{\rm as}\; |x|\uparrow 1+\ep\,.
\end{array}\right.
\end{eqnarray}
By Proposition~9(ii) in Chapter V of ~\cite{bib:leg99} we get that there exists a constant $c_d$ such that 
\begin{eqnarray}\label{pde11}
 v_{\ep}(x)&\leq &c_d(1+\ep-|x|)^{-2},\; x\in B_{1+\ep}\,.
\end{eqnarray}
Thus we get 
\begin{eqnarray*}
 \Pe_{X_0}\Bigl(\int_{{\overline B_{1+\veps}}^c}L^x\,dx =0\Bigr)&\geq& 
e^{-\int_{\overline{B}_1} c_d(1+\ep-|x|)^{-2} X_0(dx)}
\\
&\geq&e^{-\int_{\overline{B}_1} c_d(1-|x|)^{-2} X_0(dx)}
\end{eqnarray*}
and we are done by our assumptions on $X_0$, \eqref{step1}, and \eqref{step2}.

\eop

\noindent{\bf Proof of Proposition~\ref{FLeb}.} We may, and shall, assume that 
\begin{equation}\label{denssupp}\{X_0>0\}\subset \text{Supp}(X_0).
\end{equation}
By \eqref{VinftygenX0} we have (recall $d=3$)
\begin{align*}
\Pe_{X_0}(L^x=0)&=\exp\Bigl(-2\int|x-x_0|^{-2}X_0(x_0)dx_0\Bigr)\\
&\ge\exp\Bigl(-2\Vert X_0\Vert_\infty\int_0^1 cr^{-2}r^2\,dr-2\int1(|x-x_0|\ge 1)X_0(x_0)\,dx_0\Bigr)\\
&\ge \exp\Bigl(-c_1\Vert X_0\Vert_\infty-c_2X_0(1)\Bigr)=p(X_0)>0.
\end{align*}
Now use Fubini to see that 
\[\Ee_{X_0}\Bigl(\int 1(X_0(x)>0, L^x=0)\,dx\Bigr)\ge p(X_0)|\{x:X_0(x)>0\}|>0,\]
where $|A|$ denotes the Lebesgue measure of $A$. Therefore
\begin{equation}\label{poslebs}
\Pe_{X_0}(|\{x:X_0(x)>0, L^x=0)\}|>0)>0.
\end{equation}

If $B$ is an open ball which intersects Supp$(X_0)$, then $t\to X_t(B)$ is $\Pe_{X_0}$-a.s. continuous on $[0,\infty)$ (see, e.g. Corollary~6 of \cite{bib:p91}).  As $X_0(B)>0$, this implies that $\int_BL^xdx=\int_0^\infty X_t(B)dt>0$ a.s., and therefore we have $B\cap \{x:L^x>0\}\neq \emptyset$ a.s. By considering a suitable countable collection of $B$'s we see that 
\begin{equation}\label{SinR}
\text{Supp}(X_0)\subset\overline{\{x:L^x>0\}}\quad \Pe_{X_0}-\text{a.s.}
\end{equation}
It follows that (recall \eqref{denssupp})
\begin{equation}\label{poslebcont}
\{X_0>0\}\cap\{x:L^x=0\}\subset\{X_0>0\}\cap \overline{\{x:L^x>0\}}\cap\{x:L^x=0\}=\{X_0>0\}\cap F.
\end{equation}
The left-hand side of the above has positive Lebesgue measure w.p.$>0$ by \eqref{poslebs} and so the same is true of $\{X_0>0\}\cap F$.\eop

\section{Proof of Proposition~\ref{L2upperbound}}\label{sec:9}
\setcounter{equation}{0}
\setcounter{theorem}{0}
We start this section with  a series of  lemmas that will help prove the proposition.
Recall that $d=2$ or $3$.

\begin{lemma}\label{2ptvl} Let $\vec x=(x_1,x_2)\in\IR^d\times\IR^d$, with $x_1\neq x_2$ and $\vec\lambda=(\lambda_1,\lambda_2)\in[0,\infty)^2\setminus\{(0,0)\}$.  There is a positive function $C^2$ function $V(x)=V^{\vec\lambda,\vec x}(x)$ on $\IR^d\setminus\{x_1,x_2\}$ such that 
\begin{equation}\label{v2pteq} \frac{\Delta V}{2}=\frac{V^2}{2}-\sum_{i=1}^2\lambda_i\delta_{x_i}\quad\text{on }\Re^d
\end{equation}
in the distributional sense, and $\Delta V=V^2$ on $\IR^d\setminus\{x_1,x_2\}$. Moreover for all $x\in\IR^d$, 
\begin{equation}\label{LTLT2}\Ee_{\delta_x}\Bigl(\exp\Bigl(-\sum_{i=1}^2\lambda_i L^{x_i}\Bigr)\Bigr)=\exp(-V^{\vec\lambda,\vec x}(x))=\exp\Bigl(-\int 1-\exp\Bigl(-\sum_{i=1}^2\lambda_iL^{x_i}(\nu)\Bigr)\,d\N_x(\nu)\Bigr).
\end{equation}
and for some 
$c_{\ref{2ptvl}}$, 
\begin{equation}\label{VL2bound}
V^{\vec{\lambda},\vec{x}}(x)\le c_{\ref{2ptvl}}\Bigl[\sum_{i=1}^2\lambda_i g_0(x-x_i)+1\Bigr]
\end{equation}
\end{lemma}

The proof is a minor modification of that of Lemma~\ref{LTVL}. The second equality in \eqref{LTLT2} is the analogue of \eqref{vcanonical}. 

Fix $\vec\lambda$ and $\vec x$ as in the above Lemma.  Below we will always assume $x\notin\{x_1,x_2\}$.  Monotone convergence shows we may differentiate the left-hand side of \eqref{LTLT2} with respect to $\lambda_i>0$ through the integral, and so conclude that for $i=1,2$, $V_i^{\vec\lambda,\vec x}(x)=\frac{\partial}{\partial\lambda_i}V^{\vec\lambda,\vec x}(x)$ exists and 
\begin{equation}\label{V1partial}
\Ee_{\delta_x}(L^{x_i}\exp(-\sum_{i=1}^2\lambda_i L^{x_i}))=e^{-V^{\vec\lambda,\vec x}(x)}V_i^{\vec\lambda,\vec x}(x)\quad\text{for }\lambda_i>0,\lambda_{3-i}\ge 0.
\end{equation}
Repeat the above to see that $V^{\vec\lambda,\vec x}(x)$ is $C^2$ in $\lambda_1,\lambda_2>0$, and if $U^{\vec\lambda,\vec x}(x)=\frac{\partial^2 }{\partial\lambda_1\partial \lambda_2}V^{\vec\lambda,\vec x}(x)$, then
\begin{equation}\label{V12partial}
\Ee_{\delta_x}\Bigl(L^{x_1}L^{x_2}\exp\Bigl(-\sum_{i=1}^2\lambda_i L^{x_i}\Bigr)\Bigr)=e^{-V^{\vec\lambda,\vec x}(x)}
\Bigl[V_1^{\vec\lambda,\vec x}(x)V_2^{\vec\lambda,\vec x}(x)-U^{\vec\lambda,\vec x}(x)\Bigr]\quad\text{for }\lambda_1,\lambda_2>0.
\end{equation}
Next, note that by \eqref{LTLT2} we have 
\begin{equation}\label{Vcanequation}
V^{\vec\lambda,\vec x}(x)=\int 1-\exp\Bigl(-\sum_{i=1}^2\lambda_iL^{x_i}(\nu)\Bigr)d\N_x(\nu).
\end{equation}

\begin{lemma}\label{Vmonotone}
\noindent(a) $V_i^{\vec\lambda,\vec x}(x)>0$ is strictly decreasing in $\vec\lambda\in\{(\lambda_1,\lambda_2):\lambda_i>0,\lambda_{3-i}\ge 0\}$, for $i=1,2$.

\noindent(b) $-U^{\vec\lambda,\vec x}(x)>0$ is strictly decreasing in $\vec\lambda\in(0,\infty)^2$.
\end{lemma}
\paragraph{Proof.} (a) Differentiate \eqref{Vcanequation} with respect to $\lambda_i>0$ to conclude
\begin{equation}\label{V1canmeas}
V_i^{\vec\lambda,\vec x}(x)=\int L^{x_i}(\nu)\exp\Bigl(-\sum_{i'=1}^2\lambda_{i'}L^{x_{i'}}(\nu)\Bigr)\N_x(d\nu)\quad\forall\lambda_i>0,
\end{equation}
where differentiation through the integral is easily justified by Monotone Convergence.  The finiteness of the integral on the right-hand side follows from the above. This shows the claimed montonicity of $V_i^{\vec\lambda,\vec x}(x)>0$ in $\vec\lambda$ because $\N_x(L^{x_i}>0)>0$. 

\noindent(b) If $\lambda_1,\lambda_2>0$ we can take $i=1$ and differentiate \eqref{V1canmeas}
with respect to $\lambda_2$, again using Monotone Convergence to differentiate through the integral, and so conclude
\[-U^{\vec\lambda,\vec x}(x)=\int L^{x_1}L^{x_2}\exp\Bigl(-\sum_{i=1}^2 \lambda_iL^{x_i}\Bigr) d\N_x>0.\]
The stated monotonicity in $\vec\lambda$ is now clear from the above and $\N_x(L^{x_1}L^{x_2}>0)>0$.
\eop

\begin{lemma}\label{Vibounds} There is a $C_{\ref{Vibounds}}>0$ so that:\\
\noindent(a) For all $\lambda_i>0$ and $\lambda_{3-i}\ge 0$,
\[V_i^{\vec\lambda,\vec x}(x)\le \frac{2}{\lambda_i}(V^{\lambda_i}(x_i-x)-V^{\lambda_i/2}(x_i-x))\le\frac{2}{\lambda_i}\Bigl(V^\infty(x_i-x)\wedge(C_{\ref{Vibounds}}\lambda_i^{-\alpha}|x_i-x|^{-p})\Bigr).\]
(b) For all $\lambda_1,\lambda_2>0$,
\begin{align*}-U^{\vec\lambda,\vec x}(x)&\le \frac{4}{\lambda_1\lambda_2}\min_{i=1,2}(V^{\lambda_i}(x_i-x)-V^{\lambda_i/2}(x_i-x))\\
&\le\frac{4}{\lambda_1\lambda_2}\Bigl(V^\infty(x_1-x)\wedge V^\infty(x_2-x)\wedge (C_{\ref{Vibounds}}\lambda_1^{-\alpha}|x_1-x|^{-p})\wedge (C_{\ref{Vibounds}}\lambda_2^{-\alpha}|x_2-x|^{-p}))\Bigr).
\end{align*}
\end{lemma}
\paragraph{Proof.} By Proposition~\ref{Vlambdarate} it suffices to establish the first inequalities in (a) and (b).\\
(a) By symmetry take $i=1$. The monotonicity in Lemma~\ref{Vmonotone}(a) and the Fundamental Theorem of Calculus imply
\begin{align*} V_1^{\vec\lambda,\vec x}(x)\le\frac{2}{\lambda_1}\int_{\lambda_1/2}^{\lambda_1} V_1^{(\lambda'_1,\lambda_2),\vec x}(x)d\lambda'_1
&\le \frac{2}{\lambda_1}\int_{\lambda_1/2}^{\lambda_1} V_1^{(\lambda'_1,0),\vec x}(x)d\lambda'_1\\
&=\frac{2}{\lambda_1}(V^{\lambda_1}(x_1-x)-V^{\lambda_1/2}(x_1-x)).
\end{align*}

(b) Argue as above using the monotonicity in Lemma~\ref{Vmonotone}(b) to see that for $\lambda_1,\lambda_2>0$,
\begin{align*}
-U^{\vec\lambda,\vec x}(x)&\le\frac{2}{\lambda_1}\int_{\lambda_1/2}^{\lambda_1}-\frac{\partial}{\partial\lambda_1'}V_2^{(\lambda'_1,\lambda_2),\vec x}(x)\, d\lambda'_1\\
&\le\frac{2}{\lambda_1}V_2^{(\lambda_1/2,\lambda_2),\vec x}(x)\\
&\le \frac{4}{\lambda_1\lambda_2}(V^{\lambda_2}(x_2-x)-V^{\lambda_2/2}(x_2-x)),
\end{align*}
where the last line follows from part (a) with $i=2$. The first inequality now follows by symmetry. 
\eop

In what follows we will always assume $0<\veps<\min\{|x_i-x|:i=1,2\}$. Set $T_\veps^i=\inf\{t\ge0:|B_t-x_i|\le\veps\}$ and $T_\veps=T_\veps^1\wedge T^2_\veps$, and let $\cF_t$ denote the right-continuous filtration generated by the Brownian motion $B$, which starts at $x$ under $P_x$. 

\begin{lemma}\label{v1martingale} Let $\lambda_1,\lambda_2>0$ and $\veps>0$.\\
(a) $V_1^{\vec\lambda,\vec x}(B(t\wedge T_\veps))-\int_0^{t\wedge T_\veps}V^{\vec\lambda,\vec x}(B(s))V_1^{\vec\lambda,\vec x}(B(s))\,ds$ is an $\cF_t$-martingale.\\
(b) For any $t\ge 0$, $V_1^{\vec\lambda,\vec x}(x)=E_x\Bigl(V_1^{\vec\lambda,\vec x}(B(t\wedge T_\veps))\exp\Bigl(-\int_0^{t\wedge T_\veps}V^{\vec\lambda,\vec x}(B(s))\,ds\Bigr)\Bigr)$.
\end{lemma}
\paragraph{Proof.} (a) 
By Lemma~\ref{2ptvl} for $\delta>0$,
\begin{equation}\label{approxnle}
\Bigl(\frac{\Delta V^{\vec\lambda+(\delta,0),\vec x}}{2}(x)-\frac{\Delta V^{\vec\lambda,\vec x}}{2}(x)\Bigr)\delta^{-1}=\Bigl(\Bigl(\frac{V^{\vec\lambda+(\delta,0),\vec x}(x)^2}{2}\Bigr)-\Bigl(\frac{ V^{\vec\lambda,\vec x}(x)^2}{2}\Bigr)\Bigr)\delta^{-1}\rightarrow V_1^{\vec\lambda,\vec x}(x)V^{\vec\lambda,\vec x}(x),
\end{equation}
as $\delta\to 0$.  Moreover the bounds on $V_1^{\vec\lambda,\vec x}$ in Lemma~\ref{Vibounds}(a) and 
on $V^{\vec\lambda,\vec x}$ in Lemma~\ref{2ptvl}  show that the above pointwise convergence is also uniformly bounded for $x$ satisfying $|x-x_i|>\veps$.  It\^o's lemma  shows that $V^{\vec\lambda,\vec x}(B(t\wedge T_\veps))-\int_0^{t\wedge T_\veps}\frac{\Delta V^{\vec\lambda,\vec x}(B(r))}{2}\,dr$ is an $\cF_t$-martingale. Therefore if $s<t$ and $\delta>0$, 
\begin{align}\label{premartingale}
E_x\Bigl(&\frac{V^{\vec\lambda+(\delta,0),\vec x}(B(t\wedge T_\veps))-V^{\vec\lambda,\vec x}(B(t\wedge T_\veps))}{\delta}\Bigr|\cF_s\Bigr)\\
\nonumber &=E_x\Bigl(\int_0^{t\wedge T_\veps}\frac{\Delta V^{\vec\lambda+(\delta,0),\vec x}(B(r))-\Delta V^{\vec\lambda,\vec x}(B(r))}{2\delta}\,dr\Bigr|\cF_s\Bigr).
\end{align}
The left-hand side of the above approaches $E_x(V_1^{\vec\lambda,\vec x}(B(t\wedge T_\veps))|\cF_s)$ as $\delta\to 0$ (the uniform boundedness of $V_1^{\vec\lambda,\vec x}(x)$ noted below allows us
to take the limit through the conditional expectation).  
The result now follows by letting $\delta\to 0$ in the above and applying \eqref{approxnle} and
the boundedness established above to take the limits through the conditional expectation and Lebesgue integral on the right-hand side.

\noindent(b) By (a), It\^o's lemma, and boundedness of $V_1^{\vec\lambda,\vec x}$ on $\{|x-x_1|\ge\veps\}$ (from Lemma~\ref{Vibounds}),\\ $V_1^{\vec\lambda,\vec x}(B(t\wedge T_\veps))\exp\Bigl(-\int_0^{t\wedge T_\veps}V^{\vec\lambda,\vec x}(B(s))ds\Bigr)$ is an $\cF_t$-martingale and so the result follows.  
\eop

\begin{lemma}\label{UFKrep}
For all $\lambda_1,\lambda_2>0, \veps>0$, 
\begin{align}
\nonumber
-U^{\vec\lambda,\vec x}(x)=&E_x\Bigl(\int_0^{T_\veps}\prod_{i=1}^2V_i^{\vec\lambda,\vec x}(B(t))\exp\Bigl(-\int_0^t V^{\vec\lambda,\vec x}(B(s))\,ds\Bigr)\,dt\Bigr)\\
\nonumber&\ \ + E_x\Bigl(\exp\Bigl(-\int_0^{T_\veps}V^{\vec\lambda, \vec x}(B(s))\,ds\Bigr)1(T_\veps<\infty)(-U^{\vec\lambda,\vec x}(B(T_\veps))\Bigr).
\end{align}
 
\end{lemma}
\paragraph{Proof.} Using the bounds in Lemma~\ref{Vibounds} one may easily differentiate the representation for $V_1^{\vec\lambda,\vec x}(x)$ in Lemma~\ref{v1martingale}(b) with respect to $\lambda_2>0$ through
the expectation and obtain
\begin{align}\label{Uldecomp1}
-U^{\vec\lambda,\vec x}(x)&=E_x\Bigl(V_1^{\vec\lambda,\vec x}(B(t\wedge T_\veps))\exp\Bigl(-\int_0^{t\wedge T_\veps} V^{\vec\lambda,\vec x}(B(s))\,ds\Bigr)\int_0^{t\wedge T_\veps}V_2^{\vec\lambda,\vec x}(B(s))\,ds\Bigr)\\
\nonumber&\quad -E_x\Bigl(U^{\vec\lambda,\vec x}(B(t\wedge T_\veps))\exp\Bigl(-\int_0^{t\wedge T_\veps} V^{\vec\lambda,\vec x}(B(s))\,ds\Bigr)\Bigr)\\
\nonumber&\equiv I_1(t)+I_2(t).
\end{align}
Use the Markov property, then Lemma~\ref{v1martingale}(b), and then Monotone Convergence to see that 
\begin{align}
\label{eq:17_10_1}
I_1(t)&=E_x\Bigl(\int_0^{t\wedge T_\veps} V_2^{\vec\lambda,\vec x}(B(s))\exp\Bigl(-\int_0^s V^{\vec\lambda,\vec x}(B(r))\,dr\Bigr)\\
\nonumber
&\phantom{=E_x\Bigl(\int_0^{t\wedge T_\veps} }\times E_{B(s)}\Bigl(V_1^{\vec\lambda,\vec x}(B((t-s)\wedge T_\veps)\exp\Bigl(-\int_0^{(t-s)\wedge T_\veps} V^{\vec\lambda,\vec x}(B(r))dr\Bigr)\,ds\Bigr)\\
\nonumber
&=E_x\Bigl(\int_0^{t\wedge T_\veps} V_2^{\vec\lambda,\vec x} V_1^{\vec\lambda, \vec x}(B(s))
\exp\Bigl(-\int_0^s V^{\vec\lambda,\vec x}(B(r))\,dr\Bigr)\,ds\Bigr)\\
\nonumber
&\rightarrow E_x\Bigl(\int_0^{T_\veps} V_2^{\vec\lambda,\vec x} V_1^{\vec\lambda, \vec x}(B(s))
\exp\Bigl(-\int_0^s V^{\vec\lambda,\vec x}(B(r))\,dr\Bigr)\,ds\Bigr),\;\;\text{as}\; t\to\infty.
\end{align}
Lemma~\ref{Vibounds}(b) shows that $-U^{\vec\lambda,\vec x}(x)$ is uniformly bounded on
$\{x:|x-x_i|\ge\veps,\ i=1,2\}$ and $\lim_{|x|\to\infty}-U^{\vec\lambda,\vec x}(x)=0$.  Therefore by Dominated Convergence, $$I_2(t)\to E_x\Bigl(\exp\Bigl(-\int_0^{T_\veps}V^{\vec\lambda, \vec x}(B(s))\,ds\Bigr)1(T_\veps<\infty)(-U^{\vec\lambda,\vec x}(B(T_\veps))\Bigr),$$ as $t\to\infty$, and this, 
together with \eqref{eq:17_10_1} and \eqref{Uldecomp1},
 completes the proof. 
\eop

\medskip
\noindent Now we are ready to turn to the
\paragraph{Proof of Proposition~\ref{L2upperbound}} It clearly suffices to obtain the result for $\lambda\ge \lambda_1(\veps_0)$ by adjusting $C_{\ref{L2upperbound}}$.  We will set $r_\lambda=\lambda_0\lambda^{-1/(4-d)}$, where $\lambda_0=\lambda_0(d)$ will be chosen large enough below, and we may assume $\lambda>\lambda_1(\veps_0)$ is large enough so that 
\begin{equation}\label{rlbound} r_\lambda<\veps_0.
\end{equation}
Recall that 
\[T^i_{r_\lambda}=\inf\{t:|B_t-x_i|\le r_\lambda\}\text{ and }T_{r_\lambda}=T^1_{r_\lambda}\wedge T^2_{r_\lambda}.\]
As we always assume $|x_i|\ge\veps_0$, we have $T_{r_\lambda}>0$ $P_0$-a.s. by \eqref{rlbound}.  We set $\vec{\lambda}=(\lambda,\lambda)$, $\vec x=(x_1,x_2)$, and $\Delta=|x_1-x_2|$. 

Apply \eqref{V12partial} and then Lemma~\ref{Vibounds}(a) to see that 
\begin{align}\nonumber\lambda^{2\alpha+2}\Ee_{\delta_0}\Bigl(L^{x_1}L^{x_2}\exp(-\lambda\sum_{i=1}^2L^{x_i})\Bigr) &=\lambda^{2\alpha+2}e^{-V^{\vec\lambda,\vec x}(0)}\Bigl(\prod_{i=1}^2 V_i^{\vec\lambda, \vec x}(0)-U^{\vec\lambda,\vec x}(0)\Bigr)\\
\nonumber&\le c\Bigl(\prod_{i=1}^2|x_i|^{-p}-\lambda^{2+2\alpha} U^{\vec\lambda,\vec x}(0)\Bigr)\\
\label{lbound1}&\le c\Bigl(\ve_0^{-2p}-\lambda^{2+2\alpha}U^{\vec\lambda,\vec x}(0)\Bigr).
\end{align}
To bound the last term, use Lemma~\ref{UFKrep} 
to arrive at
\begin{align}\label{USMP}-\lambda^{2+2\alpha}U^{\vec\lambda,\vec x}(0)=&\lambda^{2+2\alpha}E_0\Bigl(\int_0^{T_{r_\lambda}}\prod_{i=1}^2V_i^{\vec\lambda,\vec x}(B(t))\exp\Bigl(-\int_0^t V^{\vec\lambda,\vec x}(B(s))\,ds\Bigr)\,dt\Bigr)\\
\nonumber&\ \ +\lambda^{2+2\alpha} E_0\Bigl(\exp\Bigl(-\int_0^{T_{r_\lambda}}V^{\vec\lambda, \vec x}(B(s))\,ds\Bigr)1(T_{r_\lambda}<\infty)(-U^{\vec\lambda,\vec x}(B(T_{r_\lambda}))\Bigr)\\
\nonumber&\equiv K_1+K_2.
\end{align}

We first consider $K_2$. On $\{T_{r_\lambda}<\infty\}$ we may set $x_\lambda(\omega)=B(T_{r_\lambda})$ and choose $i(\omega)$ so that $|x_i-x_\lambda|\ge\Delta/2$. By definition of $T_{r_\lambda}$, $|x_i-x_\lambda|\ge r_\lambda$, and so $|x_i-x_\lambda|\ge\frac{1}{2}(\Delta\vee r_\lambda)$.  Lemma~\ref{Vibounds}(b) and the above imply
\[-\lambda^{2+2\alpha}U^{\vec\lambda, \vec x}(x_\lambda)\le\lambda^{2+2\alpha}\lambda^{-2-\alpha}4C_{\ref{Vibounds}}2^p(\Delta\vee r_\lambda)^{-p}=c\lambda^\alpha(\Delta\vee r_\lambda)^{-p}.\]
This shows that 
\begin{equation}\label{T2bound1}
K_2\le c\lambda^\alpha(\Delta\vee r_\lambda)^{-p}\sum_{i=1}^2E_0\Bigl(1(T^i_{r_\lambda}<\infty)\exp\Bigl(-\int_0^{T^i_{r_\lambda}} V^{\vec\lambda,\vec x}(B(s))\,ds\Bigr)\Bigr).
\end{equation}

By \eqref{LTLT} and \eqref{LTLT2}, and then Proposition~\ref{Vlambdarate}(a),  we have
\begin{equation}\label{Vveclb}
V^{\vec\lambda,\vec x}(B(s))\ge V^\lambda(B(s)-x_i)\ge V^\infty(B(s)-x_i)-C_{\ref{Vlambdarate}}|B(s)-x_i|^{-p}\lambda^{-\alpha}.
\end{equation}
Use the above in \eqref{T2bound1} and then use Brownian scaling to see that for $i=1,2$, (recall \break$\tau_r=\inf\{t:|B_t|\le r\}$)

\begin{align}\label{T2bound2}
E_0\Bigl(1&(T^i_{r_\lambda}<\infty)\exp\Bigl(-\int_0^{T^i_{r_\lambda}} V^{\vec\lambda,\vec x}(B(s))\,ds\Bigr)\Bigr)\\
\nonumber&\le E_{-x_i}\Bigl(1(\tau_{r_\lambda}<\infty) \exp\Bigl(\int_0^{\tau_{r_\lambda}}C_{\ref{Vlambdarate}}\lambda^{-\alpha}|B(s)|^{-p}ds\Bigr)\exp\Bigl(-\int_0^{\tau_{r_\lambda}}2(4-d)|B(s)|^{-2}\,ds\Bigr)\Bigr)\\
\nonumber&=E_{-x_i/r_\lambda}\Bigl(1(\tau_1<\infty)\exp\Bigl(\int_0^{\tau_1} C_{\ref{Vlambdarate}}\lambda^{-\alpha}r_\lambda^{2-p}|B(s)|^{-p}ds\Bigr)\exp\Bigl(-\int_0^{\tau_1}2(4-d)|B(s)|^{-2}ds\Bigr)\Bigr)\\
\nonumber&\le\liminf_{t\to\infty}E_{-x_i/r_\lambda}\Bigl(1(\tau_1<t)\exp\Bigl(\int_0^{\tau_1\wedge t}C_{\ref{Vlambdarate}}\lambda_0^{2-p}|B(s)|^{-p}\,ds\Bigr)\exp\Bigl(-\int_0^{\tau_1\wedge t}2(4-d)|B(s)|^{-2}\,ds\Bigr)\Bigr)\\
\nonumber&=\lim_{t\to\infty}E^{(2+2\nu)}_{|x_i|/r_\lambda}\Bigl(1(\tau_1<t)\exp\Bigl(\int_0^{\tau_1\wedge t}C_{\ref{Vlambdarate}}\lambda_0^{2-p}\rho_s^{-p}\,ds\Bigr)\rho_{t\wedge\tau_1}^{-\nu+\mu}\Bigr)(|x_i|/r_\lambda)^{\nu-\mu}\\
\nonumber&=E^{(2+2\nu)}_{|x_i|/r_\lambda}\Bigl(\exp\Bigl(\int_0^{\tau_1}C_{\ref{Vlambdarate}}\lambda_0^{2-p}\rho_s^{-p}ds\Bigr)\Bigr|\tau_1<\infty\Bigr)P_{|x_i|/r_\lambda}^{(2+2\nu)}(\tau_1<\infty)(|x_i|/r_\lambda)^{\nu-\mu},
\end{align}
where in the next to last line we have used Proposition~\ref{yorthm} with $\mu$ and $\nu$ as in \eqref{munudef}, so that $p=\mu+\nu$.  Now choose $\lambda_0$ large enough so that 
\[\sqrt{2C_{\ref{Vlambdarate}}\lambda_0^{2-p}}\le \sqrt{4(4-d)}\,(\le \nu).\]
This allows us to apply Lemma~\ref{expbound2} and conclude that the right-hand side of \eqref{T2bound2}
is at most
\begin{equation}\label{T2bound3}
cP^{(2+2\nu)}_{|x_i|/r_\lambda}(\tau_1<\infty)(|x_i|/r_\lambda)^{\nu-\mu}=c(|x_i|/r_\lambda)^{-2\nu+\nu-\mu}=c|x_i|^{-p}r_\lambda^p.\end{equation}
Now insert the above bound \eqref{T2bound3} into \eqref{T2bound1} to see that 
\begin{align}\nonumber K_2\le c\lambda^\alpha(\Delta\vee r_\lambda)^{-p}\veps_0^{-p}r_\lambda^p&\le c\veps_0^{-p}\lambda^\alpha r_\lambda^{-2}\Delta^{-(p-2)}r_\lambda^p\\
\label{T2bound4}&=c\ve_0^{-p}\lambda_0^{p-2}\Delta^{-(p-2)}.
\end{align}

In view of \eqref{lbound1}, \eqref{USMP} and \eqref{T2bound4}, it remains to establish
\begin{equation}\label{T1bound0} K_1\le C(\ve_0)\Bigl(1+\Delta^{2-p}\Bigr).
\end{equation}
Let $\Delta_i=x_{3-i}-x_i$, so that $|\Delta_i|=\Delta$, and let $T'_{r_\lambda}=\inf\{t:|B(t)|\le r_\lambda\text{ or }|B(t)-\Delta_i|\le r_\lambda\}$.  Lemma~\ref{Vibounds}(a) and then \eqref{Vveclb} give us
\begin{align}\label{T1bound1}
\nonumber K_1&\le cE_0\Bigl(\int_0^{T_{r_\lambda}} \prod_{i=1}^2|B(t)-x_i|^{-p}\exp\Bigl(-\int_0^t V^{\vec\lambda,\vec x}(B(s)ds\Bigr)dt\Bigr)\\
&\le c\sum_{i=1}^2 E_{-x_i}\Bigl(\int_0^{T'_{r_\lambda}} |B(t)|^{-p}|B(t)-\Delta_i|^{-p}1(|B(t)|\le |B(t)-\Delta_i|)\\
\nonumber&\phantom{\le c\sum_{i=1}^2 E_{-x_i}\Bigl(\int_0^{T'_{r_\lambda}}}\times\exp\Bigl(\int_0^tC_{\ref{Vlambdarate}}\lambda^{-\alpha}|B(s)|^{-p}ds\Bigr)\exp\Bigl(-\int_0^t V^\infty(B(s))ds\Bigr)dt\Bigr).
\end{align}
On $\{|B(t)|\le |B(t)-\Delta_i|\}$, a simple application of the triangle inequality shows that 
\[|B(t)-\Delta_i|\ge|\Delta_i|/2=\Delta/2,\]
 and so if $\mcF_t^{|B|}=\sigma(|B(s)|,s\le t)$, then
\[E(|B(t)-\Delta_i|^{-p}1(|B(t)|\le |B(t)-\Delta_i|)|\mcF_t^{|B|})\le 2^p(|B(t)|^{-p}\wedge\Delta^{-p}).\]
Use the above in \eqref{T1bound1} to obtain
\[K_1\le c\sum_{i=1}^2 E_{-x_i}\Bigl(\int_0^{\tau_{r_\lambda}}|B(t)|^{-p}(|B(t)|^{-p}\wedge \Delta^{-p})\exp \Bigl(\int_0^tC_{\ref{Vlambdarate}}\lambda^{-\alpha}|B(s)|^{-p}ds\Bigr)\exp\Bigl(-\int_0^t\frac{2(4-d)}{|B(s)|^2}ds\Bigr)\,dt\Bigr).
\]
Brownian scaling now gives
\begin{align*}
K_1&\le c\sum_{i=1}^2 E_{-x_i/r_\lambda}\Bigl(\int_0^{\tau_1}r_\lambda^{2-2p}|B(t)|^{-p}(|B(t)|^{-p}\wedge(\Delta/r_\lambda)^{-p})\exp\Bigl(\int_0^t C_{\ref{Vlambdarate}}\lambda^{-\alpha}r_\lambda^{2-p}|B(s)|^{-p}ds\Bigr)\\
&\phantom{\le \sum_{i=1}^2 E_{-x_i/r_\lambda}\Bigl(\int_0^{\tau_1}r_\lambda^{2-2p}|B(t)|^{-p}(|B(t)|^{-p}\wedge}\times\exp\Bigl(-\int_0^t\frac{2(4-d)}{|B(s)|^2}ds\Bigr)\,dt\Bigr)\\
&=cr_\lambda^{2-2p}\sum_{i=1}^2\int_0^\infty E_{-x_i/r_\lambda}\Bigl(1(t<\tau_1)|B(t\wedge\tau_1)|^{-p}(|B(t\wedge\tau_1)|\vee(\Delta/r_\lambda))^{-p}\\
&\phantom{\le \sum_{i=1}^2 E_{-x_i/r_\lambda}\Bigl(\int_0^{\tau_1}}\times\exp\Bigl(\int_0^{t\wedge\tau_1}C_{\ref{Vlambdarate}}\lambda_0^{2-p}|B(s)|^{-p}ds\Bigr)\exp\Bigl(-\int_0^{t\wedge\tau_1}\frac{2(4-d)}{|B(s)|^2}ds\Bigr)\Bigr)\,dt.
\end{align*}
Now we may use Proposition~\ref{yorthm} with $\mu,\nu$ as above (so that $p=\mu+\nu$) to see that if
\begin{equation}\label{deltadef}\delta=C_{\ref{Vlambdarate}}\lambda_0^{2-p},
\end{equation}
then
\begin{align}
\nonumber K_1&\le cr_\lambda^{2-2p}\sum_{i=1}^2E^{(2+2\nu)}_{|-x_i|/r_\lambda}\Bigl(\int_0^{\tau_1}\rho_t^{-p}(\rho_t\vee(\Delta/r_\lambda))^{-p}\exp\Bigl(\int_0^t\delta\rho_s^{-p}ds\Bigr)\\
\nonumber&\phantom{\le cr_\lambda^{2-2p}\sum_{i=1}^2R^{(2+2\nu)}_{-x_i/r_\lambda}\Bigl(\int_0^{\tau_1}}
\times\rho_t^{-\nu+\mu}(|x_i|/r_\lambda)^{\nu-\mu} dt\Bigr)\\
\label{T1bound2}&=cr_\lambda^{2-2p+\mu-\nu}\sum_{i=1}^2|x_i|^{\nu-\mu}E^{(2+2\nu)}_{|x_i|/r_\lambda}\Bigl(\int_0^{\tau_1}\rho_t^{-p-\nu+\mu}(\rho_t\vee(\Delta/r_\lambda))^{-p}\exp\Bigl(\int_0^t\delta\rho_s^{-p}ds\Bigr)\,dt\Bigr).
\end{align}

Fix $i\in\{1,2\}$, set $x_0=x_0(i)=|x_i|^2r_\lambda^{-2}>1$ (recall \eqref{rlbound}), and let $Y_t=\rho_t^2$.
 Then under $P_{x_0}=P_{|x_i|/r_\lambda}^{(2+2\nu)}$, $Y$ satisfies 
\[Y_t=x_0+\int_0^t 2\sqrt{Y_s}dW_s+(2+2\nu)t,\]
where $W$ is a Brownian motion.  We let $\tau^Y_r$ denote the hitting time of $r$ by $Y$ and set $q=(p-2)/2$. It\^o's Lemma implies that if $M_t=-2q\int_0^{t\wedge \tau^Y_1}Y_s^{-q-\frac{1}{2}}dW_s$, then 
\begin{align*}
Y_{t\wedge\tau^Y_1}^{-q}=x_0^{-q}+M_t+2q(q-\nu)\int_0^{t\wedge \tau^Y_1}Y_s^{-q-1}ds.
\end{align*}
This implies that 
\begin{align*}
M_t-\frac{1}{2}\langle M\rangle_t&=Y_{t\wedge\tau_1^Y}^{-q}-x_0^{-q}+\int_0^{t\wedge \tau^Y_1}2q(\nu-q)Y_s^{-q-1}-2q^2Y_s^{-2q-1}ds\\
&\ge Y_{t\wedge\tau_1^Y}^{-q}-x_0^{-q}+\int_0^{t\wedge \tau^Y_1}2q(\nu-2q)Y_s^{-q-1}ds.
\end{align*}
The constant inside the integral is $(p-2)(2-\mu)>0$ and so we may choose $\lambda_0(d)$ sufficiently large so that $\delta\le (p-2)(2-\mu)$.  If $\Ecal(M)_t=\exp(M_t-\langle M\rangle_t/2)$ is the stochastic exponential of $M$ the above shows that 
\begin{align*}\Ecal(M)_t\ge \exp(-x_0^{-q})\exp\Bigl(\int_0^{t\wedge\tau^Y_1}\delta Y_s^{-q-1}ds\Bigr),
\end{align*}
and so 
\begin{equation}\label{expbnd1}
\exp\Bigl\{\int_0^{t\wedge\tau_1^Y}\delta Y_s^{-p/2}ds\Bigr\}\le e\Ecal(M)_t.
\end{equation}
Noting that $-p-\nu+\mu=-2\nu$ and recalling \eqref{T1bound2}, from the above bound we arrive at
\begin{equation}\label{T1bound2.5}
K_1\le cr_\lambda^{2-2p+\mu-\nu}\sum_{i=1}^2 |x_i|^{\nu-\mu}e\int_0^\infty E_{x_0(i)}(\Ecal(M)_t1(t<\tau^Y_1)Y_t^{-\nu}(Y_t\vee(\Delta/r_\lambda)^2)^{-p/2})\,dt.
\end{equation}
Now $M$ is a martingale with $\langle M\rangle_t=\int_0^{t\wedge\tau^Y_1}4q^2 Y_s^{-2q-1}ds\le 4q^2 t$ and so by Girsanov's 
theorem (see, e.g. Chapter IV.4 of \cite{bib:iw79}) there is a unique probability $Q_{x_0(i)}$ on $C([0,\infty),\Re_+)$ so that 
if we also let $Y$ denote the coordinate variables on this space with its generated right continuous filtration $(\cF_t)$, then for any $t\ge 0$,
$dQ_{x_0(i)}|_{\cF_t}=\Ecal(M)_t dP|_{\cF_t}$, and under $Q_{x_0(i)}$, $Y$ is the unique solution of 
\begin{equation}\label{Qsde}
Y_t=x_0(i)+2\int_0^{t\wedge \tau^Y_1}\sqrt{Y_s}\,dW_s+(2+2\nu)(t\wedge\tau^Y_1)-2(p-2)\int_0^{t\wedge\tau^Y_1} Y_s^{-q}ds,
\end{equation}
(so $Y$ is stopped when it hits $1$).  Therefore if we use $Q_{x_0}$ to also denote expectation with respect to $Q_{x_0}$, then we have
\begin{align}\label{T1bound3}
K_1&\le cr_\lambda^{2-2p+\mu-\nu}\sum_{i=1}^2 |x_i|^{\nu-\mu}eQ_{x_0(i)}\Bigl(\int_0^{\tau_1^Y}Y_t^{-\nu}(Y_t\vee(\Delta/r_\lambda)^2)^{-p/2}\,dt\Bigr)\\
\nonumber&\equiv cr_\lambda^{2-2p+\mu-\nu}\sum_{i=1}^2 |x_i|^{\nu-\mu}eJ_i.
\end{align}
We interrupt the proof of the proposition for another auxiliary result. 
\begin{lemma}\label{Qdiff} Assume $Y$ and $Q_{x_0}$ are as in \eqref{Qsde} and $1\le a\le x_0$.\\
\noindent(a) $Q_{x_0}(\tau_a<\infty)\le e^2(a/x_0)^\nu$.

\noindent(b) For $\gamma>1$,
\[Q_{x_0}\Bigl(\int_0^{\tau^Y_a}Y_t^{-\gamma}\,dt\Bigr)\le\begin{cases}\frac{e^2}{2(\gamma-1-\nu)(\gamma-1)}\cdot\frac{a^{1+\nu-\gamma}}{x_0^\nu}&\text{ if }\gamma>1+\nu\\
\frac{x_0^{1-\gamma}}{2(1+\nu-\gamma-(p-2)a^{(p-2)/2})(\gamma-1)}&\text{ if }\gamma+(p-2)a^{(p-2)/2}<1+\nu.
\end{cases}\]
\end{lemma}
\paragraph{Proof of Lemma.}(a) It is easy to check that $s(x)=-\int_x^\infty y^{-1-\nu}\exp(-2y^{1-(p/2)})\,dy$ is a scale function for $Y$ under $Q_{x_0}$ and so (recall that $x_0\ge a$),
\begin{align*}
Q_{x_0}(\tau^Y_a<\infty)=\frac{s(\infty)-s(x_0)}{s(\infty)-s(a)}&=\frac{\int_{x_0}^\infty y^{-1-\nu}\exp(-2y^{1-(p/2)})\,dy}{\int_a^\infty y^{-1-\nu}\exp(-2y^{1-(p/2)})\,dy}\\
&\le\frac{\int^\infty_{x_0}y^{-1-\nu}\,dy}{e^{-2}\int_a^\infty y^{-1-\nu}\,dy}=e^2(x_0/a)^{-\nu}.
\end{align*}
\noindent(b) Let $g(x)=-x^{1-\gamma}/(\gamma-1)$, so that $g'(x)=x^{-\gamma}$. An application of It\^o's Lemma gives
\begin{equation}\label{ito}(2\gamma-2-2\nu)\int_0^{t\wedge\tau^Y_a}Y_s^{-\gamma}\,ds=-g(Y_{t\wedge\tau^Y_a})+g(x_0)+2\int_0^{t\wedge\tau_a^Y}Y_s^{-\gamma+(1/2)}dW_s-
2(p-2)\int_0^{t\wedge\tau^Y_a}Y_s^{-\gamma-q}\,ds.
\end{equation}
{\bf Case 1.} $\gamma>1+\nu$.\\
Take means in \eqref{ito}, drop the second and last terms (both negative) on the right-hand side and then let $t\to\infty$ to conclude that 
\begin{equation}\label{Jbound2}
Q_{x_0}\Bigl(\int_0^{\tau_a^Y}Y_s^{-\gamma}\,ds\Bigr)\le \limsup_{t\to\infty}-(2(\gamma-1-\nu))^{-1}Q_{x_0}(g(Y_{t\wedge \tau_a^Y})).
\end{equation}
The drift of $Y$ in \eqref{Qsde} is bounded below by $(2+2\nu-2(p-2))1(t\le \tau_1^Y)\ge 5\times1(t\le \tau^Y_1)$, and so by a standard comparison with the square of a $5$-dimensional Bessel process we see that a.s. on $\{\tau^Y_a=\infty\}$, $Y_t\to\infty$ as $t\to\infty$. Therefore \eqref{Jbound2} implies that 
\begin{equation}\label{Jbound3}
Q_{x_0}\Bigl(\int_0^{\tau^Y_a}Y_s^{-\gamma}\,ds\Bigr)\le \frac{a^{1-\gamma}}{2(\gamma-1-\nu)(\gamma-1)}Q_{x_0}(\tau^Y_a<\infty).\end{equation}
An application of (a) now completes the proof of (b) in this case.

\medskip

\noindent{\bf Case 2.} $\gamma+(p-2)a^{(p-2)/2}<1+\nu$.\\
In this case \eqref{ito} implies
\begin{align*}
-g(x_0)-\int_0^{t\wedge\tau^Y_a}2Y_s^{-\gamma+(1/2)}dW_s&\ge2(1+\nu-\gamma)\int_0^{t\wedge\tau^Y_a} Y_s^{-\gamma}\,ds-2(p-2)a^{-q}\int_0^{t\wedge\tau_a^Y} Y_s^{-\gamma}\,ds\\
&=2(1+\nu-\gamma-(p-2)a^{-(p-2)/2})\int_0^{t\wedge\tau_a^Y}Y_s^{-\gamma}ds.
\end{align*}
Take means and let $t\to\infty$ to conclude that 
\[Q_{x_0}\Bigl(\int_0^{\tau^Y_a} Y_s^{-\gamma}ds\Bigr)\le\frac{x_0^{1-\gamma}}{2(1+\nu-\gamma-(p-2)a^{-(p-2)/2})(\gamma-1)}.\text{\eop}\]

Returning now to the proof of Proposition~\ref{L2upperbound}, fix a value of $i\in\{1,2\}$ and consider:\\
{\bf Case 1.} $|x_i|>\Delta>r_\lambda$.\\
Then $x_0(i)>b\equiv(\Delta/r_\lambda)^2>1$ and so by the strong Markov property,
\begin{equation}\label{Jbound1}
J_i\le Q_{x_0(i)}\Bigl(\int_0^{\tau^Y_b} Y_t^{-\nu-(p/2)}\,\Bigr)+Q_{x_0(i)}(\tau_b<\infty)Q_b\Bigl(\int_0^{\tau^Y_1}Y_t^{-\nu}\,dt\Bigr)(r_\lambda/\Delta)^p\equiv I_1+I_2.
\end{equation}
Consider first $I_1$. Apply Lemma~\ref{Qdiff}(b) with $\gamma=\nu+(p/2)>\nu+1$ and $a=b\in(1,x_0(i))$ to see that 
\begin{align}\label{Jbound4}I_1&\le\frac{e^2}{(p-2)(\gamma-1)}\frac{b^{1-(p/2)}}{x_0(i)^\nu}=cr_\lambda^{-2+2p+\nu-\mu}|x_i|^{-2\nu}\Delta^{2-p},
\end{align}
where to get the power on $r_\lambda$ we used the identity $p=\nu+\mu$.

Turning next to $I_2$, note that $p<3$ implies that $\nu+p-2<\nu+1$ and so we may apply the second inequality in Lemma~\ref{Qdiff}(b) with $\gamma=\nu$, $a=1$, and $x_0=b>1$, to conclude that
\[Q_b\Bigl(\int_0^{\tau^Y_1}Y_t^{-\nu}\,dt\Bigr)\le \frac{b^{1-\nu}}{2(1-(p-2))(\nu-1)}=\frac{b^{1-\nu}}{(6-2p)(\nu-1)}.\]
Next use Lemma~\ref{Qdiff}(a) with $a=b$ and the above to see that 
\begin{align}
\nonumber I_2&\le e^2(b/x_0(i))^\nu\frac{b^{1-\nu}}{(6-2p)(\nu-1)}(r_\lambda/\Delta)^p\\
\label{Jbound5}&=cr_\lambda^{-2+2p+\nu-\mu}|x_i|^{-2\nu}\Delta^{2-p}.
\end{align}
So using \eqref{Jbound5} and \eqref{Jbound4} to bound $J_i$ in \eqref{Jbound1}, and recalling \eqref{T1bound3}, we arrive at 
\begin{equation}\label{T1bound4}
K_1\le c\sum_{i=1}^2|x_i|^{-\mu-\nu}\Delta^{2-p}\le c\veps_0^{-p}\Delta^{2-p}.
\end{equation}
\noindent{\bf Case 2.} $\Delta\le r_\lambda$.\\
In this case we apply Lemma~\ref{Qdiff}(b) with $\gamma=\nu+(p/2)>\nu+1$ and $a=1$ to see that
\begin{align*}
J_i=Q_{x_0(i)}\Bigl(\int_0^{\tau_1^Y} Y_t^{-\nu-(p/2)}dt\Bigr)\le\frac{e^2}{(p-2)(\nu+(p/2)-1)}x_0(i)^{-\nu}=cr_\lambda^{2\nu}|x_i|^{-2\nu}.
\end{align*}
So in this case, by \eqref{T1bound3} we again conclude that 
\[K_1\le cr_\lambda^{2-2p+\mu+\nu}\sum_{i=1}^2|x_i|^{-\nu-\mu}\le c\veps_0^{-p}r_\lambda^{2-p}\le c\veps_0^{-p}\Delta^{2-p},\]
where in the last inequality we use our assumption that $\Delta\le r_\lambda$.\\

\noindent{\bf Case 3.} $\Delta\ge(r_\lambda\vee|x_i|)(\ge \veps_0)$.\\
Apply Lemma~\ref{Qdiff}(b) with $a=1$ and $\gamma=\nu$, for which $\gamma+(p-2)<1+\nu$, to see that 
\begin{align*}
J_i&\le(r_\lambda/\Delta)^{p}Q_{x_0(i)}\Bigl(\int_0^{\tau_1^Y}Y_t^{-\nu}\,dt\Bigr)\\
&\le  (r_\lambda/\Delta)^{p}(2(3-p)(\nu-1))^{-1}x_0(i)^{1-\nu}=cr_\lambda^{p+2\nu-2}\Delta^{-p}|x_i|^{2(1-\nu)}.
\end{align*}
So again in this case we have from \eqref{T1bound3} that
\[K_1\le cr_\lambda^{2-2p+\mu-\nu+p+2\nu-2}\Delta^{-p}\sum_{i=1}^2|x_i|^{2-p}\le c\Delta^{-p}\sum_{i=1}^2|x_i|^{2-p}\le c\veps_0^{2-2p}.\]

We have established \eqref{T1bound0} in each possible case and so the proof is complete.\eop

\end{document}